\documentclass[12pt,openany]{article}
\usepackage{epsfig}

\newenvironment{proof}{{\it Proof:\/}}{$\Box$\vskip 0.08in}
\usepackage{amssymb}

 \newtheorem{theorem}{Theorem}[section]
 \newtheorem{lemma}[theorem]{Lemma}
 \newtheorem{corollary}[theorem]{Corollary}
 \newtheorem{remark}[theorem]{Remark}
\newtheorem{exercise}[theorem]{Exercise}
\newtheorem{definition}[theorem]{Definition}
\newcommand{\mod}{{\mbox{ mod }}}

\newtheorem{problem}[theorem]{Problem}
\newtheorem{example}[theorem]{Example}

\newtheorem{proposition}[theorem]{Proposition}
\newtheorem{construction}[theorem]{Construction}

\newcommand{\sgn}{{\mbox{ sgn\ }}}
\newcommand{\nul}{{\mbox{ nul }}}
\newcommand{\lk}{{\mbox{ lk }}}

\newcommand{\kwad}{\#}

\newcommand{\pct}[1]{}

\begin{document}

\thispagestyle{empty}
\
\vspace{0.4in}
 \begin{center}
{\LARGE \bf From Goeritz matrices to quasi-alternating links}
\end{center}

\vspace*{0.2in}

\centerline{by}
 \begin{center}
                      {\LARGE \bf J\'ozef H.~Przytycki}
\end{center}

\vspace*{0.2in}

\ \\
{\LARGE \bf Introduction}\\
\ \\

Knot Theory is currently a very broad field. Even a long survey can only cover a 
narrow area. Here we concentrate on the path from Goeritz matrices to quasi-alternating links.
On the way, we often stray  from the main road and tell related stories, especially if they allow as 
to place the main topic in a historical context. For example, we mention that the Goeritz 
matrix was preceded by the Kirchhoff matrix of an electrical network. The network complexity
 extracted from the matrix corresponds to the determinant of a link. 
We assume basic knowledge of knot theory and graph theory, however, we offer 
a short introduction under the guise of a historical perspective.

\section{Short historical introduction} 

Combinatorics, graph theory, and knot theory have their common roots in 
Gottfried Wilhelm Leibniz' (1646-1716) ideas of 
{\it Ars Combinatoria}, and {\it Geometria Situs}. In Ars Combinatora, Leibniz was influenced 
by Ramon Llull (1232 -- 1315) and his combinatorial machines (Figure 1.1; \cite {Bon}). 
\\ \ \\
\begin{center}
\psfig{figure=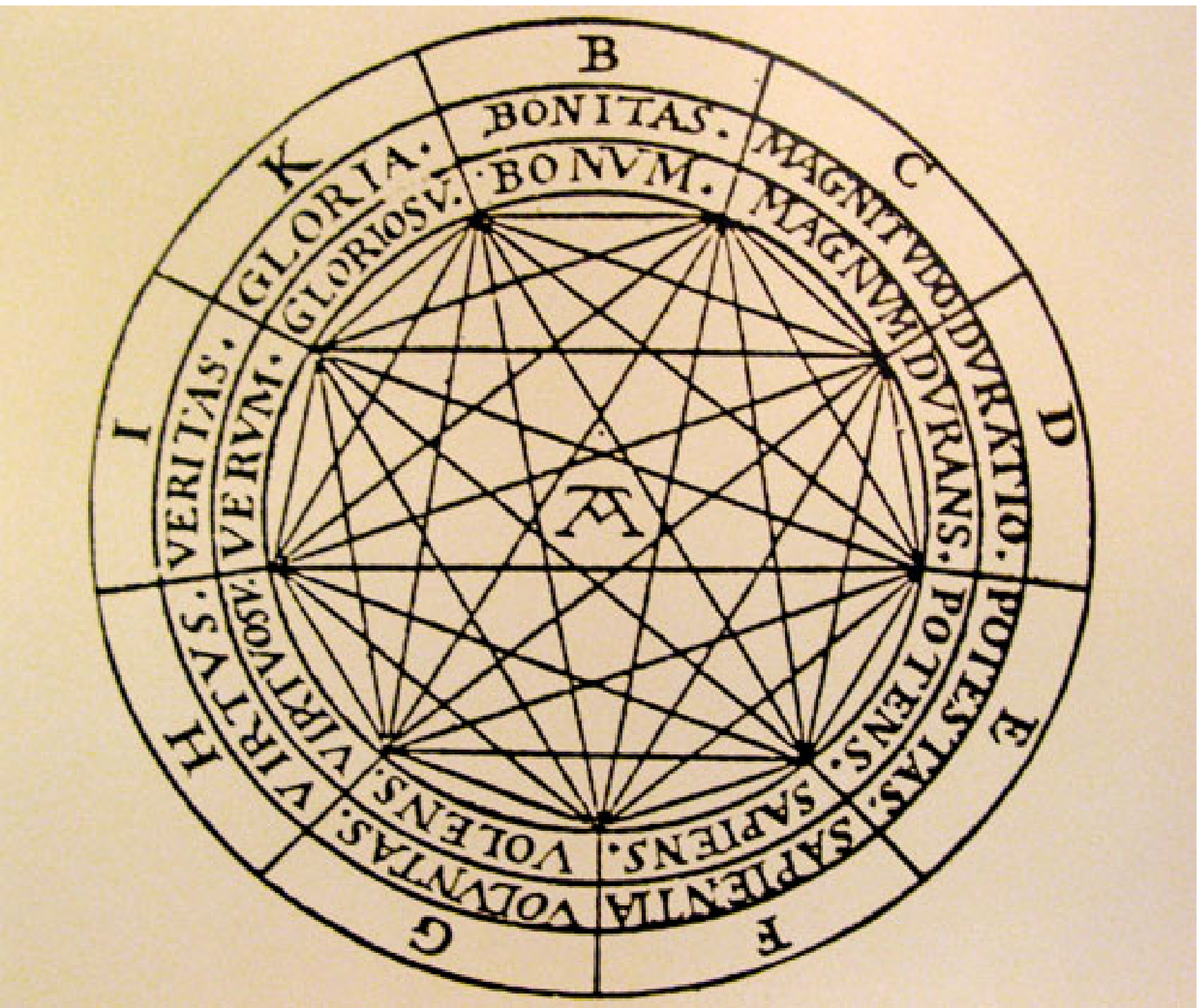,height=5.3cm}
\end{center}
\centerline{Fig. 1.1; Combinatorial machine of Ramon Llull from his Ars Generalis Ultima}
\ \\ \ \\
Geometria (or Analysis) Situs seems to be an 
invention of Leibniz. I am not aware of any  Ancient or Renaissance 
influence  (compare however \cite{P-21}. The first convincing example of geometria situs  
was proposed by Heinrich Kuhn in a letter
written in 1735 to Leonard Euler (1707-1783). Kuhn (1690-1769) was a Danzig (Gda\'nsk) mathematician 
 born in K\"onigsberg, studied at the Pedagogicum there,
and in 1733 settled in Danzig as a mathematics professor at the Academic Gymnasium (he was also 
a co-founder of the Nature Society) \cite{Janus}. Kuhn communicated to Euler the puzzle of bridges 
of K\"onigsberg, suggesting that this may be an example of geometria situs. Kuhn was communicating, 
in fact, through his friend Carl Leonhard Gottlieb Ehler (1685-1753), correspondent of Leibniz and  
future mayor of Danzig. The first letter by Ehler did not survive  but in the letter of March 9, 1736 
he writes: ``You would render to me and our friend K\"ohn a most valuable service, putting us greatly in your debt,
most learned Sir, if you would send us the solution , which you know well, to the
problem of the seven K\"onigsberg bridges, together with a proof.  It would prove to be an outstanding example
of {\bf Calculi Situs}, worthy of your great genius.  I have added a sketch of the said bridges ..."
In the reply of April 3, 1736 Euler writes ``... Thus you see,
most noble Sir, how this type of solution bears little relationship to
mathematics, and I do not understand why you expect a mathematician to produce it,
rather than anyone else, for the solution is based on reason alone, and its discovery does
not depend on any mathematical principle. Because of this,
I do not know why even questions which bear so little relationship to mathematics
are solved more quickly by mathematicians than by others. In the meantime, most noble Sir, you have assigned this
question to the {\bf geometry of position}, but I am ignorant as to what this new discipline
involves, and as to which types of problem Leibniz and Wolff expected to see expressed
in this way ... " \cite{H-W}. However when composing his famous paper on bridges of K\"onigsberg, Euler 
already agrees with Kuhn suggestion. The geometry of position figures even in the title of the 
paper {\it Solutio problematis ad geometriam situs pertinentis}.\footnote{In the paper, Euler writes: 
``The branch of geometry that deals with magnitudes has been zealously studied
throughout the past,  but there is another branch that has been almost
unknown up to now; Leibniz spoke of it first,  calling it the
``geometry of position" (geometria situs).
This branch of geometry deals with relations dependent
on position; it does not take magnitudes into considerations, nor does it
involve calculation with quantities. But as yet no satisfactory definition
has been given of the problems that belong to this geometry of position
or of the method to be used in solving them. Hence, when a problem was recently mentioned, 
which seemed geometrical but was so constructed that it did not require the measurement of 
distances, nor did calculation help at all, I had no doubt that it was concerned with the 
geometry of position--especially as its solution involved only position, and no calculation 
was of any use. I have therefore decided to give here the method which I have found for 
solving this kind of problem, as an example of the geometry of position.\
2.\ The problem, which I am told is widely known, is as follows: in K\"onigsberg in Prussia,
there is...  "\cite{Eu,B-L-W}.} 

The first paper mentioning knots from the mathematical point of view is that of 
Alexandre-Theophile Vandermonde (1735-1796) {\it Remarques sur les probl\`emes de situation} \cite{Va}.
Carl Friedrich Gauss (1777-1855) had interest in Knot Theory whole his life, starting from 1794 drawings of knots, 
the drawing of a braid with complex coordinates (c. 1820), several drawing of knots with 
``Gaussian codes", and Gauss' linking number of 1833. He did not publish anything however; this 
was left to his student Johann Benedict Listing (1808-1882) who in 1847 published his monograph 
(Vorstudien zur Topologie, \cite{Lis}). The monograph is mostly devoted to knots, graphs and combinatorics.


In the XIX century Knot Theory was an experimental science.
Topology (or geometria situs) had not developed enough to
offer tools allowing precise definitions and proofs\footnote{Listing writes in \cite{Lis}:
{\it In order to reach the level of exact science, topology will
have to translate facts of spatial contemplation into easier notion
which, using corresponding symbols analogous to mathematical ones, we
will be able to do corresponding operations following some simple rules.}} (here
Gaussian linking number is an exception). Furthermore, in
the second half of that century Knot Theory was developed
mostly by physicists (William Thomson (Lord Kelvin)(1824-1907), James Clerk Maxwell (1831-1879), 
Peter Guthrie Tait(1831-1901)) 
and one can argue that the high
level of precision was not appreciated\footnote{This may be
a controversial statement. The precision of Maxwell was different
than that of Tait and both were physicists.}.
We outline the global history of the Knot Theory in \cite{P-21} and in the second chapter 
of my book on Knot Theory \cite{P-Book}.
 In the next subsection we deal with the mathematics
developed in order to understand precisely the phenomenon of
knotting.

\subsection{Precision comes to Knot Theory}
Throughout the XIX century knots were understood as closed curves in a space
up to a natural deformation, which was described as a movement in
space without cutting and pasting. This understanding allowed scientists
(Tait, Thomas Penyngton Kirkman, Charles Newton Little, Mary Gertrude Haseman) 
to build tables of knots but didn't lead to precise methods
 allowing one to distinguish
knots which could not be practically deformed from one to another.
In a letter to O. Veblen, written in 1919, young J. Alexander expressed
his disappointment\footnote{We should remember that it was written by
a young revolutionary mathematician forgetting that he is ``standing on
the shoulders of giants." \cite{New}. In fact the invariant Alexander outlined in the letter 
is closely related to  Kirchhoff matrix, and extracted numerical invariant is equivalent to 
complexity of a signed graph corresponding to the link via Tait translation; 
see Subsection \ref{Subsection 1.4}.}:
{\it ``When looking over Tait {\it On Knots} among other things,
He really doesn't get very far. He merely writes down all the plane projections
of knots with a limited number of crossings, tries out a few transformations
that he happen to think of and assumes without proof that if he is unable
to reduce one knot to another with a reasonable number of tries, the two
are distinct. His invariant, the generalization of the Gaussian invariant ...
for links is an invariant merely of the particular projection of the knot
that you are dealing with, - the very thing I kept running up against
in trying to get an integral that would apply. The same is true of his
`Beknottednes'."}

In the famous Mathematical Encyclopedia Max Dehn and
Poul Heegaard outlined a systematic approach to topology,
in particular they precisely formulated the subject of
the Knot Theory \cite{D-H}, in 1907. To bypass the notion
of deformation of a curve in a space (then not yet well defined)
they introduced lattice knots and the precise definition of
their (lattice) equivalence.
Later Reidemeister and Alexander considered more general
polygonal knots in a space with equivalent knots related by
a sequence of $\Delta$-moves; they also explained $\Delta$-moves by elementary moves 
on link diagrams -- Reidemeister moves (see Subsection 1.6).
 The definition of Dehn and Heegaard
was long ignored and only recently lattice knots are again studied.
It is a folklore result, probably never written down in detail\footnote{It is however long 
routine exercise}, that the
two concepts, {\it lattice knots} and {\it polygonal knots}, are
equivalent.

\subsection{Lattice knots and Polygonal knots}

In this part we discuss two early XX century
definitions of knots and their
equivalence, by Dehn-Heegaard and by Reidemeister. In the XIX century
knots were treated from the intuitive point of view and 
 was P.~Heegaard in his 1898 thesis who came close to a
formal proof that there are nontrivial knots.

Dehn and Heegaard gave the following definition of a knot (or curve in
their terminology) and of equivalence of knots (which they call isotopy
of curves)\footnote{Translation from German due to Chris Lamm.}.

\begin{definition}[\cite{D-H}]\label{I.1.1}\ \\
A curve is a simple closed polygon on a cubical lattice.
It has coordinates $x_i,y_i,z_i$ and an {\it isotopy} of these curves
is given by:
\begin{enumerate}
\item[(i)]
Multiplication of every coordinate by a natural number,
\item[(ii)]
Insertion of an elementary square, when it does not interfere
with the rest of the polygon.
\item[(iii)]
Deletion of the elementary square.
\end{enumerate}
\end{definition}

Elementary moves of Dehn and Heegaard can be  summarized/explained as
follows:
\begin{enumerate}
\item[$(DH_0)$] Rescaling. We show in \cite{P-Book} that this move is a consequence of
other Dehn-Heegaard moves.
\item[$(DH_1)$] If a unit square intersects the lattice knot in exactly two
neighboring edges then we replace this edges by two other edges of the square,
as illustrated in Fig. 1.2 $(DH_1)$.
\item[$(DH_2)$] If a unit square intersects the lattice knot in exactly
one edge then we replace this edges by three other edges of the square,
as illustrated in Fig. 1.2 $(DH_2)$.
\end{enumerate}

\begin{center}
\psfig{figure=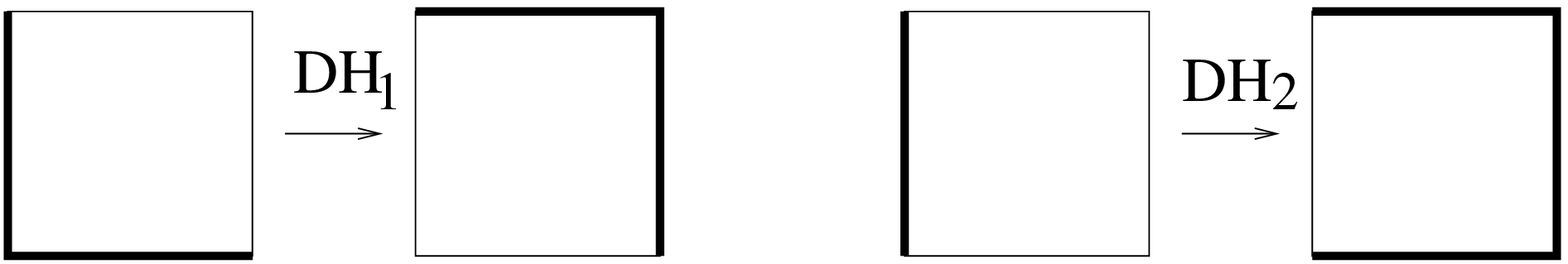,height=2.3cm}
\end{center}
\centerline{Fig. 1.2; Lattice moves $DH_1$ and $DH_2$}

In this language,  lattice knots (or links) and lattice isotopy
 are defined as follows.
\begin{definition}\label{1.2}
A lattice knot is a simple closed polygon on a cubical lattice.
Its vertices have integer coordinates $x_i,y_i,z_i$ and edges, of length one,
are parallel to one of the coordinate axis. We say that two lattice knots
are lattice isotopic if they are related by a finite sequence of
elementary lattice (``square") moves as illustrated in Fig. 1.2 (we allow
$DH_1$-move, $DH_2$-move and its inverse $DH_2^{-1}$-move).
These are moves (ii) and (iii) of Dehn and Heegaard.
\end{definition}

Below we give a few examples of lattice knots.

They can be easily coded as (cyclic) words over the
alphabet $\{x^{\pm 1},y^{\pm 1},z^{\pm 1}\}$. For example the trivial
knot can be represented by $xyx^{-1}y^{-1}$, the trefoil knot by
$x^2z^3y^2x^{-1}z^{-2}y^{-3}zx^2y^2x^{-3}y^{-1}z^{-2}$, and
the figure-eight knot \\
by $y^2z^2xy^{-3}x^2y^2z^{-1}x^{-4}y^{-2}x^3yz^2x^{-2}z^{-3}$;
see Figure 1.3.

\bigskip
\begin{center}
\psfig{figure=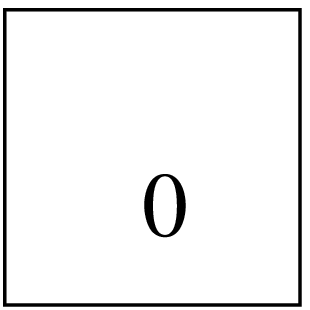,height=1.3cm}\ \ \ \ \ \ \ \ \ \ \ \
 \ \ \
\psfig{figure=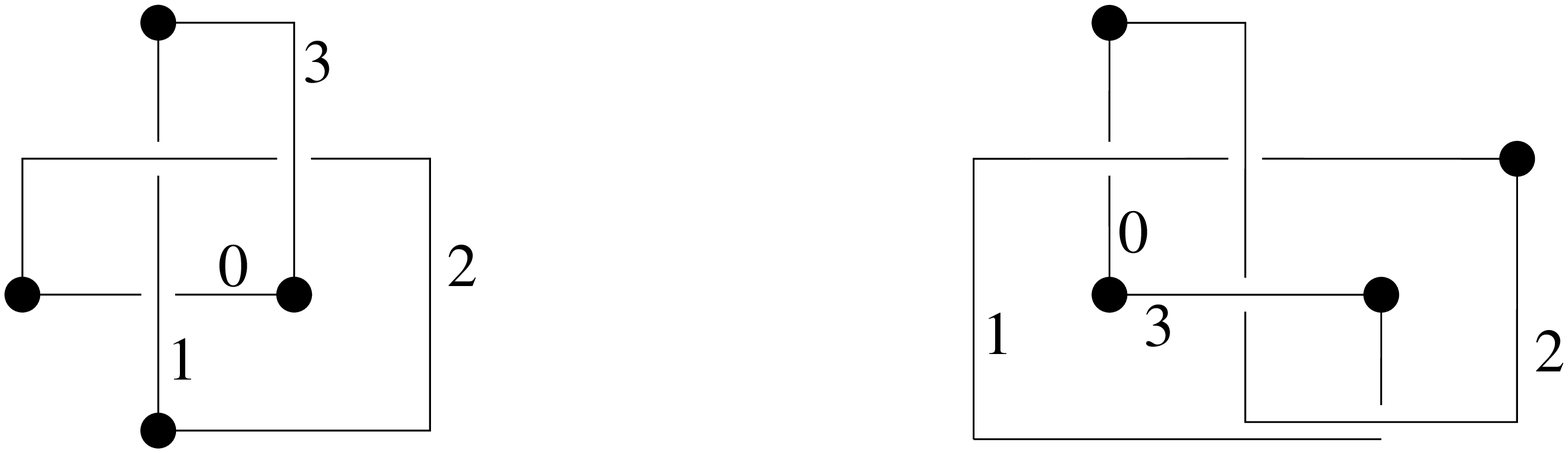,height=3cm}
\end{center}

\begin{center}
Fig. 1.3; A trivial lattice knot, with 4 edges, 4 right angles and
no changes of planes. A lattice trefoil with 24 edges, 12 right angles and
8 changes of planes. A lattice figure-eight knot with 30 edges,
14 right angles and 8 changes of planes.
The numbers are the $z$-levels and the dots are the sticks in
$z$-direction.
\end{center}

\subsection{Early invariants of links}
The fundamental problem in knot theory is\footnote{One should rather say ``was"; there are 
algorithms allowing recognition of any knots, even if very slow. Modern Knot Theory looks rather 
for structures on a space of knots or for a mathematical or physical meaning of knot invariants.} 
to be able to distinguish
non-equivalent knots. It was not achieved (even in the simple case of the
unknot and the trefoil knot) until
Jules Henri Poincar\'e (1854-1912) in his ``Analysis Situs"
paper (\cite{Po-1} 1895) laid foundations for algebraic topology.
According to W.Magnus wrote \cite{Mag}:
{\it Today,  it appears to be a hopeless task to assign priorities for the
definition and the use of fundamental groups in the study of knots,
particularly since Dehn had announced \cite{De} one of the important
results of his 1910 paper (the construction of Poincar\'e spaces with the
help of knots) already in 1907.}. Wilhelm Wirtinger (1865-1945) in his lecture
delivered at a meeting of the German Mathematical Society in 1905 outlined
a method of finding a knot group presentation (it is called now the
Wirtinger presentation of a knot group) \cite{Wir}, but examples using his method 
were given after the work of Dehn.

\subsection{Kirchhoff's complexity of a graph}\label{Subsection 1.4} 
Gustav Robert Kirchhoff (1824-1887) in his fundamental paper on electrical circuits \cite{Kir}.
published in 1847, defined the complexity of a circuit. In the language of graph theory,
this complexity of a graph, $\tau(G)$, is the number of spanning trees of $G$, that is 
 trees in $G$ which contain all vertices of $G$.  It was noted in \cite{BSST} that if $e$ is
an edge of $G$ that is not a loop  then $ \tau(G)$ satisfies the deleting-contracting relation:
$$ \tau(G) = \tau(G-e) + \tau(G/e),$$ where $G-e$ is the graph obtained from $G$ by deleting the edge $e$, 
and $G/e$ is obtained from $G$ by contracting $e$, that is identifying endpoints of $e$ in $G-e$.
The deleting-contracting relation has an important analogue in knot theory, 
usually called a skein relation  (e.g. Kauffman bracket skein relation). 
Connections were discovered only about
a hundred years later (e.g. the Kirchhoff complexity of a circuit corresponds
to the determinant of the knot or link yielded by the circuit, see the next subsection).

For completeness, and to be later to see clearly connection to Goeritz matrix in knot theory, 
let us defined the (version of) the Kirchhoff matrix of a graph, $G$, determinant of which is 
the complexity $\tau(G)$.
\begin{definition}\label{Definition H-1.3} Consider a graph $G$ with vertices $\{v_0,v_1,\ldots,v_n \}$ 
possibly with multiple-edges and loops (however loops are ignored in definitions which follows).
\begin{enumerate}
\item[(1)] The adjacency matrix of the graph $G$ is the $(n+1)\times (n+1)$ matrix $A(G)$ whose 
entries, $a_{ij}$ are equal to the number of edges connecting $v_i$ with $v_j$; we set $v_{i,i}=0$.
\item[(2)] The degree matrix $\Delta(G)$ is the diagonal $(n+1)\times (n+1)$ matrix whose $i$th entry 
is the degree of the vertex $v_i$ (loops are ignored). Thus the $i$th entry is equal to
$-\sum_{j=0}^n a_{ij}$.
\item[(3)] The Laplacian matrix $Q'(G)$ is defined to be $\Delta(G) - A(G)$; \cite{Big}.
Notice, that the sum of rows of $Q'(G)$ is equal to zero and that $Q'(G)$ is a symmetric matrix.
\item[(4)] The Kirchoff matrix (or reduced Laplacian matrix) $Q(G)$ of $G$ is obtained from $Q'(G)$ 
by deleting the first row and the first column from $Q'(G)$.
\end{enumerate}
\end{definition}
\begin{theorem}\label{Theorem H-1.4} 
$det(Q(G)) = \tau(G)$.
\end{theorem}
\begin{proof} The shortest proof, I am aware of, is by direct checking that $det(Q(G))$ satisfies 
deleting-contracting relation for any edge $e$, not a loop, that is
$$det(Q(G)) = det(Q(G-e)) + det(Q(G/e)).$$
The above equation plays an important role in showing in Section 7 that an alternating 
link is a quasi-alternating as well.
\end{proof}

\begin{example}\label{Theorem H-1.5}
Consider the graph {\parbox{1.7cm}{\psfig{figure=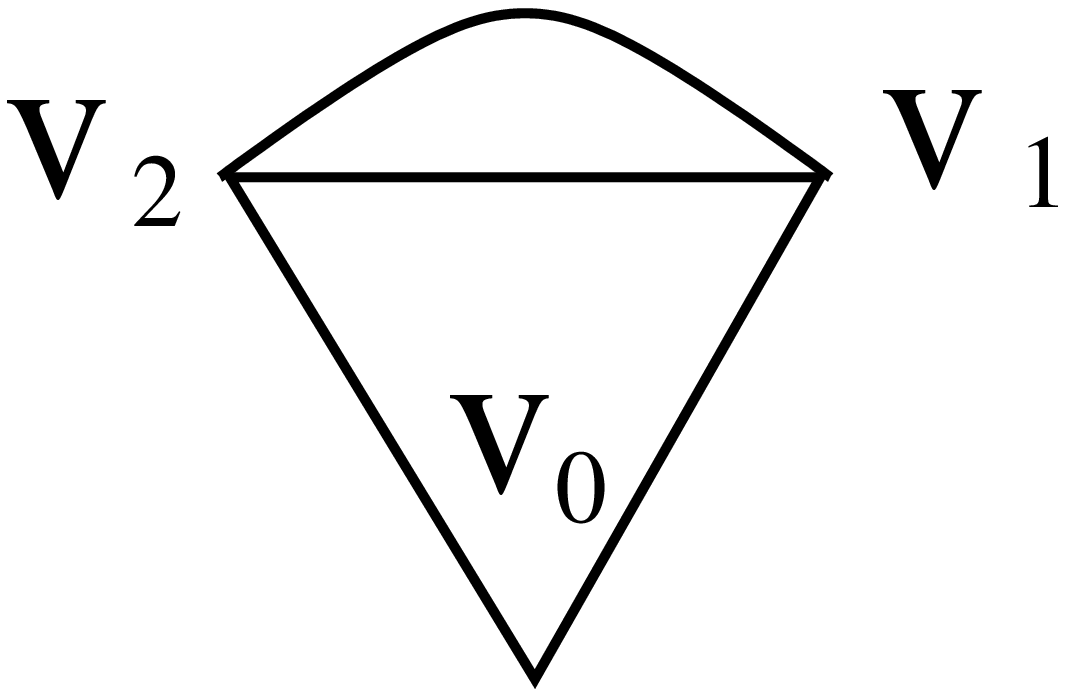,height=1.1cm}}}.\ For this graph we have:
$$  A({\parbox{0.9cm}{\psfig{figure=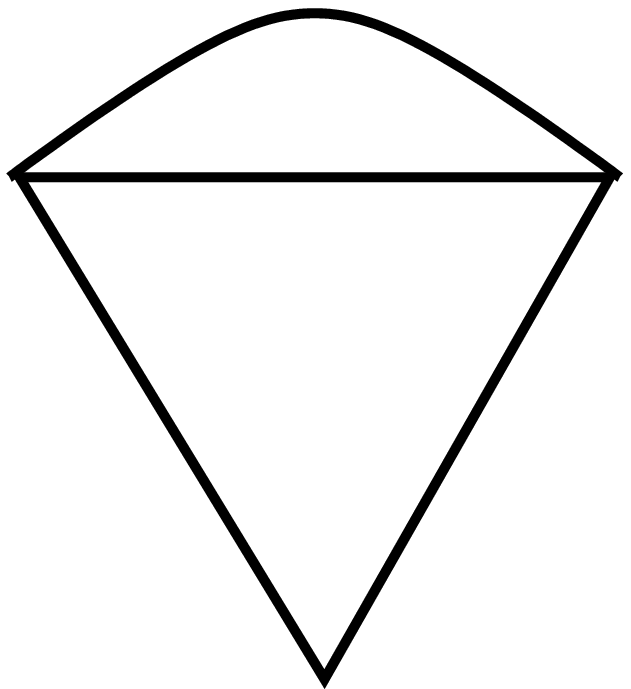,height=1.0cm}}}) = 
 \left[
\begin{array}{ccc}
 0 & 1 & 1\\
 1 & 0 & 2 \\
 1 & 2 & 0
\end{array}
\right],\ \ \Delta({\parbox{0.9cm}{\psfig{figure=Graph-H1.eps,height=1.0cm}}}) = 
\left[
\begin{array}{ccc}
 2 & 0 & 0\\
 0 & 3 & 0 \\
 0 & 0 & 3
\end{array}
\right].
$$
$$ Q'({\parbox{0.9cm}{\psfig{figure=Graph-H1.eps,height=1.0cm}}}) =
 \left[
\begin{array}{ccc}
 2 & -1 & -1\\
 -1 & 3 & -2 \\
 -1 & -2 & 3
\end{array}
\right]; \ Q({\parbox{.9cm}{\psfig{figure=Graph-H1.eps,height=1.0cm}}}) =
 \left[
\begin{array}{cc}
  3 & -2 \\
 -2 &  3
\end{array}
\right]. $$
$$ det(Q({\parbox{.9cm}{\psfig{figure=Graph-H1.eps,height=1.0cm}}}))=
det\left[
\begin{array}{cc}
  3 & -2 \\
 -2 &  3
\end{array}
\right] = 5 = \tau({\parbox{.9cm}{\psfig{figure=Graph-H1.eps,height=1.0cm}}}).
$$

As we will see in the next subsection the corresponding knot is the figure eight knot (Fig. 1.4).
\end{example}
\subsection{Tait's relation between knots and graphs.}
   Tait was the first to notice the relation between knots and
planar graphs. He colored the regions of the knot diagram alternately white
and black (following Listing) and constructed the graph by placing a vertex
inside each white region, and then connecting vertices by edges going through
the crossing points of the diagram (see Figure 1.4)\cite{D-H}.

\ \\
\centerline{\psfig{figure=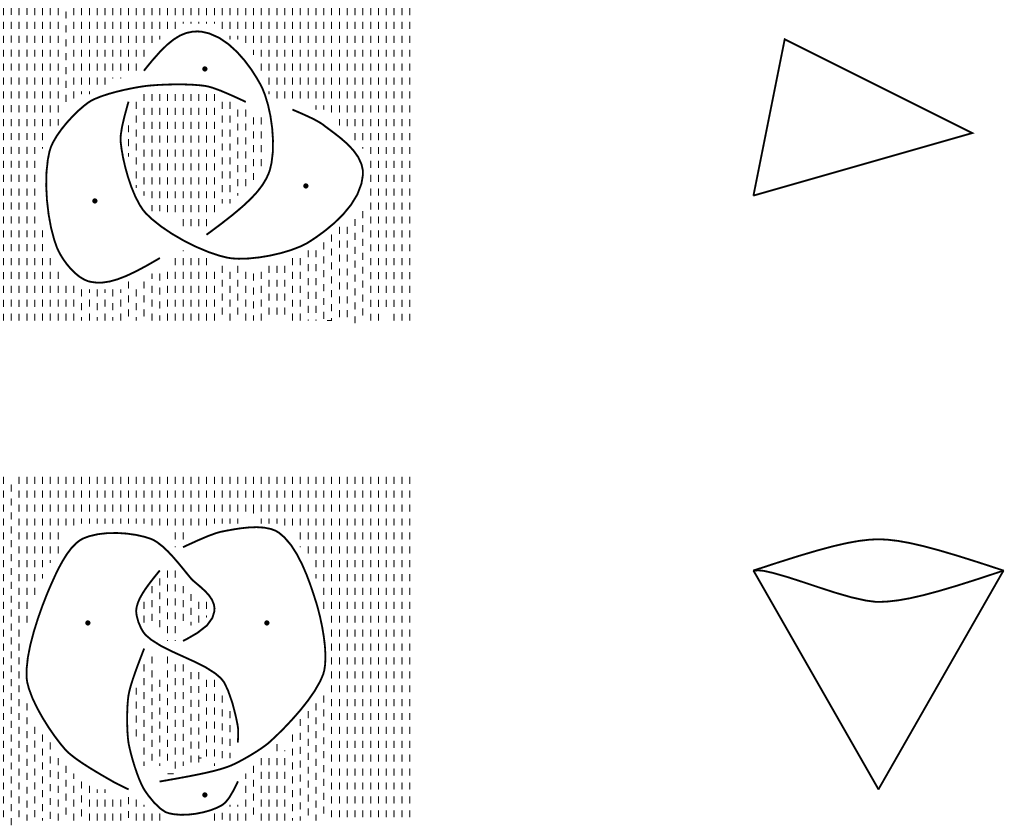,height=10.5cm}}
\begin{center}
  Figure 1.4; Tait's construction of graphs from link diagrams, according to Dehn-Heegaard
\end{center}
\ \\

It is useful to mention the Tait construction going in the opposite direction, from 
a signed planar graph, $G$ to a link diagram $D(G)$. We replace every edge of a graph by a crossing 
according to the convention of Figure 1.5 and connect endpoints along edges as in 
Figures 1.6 and 1.7.

\ \\
\centerline{{\psfig{figure=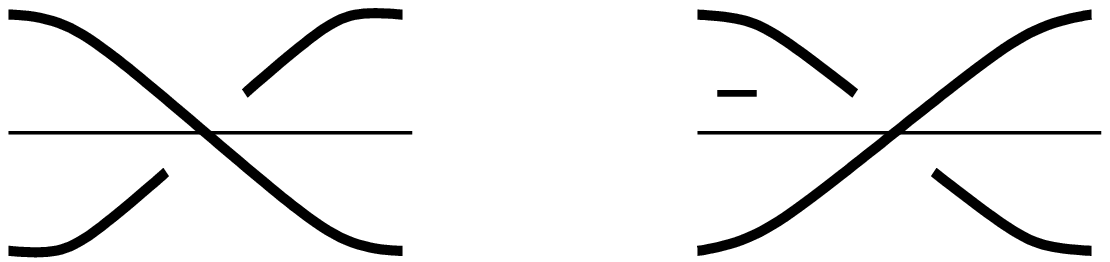,height=2.7cm}}}
\begin{center}
Fig. 1.5;\ convention for crossings of signed edges (edges without markers are 
assumed to be positive)

\end{center}

\ \\
\centerline{{\psfig{figure=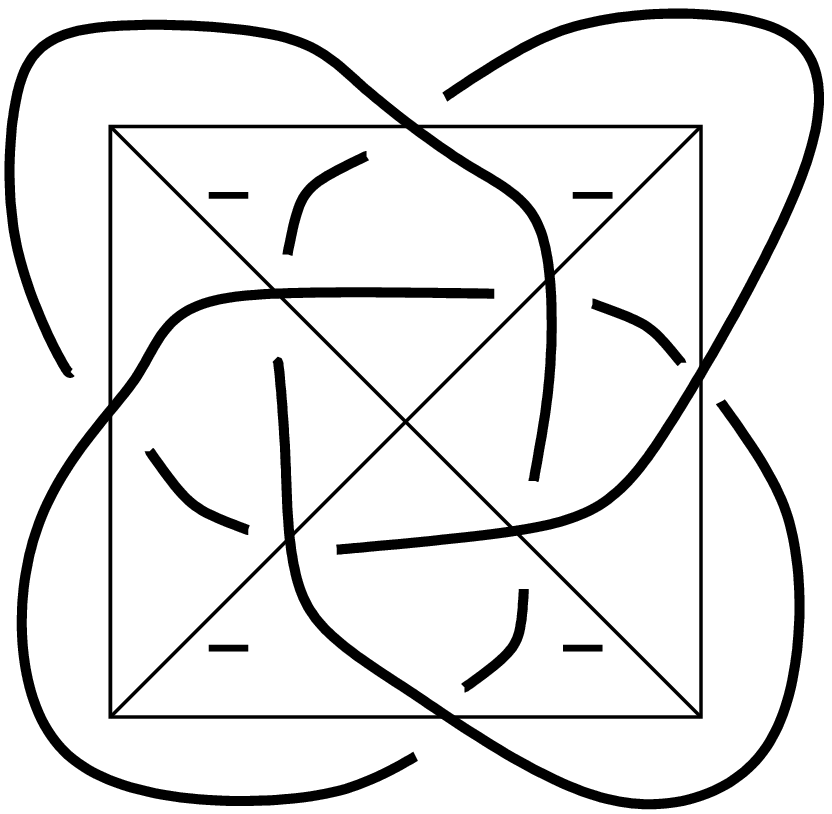,height=3.7cm}}}
\begin{center}
Fig. 1.6;\ The knot $8_{19}$ and its Tait graph ($8_{19}$ is the first in tables non-alternating knot)

\end{center}

\ \\
\centerline{{\psfig{figure=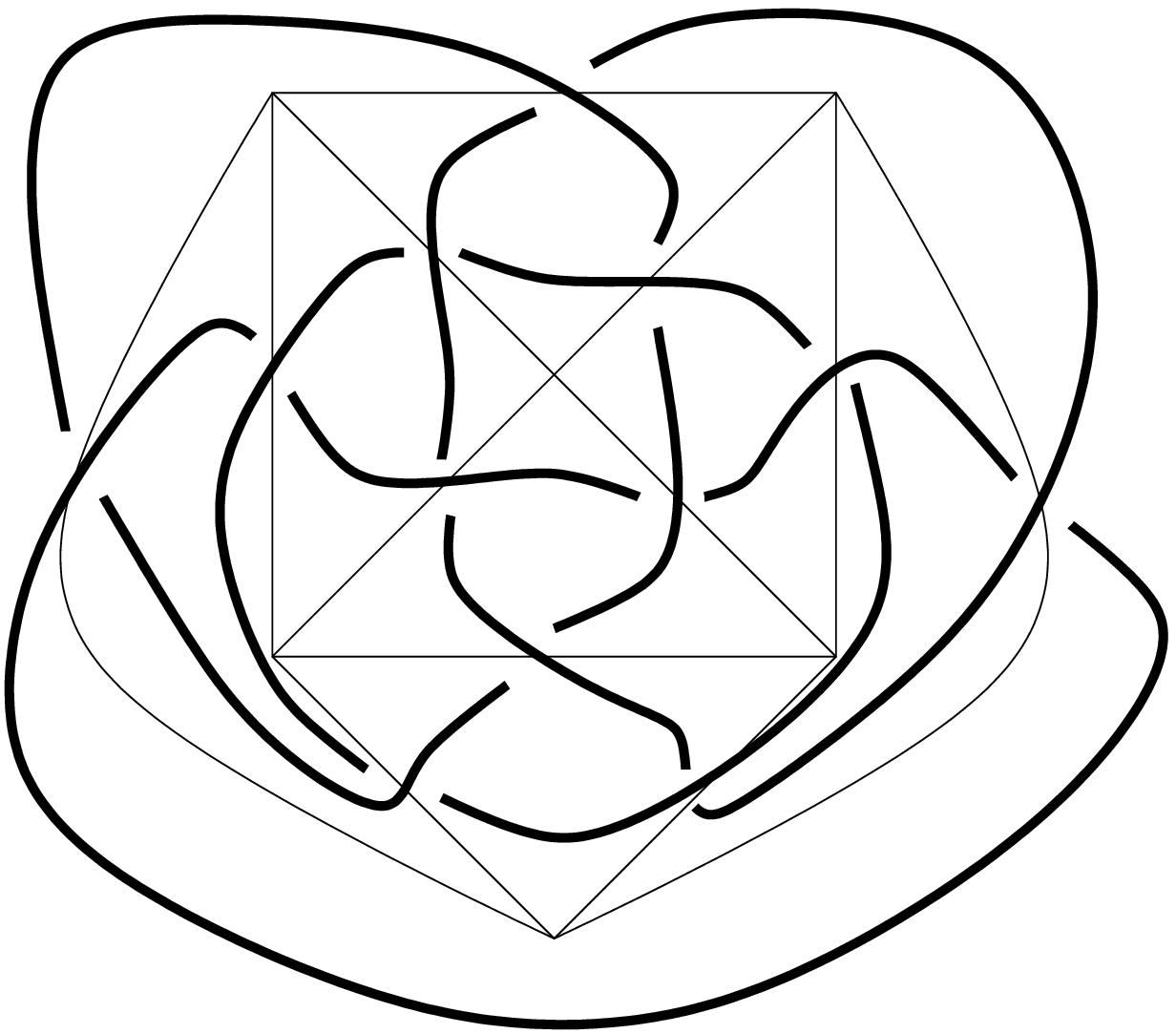,height=3.7cm}}}
\begin{center}
Fig. 1.7;\ Octahedral graph (with all positive edges) and
the associated link diagram
\end{center}

We should mention here one important observation known already to Tait (and in explicit form to Listing):
\begin{proposition}\label{Proposition 1.6}
The diagram $D(G)$ of a connected graph $G$ 
is alternating if and only if $G$ is positive (i.e. all edges of $G$ are positive) or 
$G$ is negative.
\end{proposition} 
A proof is illustrated in Figure 1.8. \\ \ \\
\centerline{{\psfig{figure=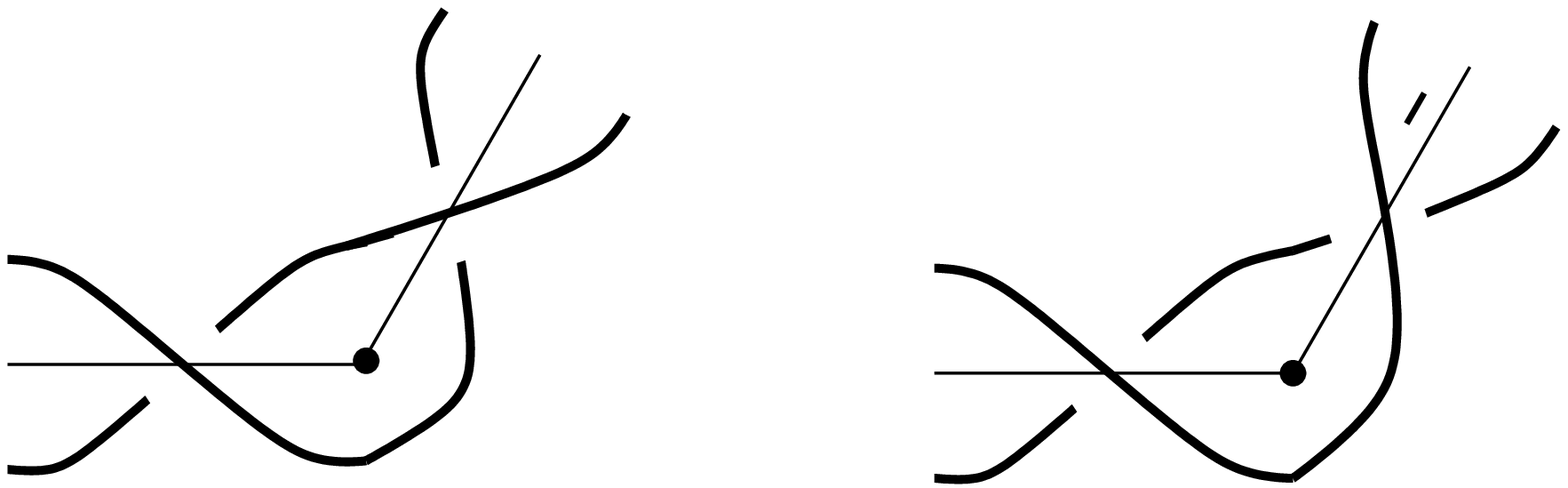,height=3.7cm}}}
\begin{center}
\ \\
Fig. 1.8;\ Alternating and  non-alternating parts of a diagram
\end{center}

\subsection{Link diagrams and Reidemeister moves}
In this part we define, after Reidemeister, a polygonal knot and link, and
$\Delta$-equivalence of knots and links. A $\Delta$-move is an
elementary deformation of a polygonal knot which intuitively agrees
with the notion of ``deforming without cutting and glueing" which
is the first underlining principle of topology. 

\begin{definition}[Polygonal knot, $\Delta$-equivalence].
\begin{enumerate}
\item[(a)]
A polygonal knot is a simple closed polygonal curve in $R^3$.
\item[(b)]
Let us assume that $u$ is a line segment (edge) in a
polygonal knot $K$ in $R^3$.
Let $\Delta$ be a triangle in $R^3$ whose boundary consists
of three line segments $u,\ v,\ w$ and such
that $\Delta \cap L=u$. The polygonal curve defined as
$K'=(K-u)\cup v\cup w$ is a new polygonal knot in $R^3$.
We say that the knot $K'$ was obtained from $K$ by a $\Delta$-move.
Conversely, we say that $L$ is obtained from $L'$ by a
${\Delta}^{-1}$-move (Fig. 1.9).
We allow the triangle $\Delta$ to be degenerate so that the vertex
$v\cap w$ is on the side $u$; in other words we allow subdivision
of the line segment $u$.\footnote{Notice, that any subdivision is a combination of three non-degenerate
$\Delta$-moves, or more precisely two $\Delta$-moves and the inverse to a $\Delta$-move:\\
\psfig{figure=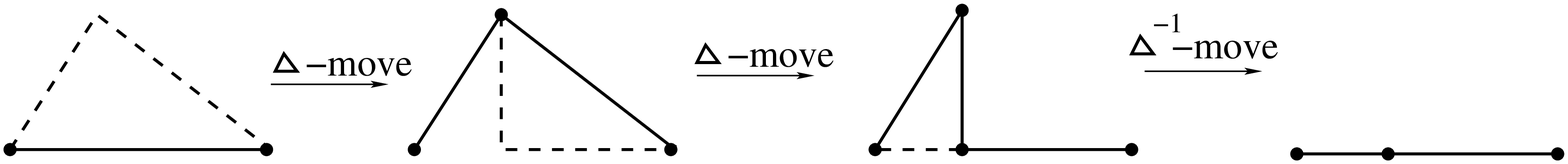,height=1.3cm}}
\item[(c)]
We say that two polygonal knots are $\Delta$-equivalent (or
combinatorially equivalent) if one can be obtained from the other
by a finite sequence of ${\Delta}$- and ${\Delta}^{-1}$-moves.
\end{enumerate}
\end{definition}

\centerline{\psfig{figure=Deltaruch.eps}}
\begin{center}
Fig. 1.9 
\end{center}

Polygonal links are usually presented by their projections to a plane. Let
$p:R^3 \to R^2$ be a projection and let $L\subset R^3$ be a link.
Then a point $P\in p(L)$ is called a multiple point (of $p$)
if $p^{-1}(P)$ contains more than one point (the number of points
in $p^{-1}(P)$ is called the multiplicity of $P$).

\begin{definition}\label{Definition H-1.7}
The projection $p$ is called regular if
\begin{enumerate}
\item
[(1)] $p$ has only a finite number of multiple points
and all of them are of multiplicity two,
\item
[(2)] no vertex of the polygonal link
is an inverse image of a multiple point of $p$.
\end{enumerate}
Thus in case of a regular projection  the parts of a diagram,
illustrated in  the figure below,  are not allowed.
\end{definition}
\vskip 0.1in
\centerline{\psfig{figure=osobliwosci.eps}
\ \ \ \psfig{figure=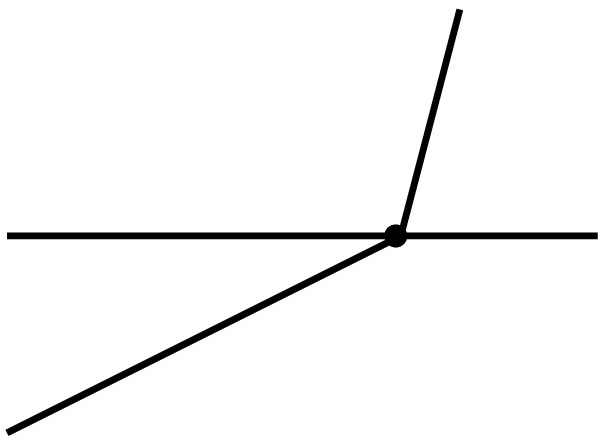,height=1.9cm}}
\ \\

Maxwell was the first person to consider the question when two projections represent 
equivalent knots. He considered some elementary moves (reminding future Reidemeister moves),
but never published his findings.

The formal interpretation of $\Delta$-equivalence of knots in terms of diagrams.
Was done by Reidemeister \cite{Re-1}, 1927, and Alexander and Briggs
\cite{A-B}, 1927.

\begin{theorem} [Reidemeister theorem]
\label{Theorem H-1.8}
\ \\
Two link diagrams are $\Delta$-equivalent\footnote{In modern Knot Theory, especially after the 
work of R.~Fox, we use 
usually the equivalent notion of ambient isotopy in $R^3$ or $S^3$. 
Two links in a 3-manifold $M$ are ambient isotopic if there is 
an isotopy of $M$ sending one link into another.}
 if and only if they are connected by a
finite sequence of Reidemeister moves $R_i^{\pm 1}, i=1,2,3$ (see Fig.~1.10)
and isotopy (deformation) of the plane of the diagram. The theorem holds also
for oriented links and diagrams. One then has to take into account
all possible coherent orientations of diagrams involved in the moves.
\end{theorem}

\centerline{\psfig{figure=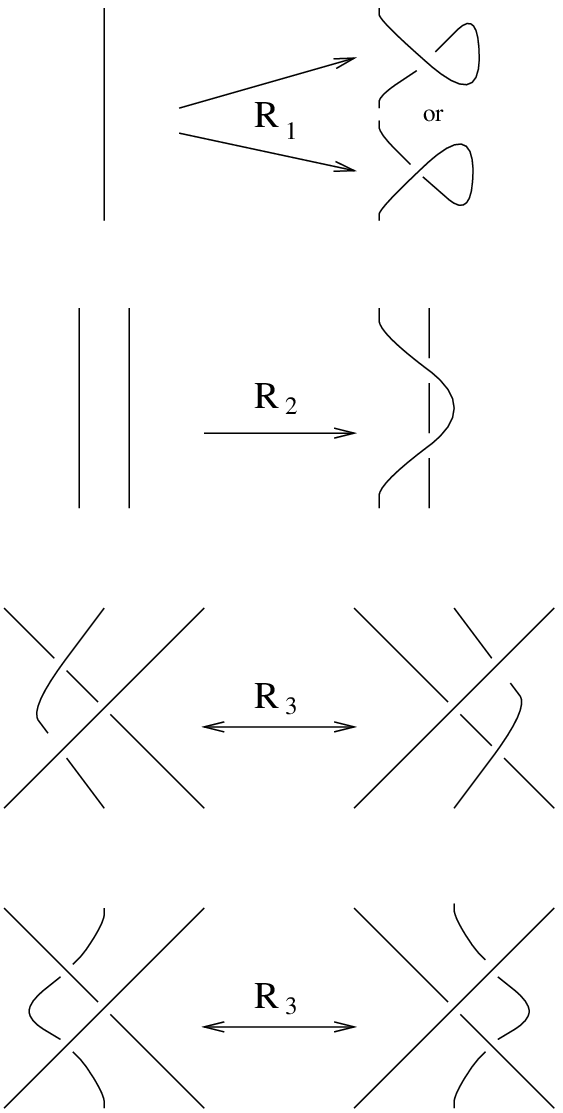}}
\begin{center}
Figure 1.10; Reidemeister moves; we draw two versions of the first and the third moves which are related by a 
mirror symmetry in the plane of the projection
\end{center}

\section{Goeritz matrix and signature of a link}\label{Section H2}
\markboth{\hfil{\sc Chapter IV. Goeritz and Seifert matrices }\hfil}
{\hfil{\sc Goeritz matrix and signature of a link}\hfil}

In the first half of XX-century combinatorial methods ruled over knot theory,
even if more topological approach was possible, for example, Reidemeister moves 
were used to prove existence of the Alexander polynomial even if purely 
topological prove using the fundamental group was possible and probably 
well understood by Alexander himself. Later, after the Second World War, to great extend 
under influence of Ralph Hartzler Fox (1913 -1973), Knot Theory was considered to be a 
part of algebraic topology with fundamental group and coverings playing an important role.
The renaissance of combinatorial methods in Knot Theory can
be traced back to Conway's paper \cite{Co-1} and bloomed after
the Jones breakthrough \cite{Jo-1} with Conway type invariants and
Kauffman approach (compare Chapter III of \cite{P-Book}). As we already mentioned,
 these had their predecessors in $1930$th \cite{Goe,Se}.
Goeritz matrix of a link can be defined purely combinatorially and
is closely related to Kirchhoff matrix of an electrical network.
Seifert matrix is a generalization of the Goeritz matrix and, even
historically, its development was mixing combinatorial and topological
methods.

In this section we start from the work of L.~Goeritz.
He showed \cite{Goe} how to associate a quadratic form
to a diagram of a link and moreover how to use this form
to get algebraic invariants of the knot
(the signature of this form, however, is not an invariant of the knot).
Later, H.~F.~Trotter \cite{Tro-1}, using Seifert form (see
Section 3), introduced another quadratic form, the signature of which 
was an invariant of links.

C.~McA.~Gordon and R.~A.~Litherland \cite{G-L} provided a unified
approach to Goeritz and Trotter forms. They showed how to use
the form of Goeritz to get (after adding a correcting factor)
the signature of a link (this signature is often 
called a classical or Trotter, or Murasugi \cite{M-10} signature of a link).

We begin with a purely combinatorial description
of the matrix of Goeritz and of the signature of a link.
This description is based on \cite{G-L} and
\cite{Tral-1}.

\begin{definition}\label{c3:1.1}
Let $L$ be a diagram of a link.
Let us checkerboard color the complement of the diagram in the projection 
plane $R^2$, that is, color in black and white 
the regions into which the plane is divided by the diagram\footnote{This 
 (checkerboard) coloring was first 
used by P.~G.~Tait in 1876/7, compare Chapter II of \cite{P-Book}, however we 
switched, after C.~Gordon, the role of white and black. We can say that Tait 
convention worked well with a blackboard, while our convention with white-board.}. 
 We assume that the unbounded region of $R^2\setminus L$
is colored white and it is denoted by $X_0$ while the other 
white regions are denoted by $X_1, \ldots, X_n$. 
Now, to any crossing, $p$, of $L$ we associate the number $\eta (p)$ 
which is either $+1$ or $-1$
according to the convention described in Fig.~2.1. \\
\ \\
\begin{center}
\begin{tabular}{c} 
\includegraphics[trim=0mm 0mm 0mm 0mm, width=.5\linewidth]
{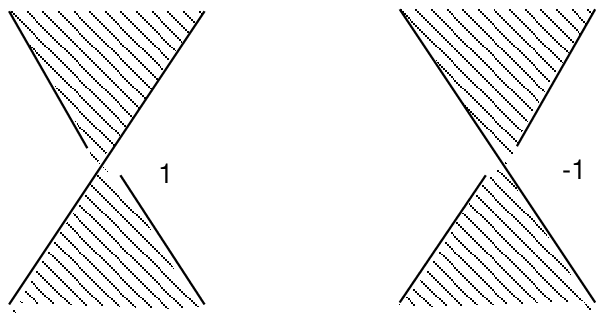}\\
\end{tabular}
\\
Fig. 2.1
\end{center}

%

Let $G' = \{ g_{i,j}\}^n_{i,j=0}$, where
$$g_{i,j} = \left\{
\begin{array}{ll}
-\sum_{p}\eta (p)
& \mbox{ for $i\neq j$, where the summation extends }\\
 \ \ &\mbox{ over crossings which connect $X_i$ and $X_j$}\\
-\sum_{k=0,1,\ldots,n; k\neq i} g_{i,k}&\mbox{ if $ i=j$}\\
\end{array}
\right.
$$

The matrix $G'=G'(L)$ is called the unreduced Goeritz matrix 
of the diagram $L$.
The reduced Goeritz matrix (or shortly Goeritz matrix) associated to 
the diagram $L$ is
the matrix $G=G(L)$ obtained by removing the first
row and the first column of $G'$.
\end{definition}

\begin{theorem}[\cite{Goe,K-P,Ky}.]\label{c3:1.2}
Let us assume that $L_1$ and $L_2$ are two diagrams of a given link.
Then the matrices $G(L_1)$ and $G(L_2)$ can be obtained one from the other
in a finite number of the following elementary 
operations on matrices:
\begin{enumerate}

\item $G\Leftrightarrow P GP^T$, where $P$ is a
matrix with integer entries and $\det P=\pm 1$.

\item 
$$G \Leftrightarrow
\left[
\begin{array}{cc}
G&0\\
0&\pm 1\\
\end{array}
\right]
$$

\item 
$$G \Leftrightarrow
\left[
\begin{array}{cc}
G&0\\
0& 0\\
\end{array}
\right]
$$
\end{enumerate}
Moreover, if $L$ is a diagrams of a knot\footnote{It suffices to 
assume that $L$ represent a non-split link, that is a link all 
projections of which are connected.}, then operations
(1) and (2) are sufficient.
\end{theorem}

\begin{corollary}\label{c3:1.3}
$|\det G|$ is an invariant of isotopy of knots called the determinant 
of a knot\footnote{Often, by the determinant of a knot one understands 
the more delicate invariant whose absolute value is equal to $|\det G|$; 
see Corollary \ref{Corollary H-2.7}. This determinant can be defined, also, 
as the Alexander-Conway or Jones polynomial at $t=-1$; compare Corollary \ref{Corollary H-6.18}.}.

\end{corollary}

A sketch of a proof of Theorem~\ref{c3:1.2}

We have to examine how a Goeritz matrix changes under Reidemeister
moves. The matrix does not depend on the orientation of the link,
let us assume, however, that the diagram $L$ is oriented.
We introduce new notation: a crossing is called of type I or II
according to Fig.~2.2. Moreover, we define
$\mu (L) = \sum \eta (p)$, where the summation is taken over crossings
of type II.  \\ 
\ \\

\begin{center}
\begin{tabular}{c} 
\includegraphics[trim=0mm 0mm 0mm 0mm, width=.5\linewidth]
{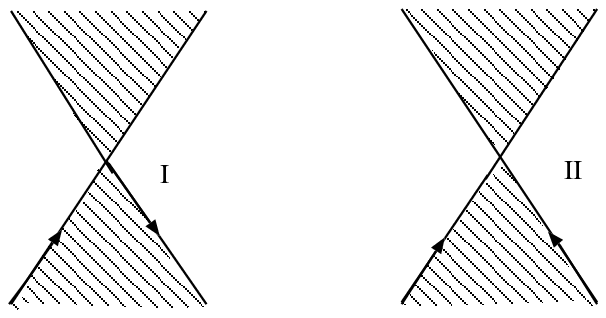}\\
\end{tabular}
\\
Fig. 2.2
\end{center}


Now let us construct a graph with vertices representing black regions
(this is the Tait's construction, however, the choice of black and white regions is reversed) 
and edges in bijection with crossings of $L$. 
Edges of the graph are in bijection with crossings of $L$: two vertices of 
the graph are joined
if and only if the respective regions meet in a 
crossing\footnote{This construction of Tait is an important motivation
for material in Chapter V of \cite{P-Book}. The constructed graph, which we denote by 
$G_b(L)$, is usually called the Tait graph of $L$ (see the first section). 
For an alternating diagram $L$ this 
graph is the same as the graph $G_{s_+}(L)$ considered in Chapter V of \cite{P-Book}. 
We often equip the edges
of $G_b(L)$ with signs: the edge corresponding to a vertex $p$ has the sign 
$\eta (p)$ (see Figure 2.1). 
The signed graph $G_b(L)$ is considered in Chapter V of \cite{P-Book}; compare also Definition 7.4.}.
Let  $B(L)$ denote the number of components of such a graph.
From now on, let $R$ be a Reidemeister move.
We denote by $G_1$ the Goeritz matrix of $L$, and  by $G_2$ the matrix 
of $R(L)$. Similarly we set
$\mu_1 = \mu (L)$, $\mu_2 = \mu (R(L))$ and also 
$\beta_1=B(L)$, $\beta_2 = B(R(L))$. We will write
$G_1\approx G_2$ if  $G_1$ and $G_2$ are in relation (1)
and  $G_1\sim G_2$ if $G_2$ can be obtained from $G_1$ 
by a sequence of relations (1)--(3).

\begin{enumerate}

\item Let us consider the first Reidemeister move $R_1$.

\begin{enumerate}

\item In the case shown in Fig.~2.3 we have: 
$\beta_1 = \beta_2$, $\mu_1 = \mu_2$ and
$G_1\approx G_2$.  \\ 
\ \\

\begin{center}
\begin{tabular}{c} 
\includegraphics[trim=0mm 0mm 0mm 0mm, width=.5\linewidth]
{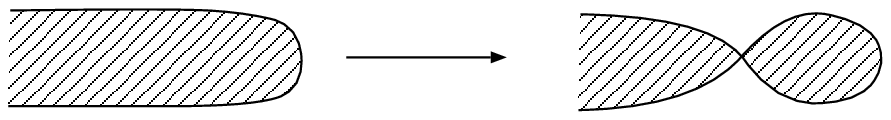}\\
\end{tabular}
\\
Fig. 2.3
\end{center}


\item In the case shown in Fig.~2.4 we have: 
$$\beta_1=\beta_2,\  \mu_2=\mu_1 +\eta(p),\ 
G_2 = \left[
\begin{array}{cc}
G_1&0\\
0&\eta(p)\\
\end{array}
\right]$$

\ \\

\begin{center}
\begin{tabular}{c} 
\includegraphics[trim=0mm 0mm 0mm 0mm, width=.5\linewidth]
{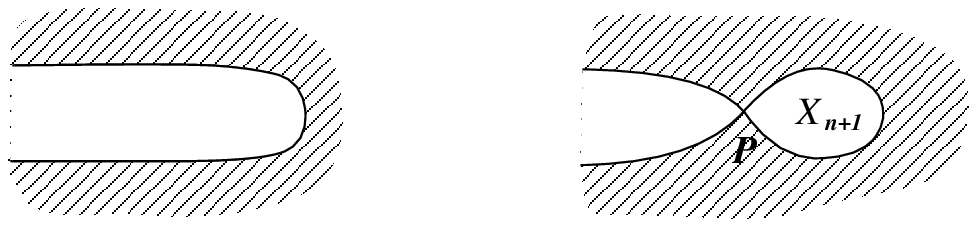}\\
\end{tabular}
\\
Fig. 2.4
\end{center}


\end{enumerate}

\item Let us consider the second Reidemeister move $R_2$. 

\begin{enumerate}

\item In the case described in Fig.~2.5 we get immediately that 
$\beta_1 = \beta_2$ and $\mu_1 = \mu_2$ 
(either both crossings are of type I or of type II
and always of opposite signs), $G_1\approx G_2$.  \\ 
\ \\

\begin{center}
\begin{tabular}{c} 
\includegraphics[trim=0mm 0mm 0mm 0mm, width=.5\linewidth]
{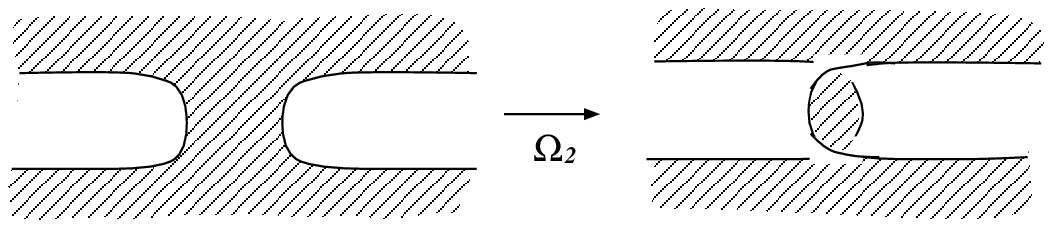}\\
\end{tabular}
\\
Fig. 2.5
\end{center}


\item In the case described in Fig.~2.6 we have to consider two
subcases. In each of them $\mu_1 = \mu_2$, since the two
new crossings are either both of type I or both of type II
and always of opposite signs: \\ 
\ \\

\begin{center}
\begin{tabular}{c} 
\includegraphics[trim=0mm 0mm 0mm 0mm, width=.5\linewidth]
{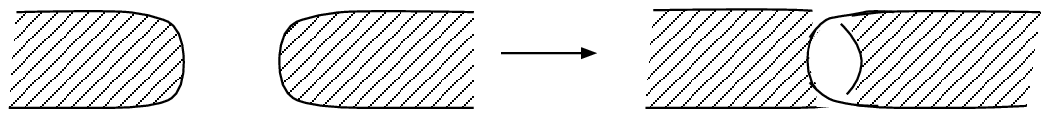}\\
\end{tabular}
\\
Fig. 2.6
\end{center}


\begin{enumerate}

\item[(i)]  $\beta_1 = \beta_2$. Then
$$G_2 \approx
\left[
\begin{array}{ccc}
G_1&&0\\
&1&\\
0&&-1\\
\end{array}
\right]
\mbox{ or }
\left[
\begin{array}{cc}
G_2&0\\
0 &1\\
\end{array}
\right]
\approx
\left[
\begin{array}{cccc}
G_1&&&0\\
&1&&\\
&&1&\\
0&&&-1\\
\end{array}
\right]
$$
We leave it for the reader to check, c.f.~\cite{K-P}. 

Both possibilities give $G_1\sim G_2$.

\item [(ii)] $\beta_2 = \beta_1 - 1$. Then we see immediately that
$$G_2 \approx \left[
\begin{array}{cc}
G_1& {\bf 0}\\
{\bf 0} &0\\
\end{array}
\right].
$$
\end{enumerate}
\end{enumerate}

\item Let us consider the Reidemeister move $R_3$ (Fig.~2.7).  \\ 
\ \\

\begin{center}
\begin{tabular}{c} 
\includegraphics[trim=0mm 0mm 0mm 0mm, width=.5\linewidth]
{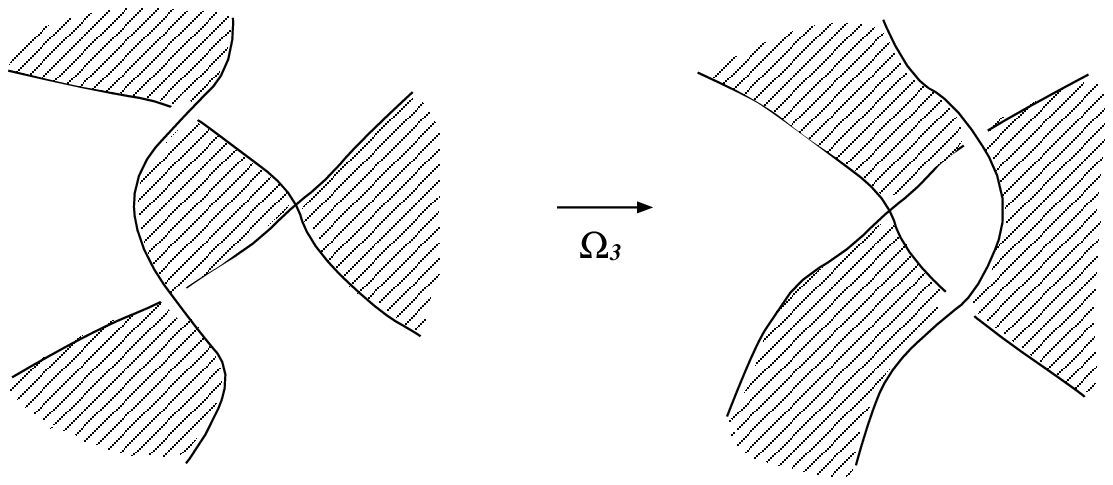}\\
\end{tabular}
\\
Fig. 2.7
\end{center}


We see immediately that $\beta_1 = \beta_2$. Next we should 
consider different orientations of arcs participating in $R_3$
and two possibilities for the crossing $p$. 
However, we will get always $\mu_2 = \mu_1+\eta (p)$ and 
$$
G_2\approx
\left[
\begin{array}{cc}
G_1&0\\
0&\eta (p)\\
\end{array}
\right].
$$
We leave it for the reader to check (c.f.~ \cite{Goe} and
\cite{Re-2}).
\end{enumerate}

This concludes the proof of Theorem \ref{c3:1.2}.

\begin{corollary}\label{Corollary H-2.4}
\begin{enumerate}
\item[(1)] For a link $L$ 
let us define $\sigma (L)= \sigma (G(L)) -\mu (L)$, where $\sigma (G(L))$
is the signature of the Goeritz matrix of $L$.
Then $\sigma (L)$ is an invariant of the link $L$, 
called the signature of the link; compare Corollary \ref{Corollary H-2.7}
 and Definition \ref{Definition IV.5.7}.

\item[(2)]
 Let us define $\nul (L) =\nul (G(L))+\beta (L) -1$, where $\nul (G(L))$
is the nullity (i.e.~the difference between the dimension and the rank)
of the matrix $G(L)$. Then $\nul (L)$ is an invariant 
of the link $L$ and we call it the nullity (or defect) of the link.
\end{enumerate}
\end{corollary}
Proof. It is enough to apply Theorem \ref{c3:1.2} to see that
$\sigma (L)$ and $\nul (L)$ are invariant with respect to
Reidemeister moves.

L.~Traldi \cite{Tral-1} introduced a modified matrix of an oriented link,
the signature and the nullity of which are invariants of the link.

\begin{definition}\label{c3:1.5}
Let $L$ be a diagram of an oriented link.
Then we define the generalized Goeritz matrix 
$$H(L) = \left[
\begin{array}{ccc}
G&&\bigcirc\\
&A&\\
\bigcirc&&B\\
\end{array}
\right],$$
where $G$ is a Goeritz matrix of $L$, and the matrices $A$ and $B$
are defined as follows. The matrix $A$ is diagonal of dimension 
equal to the number of type II crossings and the diagonal entries equal
to $-\eta(p)$, where $p$'s are crossings of type II. The matrix $B$ is
of dimension $\beta(L)-1$ with all entries equal to 0.
\end{definition}

\begin{lemma}[\cite{Tral-1}.]\label{c3:1.6}
If $L_1$ and $L_2$ are diagrams of two isotopic oriented links
then $H(L_1)$ can be obtained from $H(L_2)$ 
by a sequence of the following elementary equivalence operations:

\begin{enumerate}

\item $H \Leftrightarrow P HP^T$, where $P$ is a matrix with integer 
entries and with $\det P =\pm 1$,

\item $$ H \Leftrightarrow
\left[
\begin{array}{ccc}
H&& \bigcirc\\
&1&\\
\bigcirc &&-1\\
\end{array}
\right].
$$
\end{enumerate}
\end{lemma}

Proof. Lemma \ref{c3:1.6} follows immediately from the proof
of Theorem \ref{c3:1.2}.

\begin{corollary}\label{Corollary H-2.7}
The determinant $\det (iH(L))$ ($i=\sqrt{-1}$)
is an isotopy invariant of a link $L$, called the determinant of 
the link, $Det_L$.
Moreover, $\sigma(H(L)) = \sigma (L)$ and $\nul (H(L)) = \nul (L)$. 
\end{corollary}

The proof follows immediately from Lemma \ref{c3:1.6} and
from the proof of Theorem \ref{c3:1.2}.

\begin{example}\label{c3:1.8}
Consider a torus link of type $(2,k)$, we denote it by
$T_{2,k}$. It is a knot for odd $k$ and a link of two components 
for $k$ even; see Fig.~2.8. \\
\ \\

\begin{center}
\begin{tabular}{c} 
\includegraphics[trim=0mm 0mm 0mm 0mm, width=.5\linewidth]
{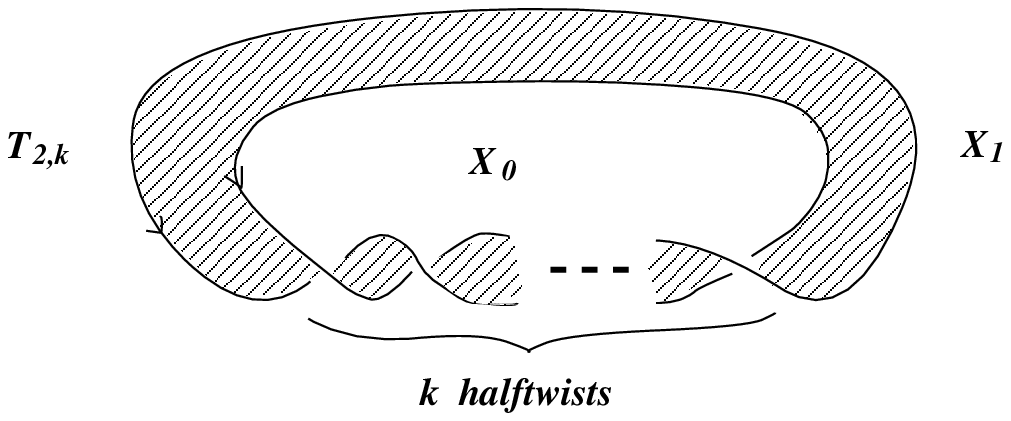}\\
\end{tabular}
\\
Fig. 2.8
\end{center}

%

The matrix $G'$ of $T_{2,k}$
 is then equal to $\left[\begin{array}{cc}k&-k\\-k&k\\
\end{array}\right]$, and thus Goeritz matrix of the link is $G = [k]$. 
Moreover, $\beta = 1$ and $\mu = k$ because all crossings are of type II.
Therefore, for $k\neq 0$, $\sigma (T_{2,k}) = \sigma (G)-\mu = 1-k$ and $\nul (T_{2,k}) =
\nul (G) = 0$. The generalized Goeritz matrix $H$ of the knot
$T_{2,k}$ is of dimension $k+1$ and it is equal to

$$ H= 
\left[
\begin{array}{ccccc}
k&&&&\bigcirc\\
&-1&&&\\
&&-1&&\\
&&&\ddots&\\
\bigcirc&&&&-1\\
\end{array}
\right]
$$
Therefore $Det_L= \det (iH) = (-1)^k  i^{k+1}k=i^{1-k}k$. Notice also that 
$i^{\sigma (T_{2,k})} = \frac{Det_{T_{2,k}}}{|Det_{T_{2,k}}|}$; compare 
Exercise \ref{Exercise 1.10}.

Let us note that if we connect black regions 
of the plane divided by the diagram 
of the link by half-twisted bands
(as indicated in Fig.~2.9)
then we get a surface in $R^3$ (and in $S^3$), 
the boundary of which is the given link; we denote this surface by $F_b$, 
and call the Tait surface of a link diagram; compare Definition \ref{Definition PPS}.
If, for some checkerboard coloring of the plane, the constructed surface
has an orientation which yields 
the given orientation of the link then
this oriented diagram is called a special diagram.\\
\ \\

\begin{center}
\begin{tabular}{c} 
\includegraphics[trim=0mm 0mm 0mm 0mm, width=.5\linewidth]
{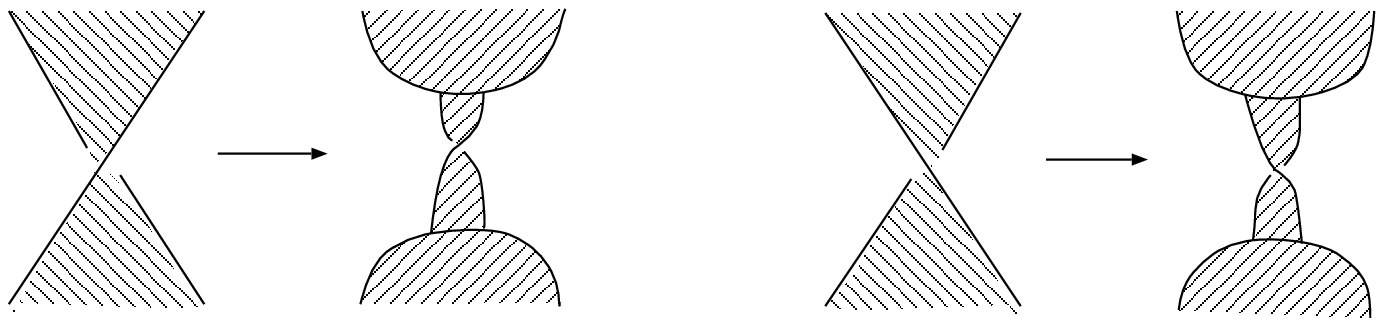}\\
\end{tabular}
\\
Fig. 2.9
\end{center}

%

\end{example}

\begin{exercise}\label{Exercise 1.9}
Prove that an oriented diagram of a link is special if and only if
all crossings are of type I for some checkerboard coloring of the plane.
Conclude from this that for a special diagram $D$, we have 
$\sigma (D)= \sigma (G(D))$. 
\end{exercise}

\begin{exercise}\label{Exercise 1.10}
Show that any oriented link has a special diagram. Conclude from this that 
for any oriented link $L$ one has \\
$Det_{L} = i^{\sigma (L)}|Det_{L}|$; compare Lemma 6.16 and Corollary \ref{Corollary H-6.18}.

\end{exercise}

Assume now that $L_0$ is a sublink of an oriented link $L$. Let $L'$ be 
an oriented link obtained from $L$ by changing the orientation of $L_0$ 
to the opposite orientation. Let $D_L$ be a diagram of $L$ and 
define $lk(L-L_0,L_0)$ be defined a $\sum_{p}\sgn p$ where the sum is taken 
over all crossings of the diagram of $L-L_0$ and $L_0$ (as subdiagrams of 
$L_D$. This definition does not depend on the choice of $D_L$, as checked 
using Reidemeister moves and agrees with the standard notion of linking number 
as defined recalled in the next section.

From Corollary \ref{Corollary H-2.4} and Corollary \ref{Corollary H-2.7},  we obtain.
\begin{proposition}[\cite{M-11}]\label{Proposition 1.11}\ 
\begin{enumerate}
\item[(i)] $Det_{L'} = (-1)^{lk(L-L_0,L_0)}Det_L$.
\item[(ii)] $\sigma(L') = \sigma(L) + 2lk(L-L_0,L_0)$.
\item[(ii)] $\sigma(L) + lk(L)$ is independent on orientation of $L$.
\end{enumerate}
\end{proposition} 
\begin{proof} The derivation of formulas is immediate but it is still 
instructive to see how Corollary 2.11(ii) follows from Corollary 2.4(1):
$$\sigma(L')= \sigma(G(L'))-\mu (L')=\sigma(G(L)) - \mu (L')=\sigma(L)+\mu (L) - \mu (L') = 
\sigma(L) + 2lk(L-L_0,L_0).$$
\end{proof}
Recall (\cite{P-2}) that an $n$-move is a local change of an unoriented 
link diagram described in Figure 2.10.

\begin{center}
\psfig{figure=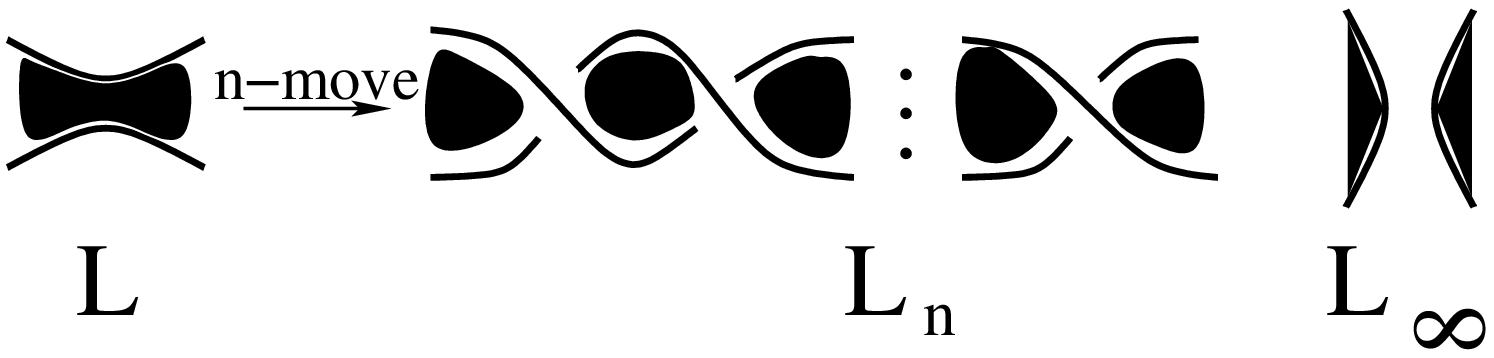,height=3.4cm}
\end{center}
\centerline{Fig. 2.10; $L_n$ obtained from $L=L_0$ by an $n$-move, and $L_{\infty}$}
\ \\

When computing and comparing Goeritz matrices of $L=L_0$, $L_n$ and 
$L_{\infty}$ we can assume that black regions are chosen as in Figure 2.10 
and that the white region $X$ in $R^2-L_{\infty}$ is divided into two 
regions $X_0$ and $X_1$ in $R^2-L$.

\begin{lemma}\label{Lemma 1.12}
$G(L_n)=\left[
  \begin{array}{cc}
  G(L_{\infty}) & \alpha \\
   \alpha^T & q+n
   \end{array}
   \right],$

\end{lemma} 

\begin{corollary}\label{Corollary IV 1.13}
$$(i)\ \ \ Det G(L_n) - Det G(L_0) = n Det G(L_{\infty}),$$
$$(ii) \ \ \ \sigma(G(L_0)) \leq \sigma(G(L_n)) \leq \sigma(G(L_0)) +2,\ 
n\geq 0. $$
(iii)\ \  $|\sigma(G(L_n)) - \sigma(G(L_{\infty}))|\leq 1$.
Furthermore, $\sigma(G(L_n)) = \sigma(G(L_{\infty}))$ if and only if 
$rank G(L_n) = rank G(L_{\infty})$ or $rank G(L_n) = rank G(L_{\infty})+2$.
\end{corollary} 

If we orient $L=L_0$ we can use Corollary 2.13(ii) to obtain 
very useful properties of signature of $L$ and $L_n$.

\begin{corollary}[\cite{P-2}]\label{Corollary IV 1.14}
\begin{enumerate}
\item[(i)] Assume that $L_0$ is oriented in such a way that its 
strings are parallel. $L_n$ is said to be obtained from $L_0$ by 
a $t_n$-move 
(\parbox{3.5cm}{\psfig{figure=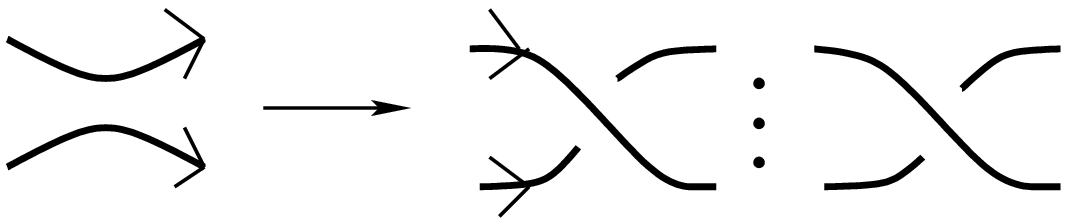,height=0.7cm}}); then
$$n-2 \leq \sigma(L_0)) - \sigma(L_n) \leq n$$
\item[(ii)] Assume that $L_0$ is oriented in such a way that its
strings are anti-parallel and that $n=2k$ is an even number.
$L_{2k}$ is said to be obtained from $L_0$ by a $\bar t_{2k}$-move 
(\parbox{4.8cm}{\psfig{figure=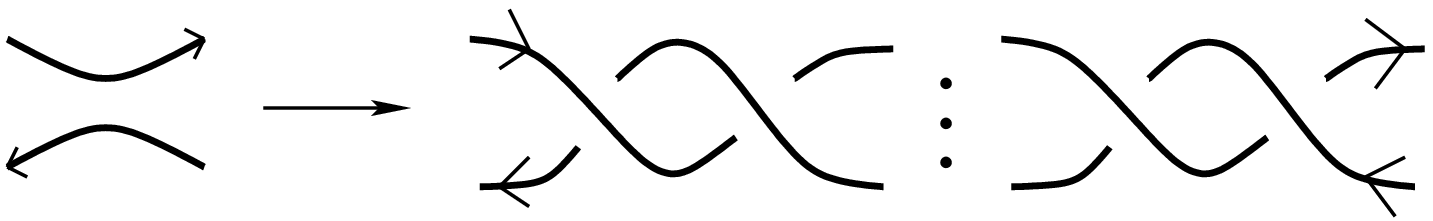,height=0.7cm}}); then 
$$ 0 \leq \sigma(L_{2k}) - \sigma(L_{0}) \leq 2.$$
\item[(iii)] (Giller \cite{Gi})  
$$ 0 \leq \sigma(
L_{\parbox{0.5cm}{\psfig{figure=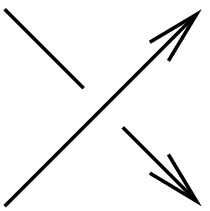,height=0.5cm}}})
 - \sigma(L_{\parbox{0.5cm}{\psfig{figure=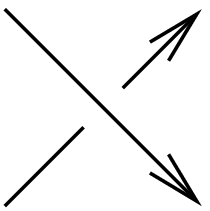,height=0.5cm}}}) \leq 2 $$
\end{enumerate}
\end{corollary} 

\begin{proof} (i) All new crossings of $L_n$ are of type II (we use shading 
of Figure 2.10), thus $\mu(L_n) - \mu(L_0)=n$. Therefore by Corollary 2.13(ii) 
we have $n-2 \leq \sigma(G_{L_0})-\mu(L_0) - (\sigma(G_{L_n})-\mu(L_n) \leq n$, 
and Corollary 2.14(i) follows by Corollary 2.4.\\
(ii) In this case $\mu(L_{2k}) = \mu(L_0)$ thus (ii) follows from 
Corollary 2.13(ii). The generalization of Corollary 2.14(ii) to 
Tristram-Levine signatures is given in Corollary 
\ref{Corollary H-6.9}(ii). \\
 (iii) follows from (i), or (ii) for $n = 2$.
\end{proof}

 We finish the section with an example of computing a close form for
the determinant of the family of links called Turk-head links.
We define the $n$th Turk-head link, $Th_n$ as the closure of the 3-braid
$(\sigma_1\sigma_2^{-1})^n$ (see Figure 2.11 for $Th_6$).\footnote{$Th_0$ is the trivial link of 3 components, 
$Th_1$ the trivial knot, $Th_2$ the figure eight knot ($4_1$), $Th_3$ the Borromean rings ($6^3_2$), 
$Th_4$, the knot $8_{18}$, $Th_5$ the knot $10_{123}$, $Th_6$ the link $12^3_{474}$ (that is 474th link 
of $12$ crossings and $3$ components in unpublished M.~Thistlethwaite tables; compare \cite{This-1}), and 
 $Th_7$ and $Th_8$ are the knots $14_{a19470}$ and $16_{a275159}$, respectively, 
in Thistlethwaite (Knotscape) list.}.

\begin{example}\label{Example IV.1.15}
We compute that
$$Det_{Th_n} = (\frac{3+\sqrt 5}{2})^n +  (\frac{3-\sqrt 5}{2})^n -2,$$
or it can be written as
$Det_{Th_n} = T_n(3) -2$, where $T_i(z)$ is the Chebyshev
(Tchebycheff) polynomial of the first kind\footnote{$T_n(3)$ is often named the Lucas 
number; more precisely $T_n(3)= {\L}_{2n}$, where ${\L}_0=2$, ${\L}_1=1$ and ${\L}_n={\L}_{n-1}+{\L}_{n-2}$ 
as ${\L}_n= 3{\L}_{n-2} - {\L}_{n-4}$.}:
$$ T_0=2,\ \ T_1=z, \ \ T_i=zT_{i-1} - T_{i-2}.$$
In particular, $Det_{Th_2} = 5$, $Det_{Th_3} = 16$, $Det_{Th_4} = 45$, $Det_{Th_5} = 121$, 
$Det_{Th_6} = 320$, $Det_{Th_7} = 841$, and $Det_{Th_8} = 2205$; compare \cite{Sed,Mye}.
\end{example}
To show the above formulas, consider the (unreduced) Goeritz matrix related
to the checkerboard coloring of the diagram of $Th_n$ as shown in Figure 2.11
(we have here $z = 3$ and we draw the case of $n = 6$).
$$G'(Th_{6})=\left[
  \begin{array}{ccccccc}
 -n  &  1  &  1 &  1 &  1  &  1  &  1 \\
  1  & -z  &  1 &  0 &  0  &  0  &  1 \\
  1  &  1  & -z &  1 &  0  &  0  &  0 \\
  1  &  0  &  1 & -z &  1  &  0  &  0 \\
  1  &  0  &  0 &  1 & -z  &  1  &  0 \\
  1  &  0  &  0 &  0 &  1  & -z  &  1 \\
  1  &  1  &  0 &  0 &  0  &  1  & -z 
   \end{array}
   \right],$$

    By crossing the first row and column of $G'(Th_n)$ we obtain the Goeritz
matrix of $Th_n$ which is also the circulant matrix with the first row
$(-z, 1, 0, ..., 0, 1)$ ($z = 3$ and $n = 6$ in our concrete case):

$$G(Th_{6})=\left[
  \begin{array}{cccccc}
   -z  &  1 &  0 &  0  &  0  &  1 \\
    1  & -z &  1 &  0  &  0  &  0 \\
    0  &  1 & -z &  1  &  0  &  0 \\
    0  &  0 &  1 & -z  &  1  &  0 \\
    0  &  0 &  0 &  1  & -z  &  1 \\
    1  &  0 &  0 &  0  &  1  & -z
   \end{array}
   \right],$$

To compute the determinant of the circulant matrix $CM_n(z)$ of the size
$n\times n$ and the first row $(-z, 1, 0, ..., 0, 1)$ we treat each row
as a relation and find the structure of the $Z[z]$ module generated
by columns (indexed by ($e_0 , e_1 , ..., e_n$)).
Thus we have $n$ relations of the form $e_k = ze_{k-1} - e_{k-2}$,
where $k$ is taken modulo $n$. The relation recalls the relation
of Chebyshev polynomials, and
in fact we easily check that $e_k = S_{k-1}(z)e_1 - S_{k-2}(z)e_0$,
where $S_k(z)$ is the Chebyshev polynomial of the second kind:
$$S_0=1,\ \ S_1=z, \ \ S_i=zS_{i-1} - S_{i-2}.$$
Thus we can eliminate all vectors (columns) $e_k$ except $e_0$ and $e_1$,
 and we are left with two equations $e_0 = e_n = S_{n-1}e_1 - S_{n-2}e_0$ ,
and $e_1 = e_{n+1} = S_{n}e_1 - S_{n-1}e_0$.
Thus, our module can be represented by the $2\times 2$ matrix
$$
\left[
  \begin{array}{cc}
S_{n-1} & 1-S_n \\
S_{n-2} +1 & -S_{n-1} 
   \end{array}
   \right],$$

We conclude that, up to a sign, $det CM_n(z)$ is equal to the determinant
of our $2\times 2$ matrix, that is
$S_n - S_{n-2} - 1 - S^2_{n−1} + S_n S_{n-2}$. To simplify this
expression let us use the substitution $z=a+a^{-1}$. Then
$S_n(z) = a^n + a^{n-2}+ \ldots  a^{2-n}+  a^{-n}=
\frac{a^{n+1}- a^{-n-1}}{a- a^{-1}}$, and $T_n(z) = a^n + a^{-n}$.
Therefore, $S_n - S_{n-2} - 1 - S^2_{n-1} + S_n S_{n-2}=
S_n- S_{n-2} - 1 - ((\frac{a^{n}- a^{-n}}{a- a^{-1}})^2-
(\frac{a^{n+1}- a^{-n-1}}{a- a^{-1}})(\frac{a^{n-1}- a^{-n+1}}{a- a^{-1}}))=$

$S_n- S_{n-2} - 1 - (\frac{(a^{n}- a^{-n})^2 - (a^{n+1}- a^{-n-1})
(a^{n-1}- a^{-n+1})}{(a- a^{-1})^2})= S_n - S_{n-2} - 2 =
a^n + a^{-n} -2 = T_n(z)-2$.\\
By comparing the maximal power of $z$ in $det CM_n(z)$ and $T_2(z)-2$, we
get that $det CM_n(z) = (-1)^n(T_n(z)-2)$. For $z=3$ we have $a+a^{-1}=3$,
thus $a= \frac{3\pm \sqrt 5}{2}$ so we can choose $a=\frac{3+ \sqrt 5}{2}$
and $a^{-1} = \frac{3 - \sqrt 5}{2}$, and thus
$T_n(3)=  (\frac{3+\sqrt 5}{2})^n +  (\frac{3-\sqrt 5}{2})^n$.

Because, $Th_n$ is
an amphicheiral link, its signature is equal to $0$ and
$$Det_{Th_n}= i^{\sigma(Th_n)}|det CM_n(3)|= T_n(3)-2 =
 (\frac{3+\sqrt 5}{2})^n +  (\frac{3-\sqrt 5}{2})^n -2.$$

\begin{center}
\psfig{figure=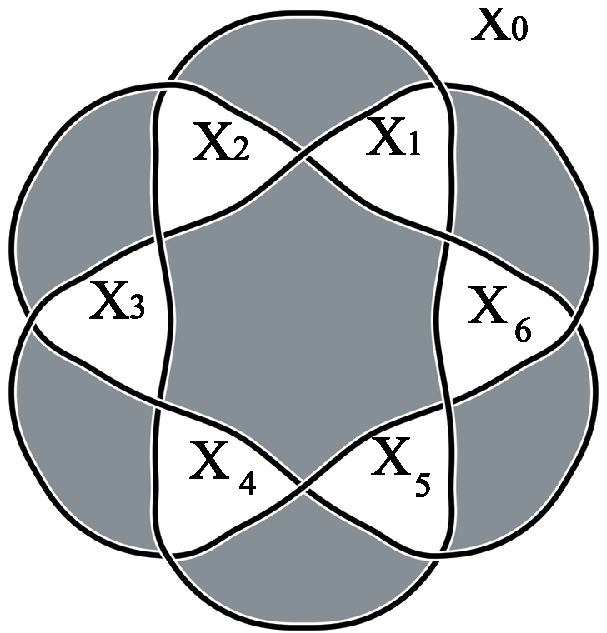,height=4.1cm}
\end{center}

\centerline{Figure 2.11; The Turk-head link $Th_6$ and its checkerboard coloring }
\ \\
We computed the determinant of the circulant\footnote{Recall, that the
circulant $n\times n$ matrix, satisfies
$a_{i,j}= a_{i-1,j-1}=\ldots a_{1,j-i+1}$, $0\leq i,j\leq n-1$. Such a matrix
has (over $C$) $n$ different eigenvectors:
$(1,\omega,\omega^2,...,\omega^{n-1})$,
where $\omega$ is any $n$th root of unity ($\omega^n=1$). The corresponding
eigenvalues are $\lambda_{\omega}=\sum_{i=0}^{n-1}\omega^i a_{1,i}$.
Thus Example \ref{Example IV.1.15} leads to a curious identity
$\Pi_{i=0}^{n-1}(\omega^i +\omega^{-i}-z) = det CM_n = (-1)^n(T_n(z)-2)$ for
any primitive $n$th root of unity $\omega$.}
matrix $CM_n$ for a general
variable $z$ and till now used it only for $z=3$, we see in the next exercise
that the matrix has knot theory interpretation for any rational number $z$.

\begin{exercise}\label{Exercise IV.1.16}
Consider the ``braid like" closure of the tangle $(\sigma_2^{-\frac{1}{a}}
\sigma_1^b)^n$ for any integers $a$ and $b$, (see Figure 2.12 for
$(\sigma_2^{-\frac{1}{3}}\sigma_1^3)^4$). Show that the
determinant of the link satisfies the formula\footnote{It is also 
the formula for the number of spanning trees of the generalized 
wheel, $W_{a,b,n}$, which is the Tait graph of the closure of 
$(\sigma_2^{-\frac{1}{a}} \sigma_1^b)^n$ ($W_{3,3,4}=$
\parbox{1.0cm}{\psfig{figure=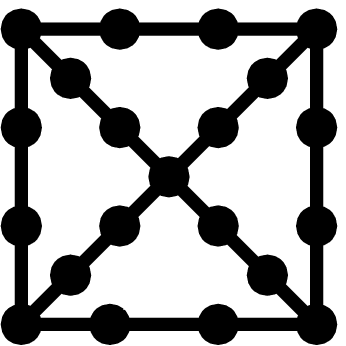,height=0.8cm}};
compare Chapter V of \cite{P-Book}).}
$$|Det_{(\sigma_2^{-\frac{1}{a}}\sigma_1^b)^n}|= |b^n det CM_n(2+
\frac{a}{b})|= |b^n (T_n(2+\frac{a}{b})-2)|.$$
\end{exercise}

\begin{center}
\psfig{figure=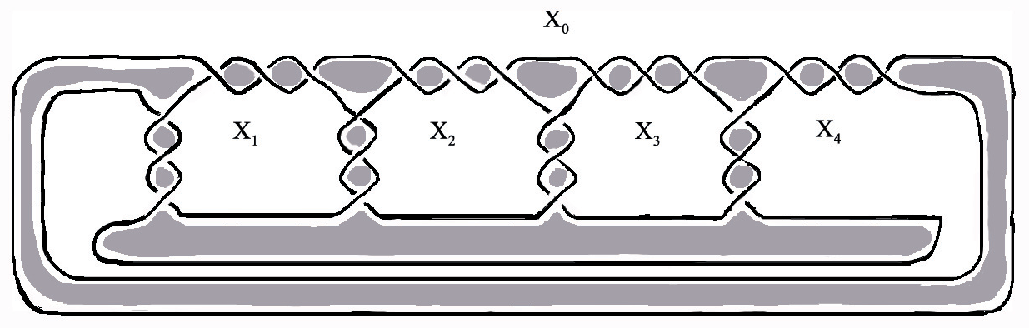,height=3.6cm}
\end{center}

\centerline{Fig. 2.12; The closure of the tangle
$(\sigma_2^{-\frac{1}{3}}\sigma_1^3)^4$ }

\section{Seifert surfaces}\label{Section H-3}

It was first demonstrated by P.~Frankl and L.~Pontrjagin in 1930 \cite{F-P} 
that any knot bounds an oriented surface\footnote{According to \cite{F-P}: 
``The Theorem... [was] found by both authors independently from each other. 
In what follows, the Frankl's form of the proof is presented." One should add 
that Seifert refers in \cite{Se} to the Frankl-Pontrjagin paper 
and says that they use a different method.}.
  H.~Seifert found a very simple 
construction of such a surface \cite{Se}
 and developed several applications of the 
surface, named now {\it Seifert surface} (also, infrequently, Frankl-Pontrjagin 
surface)\footnote{Kauffman in \cite{K-3,K-8} uses the term  
{\it Seifert surface} to describes the 
surface obtained from an oriented link diagram by the Seifert algorithm 
(Construction \ref{c3:2.4}), and the term spanning surface for an oriented 
surface bounding a link (our Seifert surface of Definition \ref{c3:2.1}). In 
\cite{Bol} the name {\it Frankl surface} is used for any, oriented 
or unoriented spanning surface.}.
\begin{definition}\label{c3:2.1}
A Seifert surface of a link $L\subset S^3$ is a compact,
connected, orientable 2-manifold 
$S\subset S^3$ such that $\partial S = L$.
\end{definition}
For example: a Seifert surface of a trefoil knot is pictured in Fig.~3.1.
If the link $L$ is oriented 
then its Seifert surface $S$ is assumed to be oriented
so that its orientation agrees with that of $L$.

\begin{center}
\begin{tabular}{c} 
\includegraphics[trim=0mm 0mm 0mm 0mm, width=.3\linewidth]
{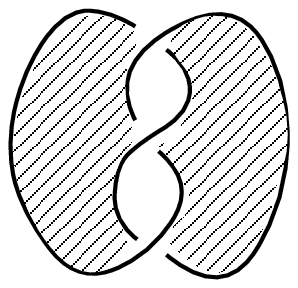}\\
\end{tabular}
\\
Fig. 3.1\\
\end{center}
%

\begin{definition}\label{c3:2.2}
The genus of a link $L\subset S^3$ is the minimal
genus of a Seifert surface of $L$.
\end{definition}

The genus is an invariant (of ambient isotopy classes) of
knots and links. The following theorem provides that it is well
defined.

\begin{theorem}(Frankl-Pontrjagin-Seifert)\label{c3:2.3}
Every link in $S^3$ bounds a Seifert surface.
If, moreover, the link is oriented then there exists a Seifert
surface, an orientation of which determines the orientation
of its boundary coinciding with that of $L$.
\end{theorem}

\begin{construction}(Seifert)\label{c3:2.4}
Consider a fixed diagram $D$ of an an oriented link $L$ in $S^3$. 
In the diagram there are two types of crossings, in a neighborhood
of each of crossings we make a modification of the link (called smoothing)
according to Fig.~3.2.

\begin{center}
\begin{tabular}{c} 
\includegraphics[trim=0mm 0mm 0mm 0mm, width=.5\linewidth]
{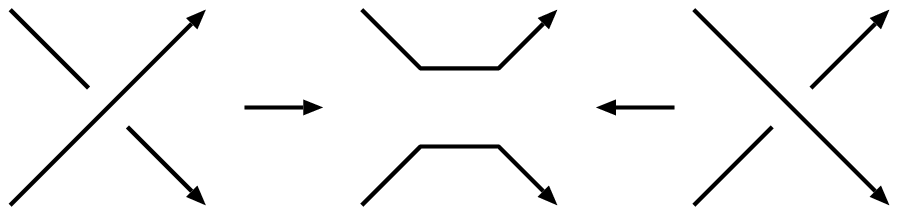}\\
\end{tabular}
\\
Fig. 3.2\\
\end{center}
%

After smoothing all crossings of $D$ we obtain a family of disjoint 
oriented simple closed curves in the plane, called by R.~Fox, Seifert circles and 
denoted by $D_{\vec s}$.
Each of the curves of $D_{\vec s}$ bounds a disk in the plane;
the disks do not have to be disjoint (they can be nested).
Now we make the disks disjoint by pushing them slightly up above the 
projection plane. We start with the innermost disks (that is disks without any 
other disks inside) and proceed outwards (i.e. if $D' \subset D$ then 
$D'$ is pushed above $D$; see Fig.~3.3.

\begin{center}
\begin{tabular}{c} 
\includegraphics[trim=0mm 0mm 0mm 0mm, width=.8\linewidth]
{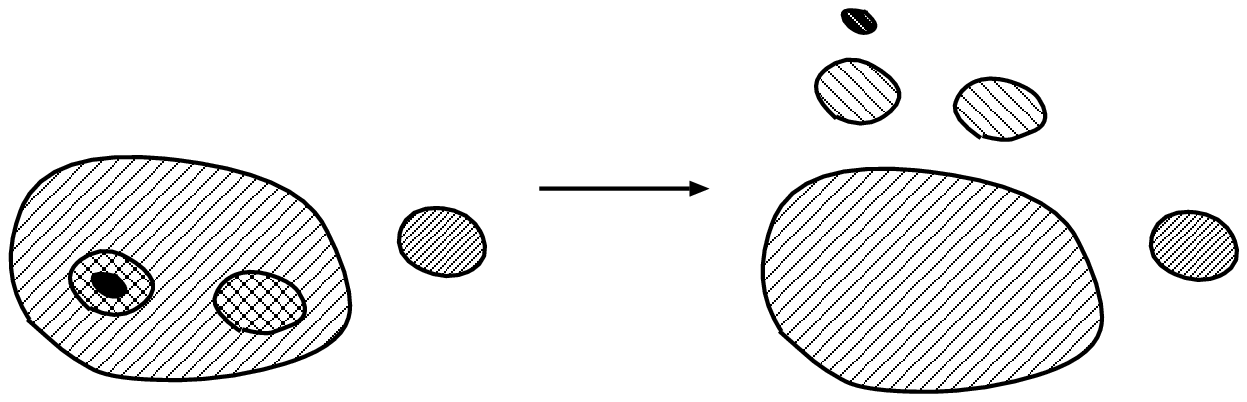}\\
\end{tabular}
\\
Fig. 3.3\\
\end{center}
%

The disks are two-sided so we can assign the sings $+$ and $-$ to each 
of the sides of a disk according to the following convention: 
the sign of the ``upper'' side of the disk is $+$ (respectively, $-$)
if its boundary is oriented  counterclockwise (respectively, clockwise),
see Fig.~3.4.

\begin{center}
\begin{tabular}{c} 
\includegraphics[trim=0mm 0mm 0mm 0mm, width=.5\linewidth]
{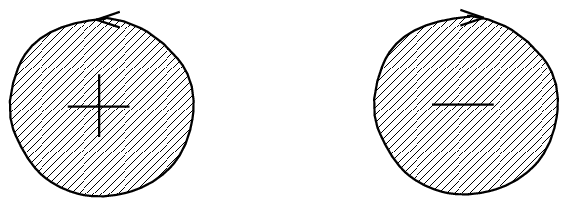}\\
\end{tabular}
\\
Fig. 3.4\\
\end{center}
%

Now we connect the disks together at the original crossings of the diagram $D$
by half-twisted bands so that the 2-manifold which we obtain 
has $L$ as its boundary, see Fig.~3.5.

\begin{center}
\begin{tabular}{c} 
\includegraphics[trim=0mm 0mm 0mm 0mm, width=.8\linewidth]
{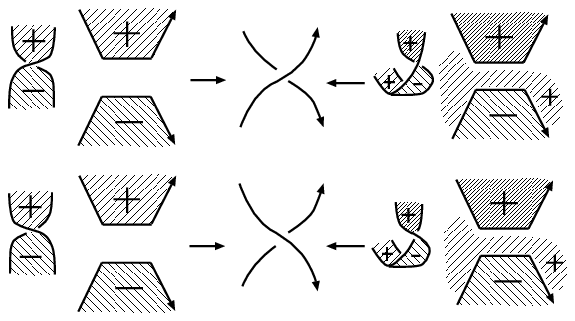}\\
\end{tabular}
\\
Fig. 3.5; Seifert surface around a crossing
\end{center}
%

Since the ``$+$ side'' is connected to another ``$+$ side''
it follows that the resulting surface is orientable.
Moreover, this surface is connected if the projection of 
the link is connected (for example if $L$ is a knot). 
If the surface is not connected then we join its components by tubes
(see Fig.~3.6) in such a way that the orientation of components is preserved.

\begin{center}
\psfig{figure=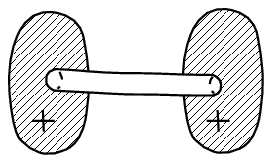,height=4.3cm}
\end{center}

\begin{center}
Fig. 3.6\\
\end{center}
%

\end{construction}

\begin{remark}\label{c3:2.5}
If the link $L$ has more than one component then the Seifert surface,
which we constructed above, depends on the orientation
of components of $L$. This can be seen on the example of a torus link 
of type $(2,4)$,
see Fig.~3.7.

\begin{center}
\psfig{figure=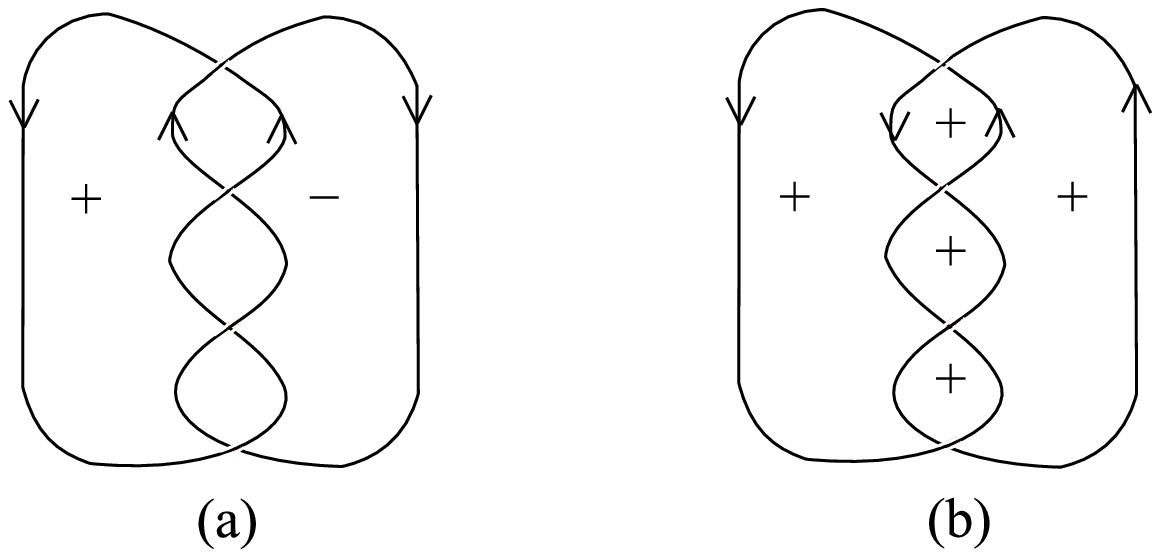,height=4.3cm}
\end{center}
\centerline{Fig. 3.7; 
different orientations result in different Seifert surfaces}

%
\ \\
\ \\
The Seifert surface from Fig.~3.7(a) has genus 1 while 
the surface from Fig.~3.7(b) has genus 0. Therefore the link $L$ 
has genus 0 (as an unoriented link).
\end{remark}

\begin{corollary}
If a projection of a link $L$ is connected (e.g. if~$L$ is a knot)
then the surface, from the Seifert Construction 
\ref{c3:2.4}, is unknotted, that is, its complement in $S^3$ 
is a handlebody. The genus of the handlebody is equal to $c+1-s$ and 
the Euler characteristic is equal to $s-c$, 
where $c$ denotes the number of crossings 
of the projection and $s$ the number of Seifert circles.
\end{corollary}

\begin{proof} The complement in $S^3$ of the plane projection of $L$  
is a 3-disk with $c+1$ handles  (the projection of $L$ cuts 
the projection plane (or 2-sphere) into $c+2$ regions). 
Furthermore adding $s$ 2-disks in the construction
of the Seifert surface we cut $s$ of the handles thus the result remains
a 3-disk with $c+1-s$ handles. The Euler characteristic of obtained 
handlebody is equal to $1-(c+1-s)= s-c$.
\end{proof}

\begin{corollary}\label{c3:2.7}
A knot $K$ in $S^3$ is trivial if and only if
its genus is equal to $0$.
\end{corollary}

\begin{exercise}\label{Exercise H-3.8}
Let $L$ be a link with $n$ components and $D_L$ its diagram. Moreover, 
let $c$ denote the number of crossings in $D_L$
and let $s$ be the number of Seifert circles. 
Prove that the genus of the resulting Seifert surface
is equal to:
$$ \mbox{genus}(S) = p-\frac{s+n-c}{2},$$ 
where $p$ is the number of connected components of the projection of $L$.\\
Check that the Euler characteristic of $S$, for $p=1$, is equal to 
$s-c$ so it agrees with the Euler characteristic of handlebody described in 
Corollary 3.6. 
\end{exercise}

Suppose that the solid torus $V_K$ is a closure of a regular neighborhood of
a knot $K$ in $S^3$ and set $M_K = S^3 - \mbox{int }V_K$ (note
that $M_K$ is homotopy equivalent to the knot complement).
Let us write Mayer-Vietoris sequence for the pair  $(M_K,V_K)$:
$$ 0=H_2(S^3) \rightarrow H_1(\partial M_K)\rightarrow H_1(M_K)\oplus H_1(V_K)
\rightarrow H_1(S^3)=0.$$
For a torus $\partial M_K$ and the solid torus $V_K$ homology are 
$Z\oplus Z$ and $Z$, respectively. Therefore  $H_1(M_K) = Z$ and 
it is generated by a meridian in
$\partial M_K =\partial V_K$, where by the meridian we understand 
a simple closed curve in $\partial V_k$ which bounds a disk in $V_K$.
We denote the meridian by $m$. 
A simple closed curve on $\partial M_K$ which generates
$\ker(H_1(\partial M_K)\rightarrow H_1(M_K))$ 
is called longitude and it is denoted by $l$. 
If $S^3$ and $K$ are oriented then the longitude is orientated
in agreement with the orientation of $K$. 
Subsequently, the meridian is given  
the orientation in such a way that the pair $(m,l)$
induces on $\partial V_K$ the same orientation as the one induced by 
the solid torus $V_K$, which inherits its orientation from $S^3$. 
Equivalently, the linking number of $m$ and $K$ is equal to $1$ 
(compare Section 5).
Similar reasoning allows us also to conclude:
\begin{proposition} For any link $L$ in $S^3$ the first homology of 
the exterior of $L$ in $S^3$ is freely generated by meridians of 
components of $L$. In particular, $H_1(S^3-L)=Z^{com(L)}$. 
\end{proposition}

We also can use the Mayer-Vietoris sequence to find the homology 
of the exterior the Seifert surface in $S^3$. Let $F_L$ be a Seifert 
surface of a link $L$ and $F'$ its restriction to $M_L= S^3 -int V_L$.
Let $V_{F'}$ be a regular neighborhood of $F'$ in $M_L$. 
Because $F'$ is orientable $V_{F'}$ is a product $F'\times [-1,1]$ with 
$F^+ = F'\times \{1\}$ and $F^-=F'\times \{-1\}$. The boundary, 
$\partial V_{F'}$ is homeomorphic to $F^+$ and $F^-$ glued together  
naturally along their boundary. Now let us apply the Mayer-Vietoris sequence 
to $V_{F'}$ and $S^3- int V_{F'}$. We get:
$$ 0=H_2(S^3) \rightarrow H_1(\partial V_{F'}) \stackrel{(i_1,-i_2)}{\to}
H_1(S^3- int V_{F'})\oplus H_1(V_{F'}) \rightarrow H_1(S^3)=0.$$
where $i_1$ and $i_2$ are induced by embeddings. Clearly, 
$H_1(S^3- int V_{F'})$ is isomorphic to the kernel of $i_2$. 
We can easily identify the elements $x^{+} - x^{-}$ 
as elements of the kernel,  for any $x$ a cycle in $F'$. 
In the case of $L$ being a knot, these elements generate 
the kernel.

\begin{corollary}\label{Corollary H-3.10}
 The homology groups, $H_1(F_L)$ and $H_1(S^3-F_L)$ are 
isomorphic to $Z^{2g+com(L)-1}$, where $g$ is the genus of $F_L$ and $com(L)$ 
is the number of components of $L$ thus also the number of boundary 
components of $F_L$. 
Compare Theorem \ref{Theorem H-3.12}.
\end{corollary}
\begin{corollary}\label{Corollary 3.11}
 Let $x_1,...,x_{2g}$ be a basis of $H_1(F_K)$ where 
$F_K$ is the Seifert surface of a knot $K$. Then $x_1^+-x_1^-,\ldots , 
x_{2g}^+-x_{2g}^-$ form  a basis of $H_1(S^3-K)$.
\end{corollary}

With some effort we can generalize Corollary \ref{Corollary H-3.10}
to get the following result which is
a version of Alexander-Lefschetz duality\footnote{Let us
recall that Alexander duality gives us an isomorphism
$\tilde H^i(S^n - X) = \tilde H_{n-i-1}(X)$ for a compact
subcomplex $X$ of $S^n$
and that on the free parts of homology the Alexander isomorphism induces
a nonsingular form
$\beta: \tilde H_i(S^n - X) \times \tilde H_{n-i-1}(X) \to Z$,
where $\tilde H$ denotes reduced (co)homology.}
(see \cite{Li-12} for an elementary proof).

\begin{theorem}i\label{Theorem H-3.12}
Let $F$ be a Seifert surface of a link, then
$H_1(S^3 - F)$ is isomorphic to $H_1(F)$ and there is a nonsingular
bilinear form $$\beta: H_1(S^3 - F) \times H_1(F) \to Z$$
given by $\beta(a,b) = \lk(a,b)$, where $\lk(a,b)$ is defined to be the intersection 
number of $a$ and a $2$-chain whose boundary is $b$ (see Chapter 5).
\end{theorem}

\section{Connected sum of links.}\label{Section H-4}

\begin{definition}\label{Definition H-4.1}
Assume that $K_1$ and $K_2$ are oriented knots in $S^3$.
A connected sum of knots, $K=K_1\kwad K_2$,
is a knot $K$ in $S^3$ obtained in the following way:

First, for $i=1,\ 2$,
choose a point $x_i\in K_i$ and its regular neighborhood
$C_i$ in the pair $(S^3,K_i)$.
Then, consider a pair 
$((S^3-\mbox{int } C_1\cup_\varphi S^3 - \mbox{int }C_2), 
(K_1-\mbox{int }C_1\cup_\varphi K_2 - \mbox{int }C_2))$, 
where $\varphi$ is an orientation reversing homeomorphism
$\partial C_1\rightarrow\partial C_2$ 
which maps the end of $K_1\cap (S^3-\mbox{int }C_1)$ 
to the beginning of $K_2\cap(S^3-\mbox{int }C_2)$
(and vice versa). (Notice that notions of beginning
and end are well defined because $K_1$ and $K_2$ are oriented.)
We see that $(S^3-\mbox{int }C_1)\cup_\varphi (S^3-\mbox{int }C_2)$ 
is a 3-dimensional 
sphere and $(K_1-\mbox{int }C_1)\cup_\varphi (K_2-\mbox{int }C_2)$ 
is an oriented knot.
\end{definition}

\begin{lemma}\label{Lemma H-4.2}
The connected sum of knots is a well defined, commutative and
associative operation in the category of oriented knots
in $S^3$ (up to ambient isotopy).
\end{lemma}

A proof of the lemma follows from two theorems in PL topology
 which we quote without proofs.

\begin{theorem}\label{Theorem H-4.3}
Let $(C, I)$ be a pair consisting of a
3-cell $C$ and 1-cell $I$ which is properly embedded and unknotted
in $C$ (i.e. the pair $(C, I)$ is homeomorphic 
to $(\bar B(0,1),[-1,1])$ where $(\bar B(0,1)$ is the closed unit ball 
in $R^3$ and $[-1,1]$ is the interval $(x,0,0)$ parameterized by $x\in [-1,1]$. 
Respectively, let $(S^2,S^0)$ be a pair consisting
of the 2-dimensional sphere and two points on it.
Then any orientation preserving homeomorphism of $C$ (respectively, $S^2$)
which preserves $I$ and is constant on $\partial I$ (respectively,
it is constant on $S^0$) is isotopic to the identity.
\end{theorem}

\begin{theorem}\label{Theorem H-4.4}
Let $K$ be a knot in $S^3$ and let $C'$ and $C''$ be two regular
neighborhoods in the pair $(S^3,K)$ of two points on $K$.
Then there exists an isotopy $F$ of the pair $(S^3,K)$
which is constant outside of a regular neighborhood of $K$
and such that  $F_0=\mbox{ Id}$ and $F_1(C') = C''$. 
\end{theorem}

\begin{remark}\label{Remark H-4.5}
In the definition of the connected sum of knots we assumed that
the homeomorphism $\varphi$ reverses the orientation. 
This assumption is significant as the following example shows.

Let us consider the right-handed trefoil knot $K$
(i.e.~a torus knot of type $(2,3)$), Fig.~3.1. 
Let $\overline{K}$ be the mirror image of $K$ (i.e.~a torus knot of
type $(2,-3)$). Then $K\kwad\overline{K}$ is the square knot
while  $K\kwad K$ is the knot ``Granny'' and these two knots are not
equivalent. To distinguish them it is enough to compute 
their signature or the Jones polynomial, or Homflypt (Jones-Conway) polynomial.
\end{remark}

\begin{theorem}(Schubert \cite{Sch-2})\label{Theorem H-4.6}
Genus of knots in $S^3$ is additive, that is
$$g(K_1\kwad K_2) = g(K_1) + g(K_2).$$
\end{theorem}

A proof of the Schubert  theorem can be found in e.g.~\cite{J-P,Li-12}.

\begin{corollary}\label{Corollary H-4.7}
Any knot in $S^3$ admits a decomposition into a finite
connected sum of prime knots, i.e.~knots which are not 
connected sums of non-trivial knots.
\end{corollary}

In fact Schubert \cite{Sch-1} showed that the prime decomposition of knots 
is unique up to order of factors; in other worlds, knots with connected 
sum form a unique factorization commutative monoid.

\begin{corollary}\label{Corollary H-4.8}
The trefoil knot is non-trivial and prime.
\end{corollary}

Proof. The trefoil knot is non-trivial therefore its genus is positive 
(Corollary \ref{c3:2.7}). Fig.~3.1 demonstrates that the genus is equal to 1.
Now primeness follows 
from Theorem \ref{Theorem H-4.6} and Corollary \ref{Corollary H-4.7}.

Similarly as for knots, 
the notion of connected sum can be extended to oriented
links. It, however, depends on the choice of the components which
we glue together. The weak version of the unique factorization of links 
with respect to connected sum was proven by Youko Hashizume \cite{Hash}.

\section{Linking number; Seifert forms and matrices.}\label{Section H-5}

We start this Section by introducing the linking number $\lk (J,K)$ 
for any pair of disjoint oriented knots $J$ and $K$. Our definition 
is topological and we will show that it agrees with the diagrammatic 
definition considered before.  We use the notation introduced
right after the Exercise \ref{Exercise H-3.8}.

\begin{definition}\label{c3:4.1}
The linking number $\lk (J,K)$ is an integer such that
$[J]=\lk(J,K) [m]$, where $[J]$ and $[m]$ 
are homology classes of 
the oriented curve $J$ and the meridian $m$
of the oriented knot $K$, respectively.
\end{definition}

\begin{lemma}\label{c3:4.2}
Let $S\subset S^3-\mbox{int }V_K$ be a Seifert surface
of a knot $K$ (more precisely, its restriction
to $S^3-\mbox{int }V_K$), such that its orientation determines
the orientation of $\partial S$ compatible with that of the longitude 
$l$. Then $\lk (J,K)$ is equal to the algebraic intersection number
of $J$ and $S$. 
\end{lemma}

Proof. First, let us recall the convention for the orientation 
of the boundary of an oriented manifold $M$.
For $x\in \partial M$ we consider a basis $(v_2,\ldots,v_n)$ of the tangent 
space 
$T_x\partial M$ together with the normal $\overline{n}$ of $\partial M$ in $M$
which is directed outwards. Then,
$v_2,\ldots,v_n$ defines orientation of $T_x\partial M$
if $\overline{n}, v_2,\ldots,v_n$ defines the orientation of $T_x M$.
Returning to the proof of \ref{c3:4.2} we note that the meridian
$m$ intersects the Seifert surface $S$ exactly at one point.
Moreover, the algebraic intersection number of $m$ and $S$ is $+1$, according
to our definition of the orientation of $S$.
Thus, if the algebraic intersection number of $J$ and $S$ is equal to $i$,
then $[J]= i[m]$, that is $i =\lk (J,M)$, which concludes the proof.

\begin{lemma}\label{c3:4.3}
Let us consider a diagram of a link $J\cup K$ consisting of two
disjoint oriented knots $J$ and $K$. 
We assume that the orientation of $S^3 = R^3\cup\infty$ is induced by the 
orientation of the plane containing the diagram of $J\cup K$
and the third axis which is directed upwards
Now, to any crossing of the diagram where $J$ passes under $K$
we assign $+1$ in the case of  
\parbox{1.1cm}{\psfig{figure=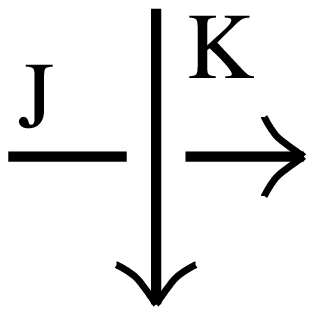,height=0.7cm}} 
 and $-1$ in the case of  
\parbox{1.1cm}{\psfig{figure=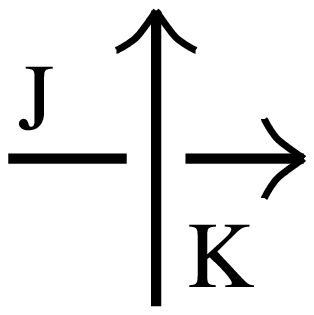,height=0.7cm}}.
Then the sum of all numbers assigned to such crossings is equal to
linking number $\lk (J,K)$.
\end{lemma}

Proof. Let us consider a Seifert surface of the knot $K$
constructed from the diagram of $K$, as described
in Construction \ref{c3:2.4}. We may assume that the knot $J$ is placed
above this surface, except small neighborhoods of the crossings
where $J$ passes under $K$. We check now that the sign of the intersection
of $J$ with this surface coincides with the number that 
we have just assigned to such a crossing.

\begin{exercise}\label{c3:4.4}
Show that $\lk (J,K) = \lk (K,J) = -lk(-K,J)$, where $-K$ denotes
the knot $K$ with reversed orientation.

Hint. Apply Lemma \ref{c3:4.3}.
\end{exercise}

The linking number may be defined for any two disjoint 1-cycles
in $S^3$. For example, as a definition we may use the condition
from Lemma \ref{c3:4.2}. That is, if $\alpha$ and $\beta$ 
are disjoint 1-cycles in $S^3$ then $\lk (\alpha,\beta)$ 
is defined as the intersection number of $\alpha$ with a 2-chain
in $S^3$ whose boundary is equal $\beta$. 

\begin{exercise}\label{c3:4.5}
Prove that $\lk (\alpha, \beta)$ is well defined,
that is, it does not depend on the 2-chain whose boundary is $\beta$.

\end{exercise}

\begin{exercise}\label{c3:4.6}
Show that $\lk$ is symmetric and bilinear, i.e.~$\lk (\alpha,\beta) = \lk(\beta,\alpha)$
and $\lk(\alpha,n\beta)=n\cdot\lk(\alpha,\beta)$, and if a cycle $\beta '$ 
is disjoint from $\alpha$ then 
$\lk (\alpha,\beta+\beta ') = \lk(\alpha,\beta)+lk(\alpha,\beta ')$. 
\end{exercise}

\begin{exercise}\label{c3:4.7}
Prove that, if $\beta$ and $\beta '$ are homologous
in the complement of $\alpha$, then 
$\lk (\alpha, \beta) = \lk(\alpha,\beta ')$. 
\end{exercise}

Now we define a Seifert form of a knot or a link.
Let $S$ be a Seifert surface of a knot or a link $K$, then
$S$ is a two-sided surface in $S^3$.
Let $S\times[-1,1]$ be a regular neighborhood of
$S$ in $S^3$. For a 1-cycle $x$ in $\mbox{ int }S$
we can consider a cycle $x^+$ (respectively $x^-$)
in $S\times\{1\}$ (respectively $S\times\{-1\}$) 
which is obtained by pushing the cycle $x$ to 
$S\times\{1\}$ (respectively, to $S\times\{-1\}$).
(We note that the sides of $S$ are uniquely defined
by the orientations of $K$ and $S^3$.) 

\begin{definition}\label{c3:4.8}
The Seifert form of the knot $K$ is a function 
$f: H_1(\mbox{int }S)\times  H_1(\mbox{int }S)\rightarrow Z$
such that $f(x,y) = \lk(x^+,y)$. Similarly we define
a Seifert form of an oriented link $L$ using an oriented
Seifert surface $S$ of $L$.
\end{definition}

\begin{lemma}\label{c3:4.9} 
The function $f$ is a well defined bilinear form
on the $Z$-module (i.e. abelian group) $H_1(\mbox{int }S)$.
\end{lemma}

Proof. The result follows from Exercises  \ref{c3:4.6} and \ref{c3:4.7}.

\begin{definition}\label{c3:4.10}\ \\
By a Seifert matrix $V=\{v_{i,j}\}$ in a basis 
$e_1,e_2,...,e_{2g+com(L)-1}$ of $H_1(S)$ 
we understand the matrix of $f$ in this basis, that is
$$ v_{i,j}  = \lk (e^+_i,e_j).$$
Then, for $x,y\in H_1(S)$ we have $f(x,y) = x^TVy$.  We use the convention 
that coefficients of a vector are written as a column 
matrix\footnote{Our notation agrees with that of Kauffman \cite{K-3},
 Kawauchi \cite{Kaw-4}, and \cite{J-P} but in the books 
by Burde and Zieschang \cite{B-Z}, Lickorish \cite{Li-12}, 
Rolfsen \cite{Ro-1}, Livingston \cite{Liv}, and Murasugi \cite{M-9} 
the convention is the opposite, that is $f(x,y)= lk(x,y^+)$.}.
\end{definition}

Notice, that a change of the basis in $H_1(S)$ 
results in the change of the matrix $V$ to a similar 
matrix $P^T VP$, where $det P =\pm 1$.

\begin{example}\label{Example H-5.11}
The Seifert matrix of a Seifert surface 
of the right-handed trefoil knot, computed in the basis $[\alpha], [\beta]$,
is equal to 
$\left[\begin{array}{cc}-1&0\\1&-1\\ \end{array}\right]$ (see
Fig.~5.1).

\begin{center}
\begin{tabular}{c} 
\includegraphics[trim=0mm 0mm 0mm 0mm, width=.35\linewidth]
{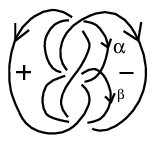}\\
\end{tabular}
\\
Fig. 5.1
\end{center}

%

\end{example}

\begin{example}\label{Example H-5.12}
The Seifert matrix of a Seifert surface of the
figure-eight knot, computed in the basis
$[\alpha], [\beta]$ is equal to
$\left[\begin{array}{cc}1&-1\\0&-1\\ \end{array}\right]$ (see
Fig.~5.2).

\begin{center}
\begin{tabular}{c} 
\includegraphics[trim=0mm 0mm 0mm 0mm, width=.35\linewidth]
{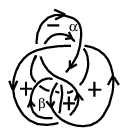}\\
\end{tabular}
\\
Fig. 5.2
\end{center}

%

\end{example}
 With some practice one should be able to find Seifert form efficiently and 
we encourage a reader to compute more examples and develop some rules; 
for example if $\alpha$ is a simple closed curve on $S$ and on the plane then 
$lk(\alpha^+,\alpha)= -\frac{1}{2}\sum \sgn p$ where the sum 
is taken over all crossings of the diagram traversed by $\alpha$.  
We illustrate it by one more example, the Seifert matrix of a pretzel knot. 
The computation is almost the same as in the trefoil case as the genus 
of the surface is equal to $1$ and three crossings of the right-handed trefoil knot $\bar 3_1$ are 
replaced by $2k_1+1$, $2k_2+1$, and $2k_3+1$, respectively.

\begin{example}\label{c3:4.11} Let $P_{n_1,n_2,n_3}$ denote the 
pretzel link of type $(n_1,n_2,n_3)$ (compare Fig. 5.3).
The Seifert matrix of a Seifert surface of the pretzel knot 
$P_{2k_1+1,2k_2+1,2k_3+1}$, computed in the basis $[\alpha], [\beta]$,
is equal to
$\left[\begin{array}{cc}-k_1-k_2& k_2\\k_2+1&-k_1-k_2\\ \end{array}\right]$ 
(see Fig.~5.3).

\begin{center}
\begin{tabular}{c}
\includegraphics[trim=0mm 0mm 0mm 0mm, width=.35\linewidth]
{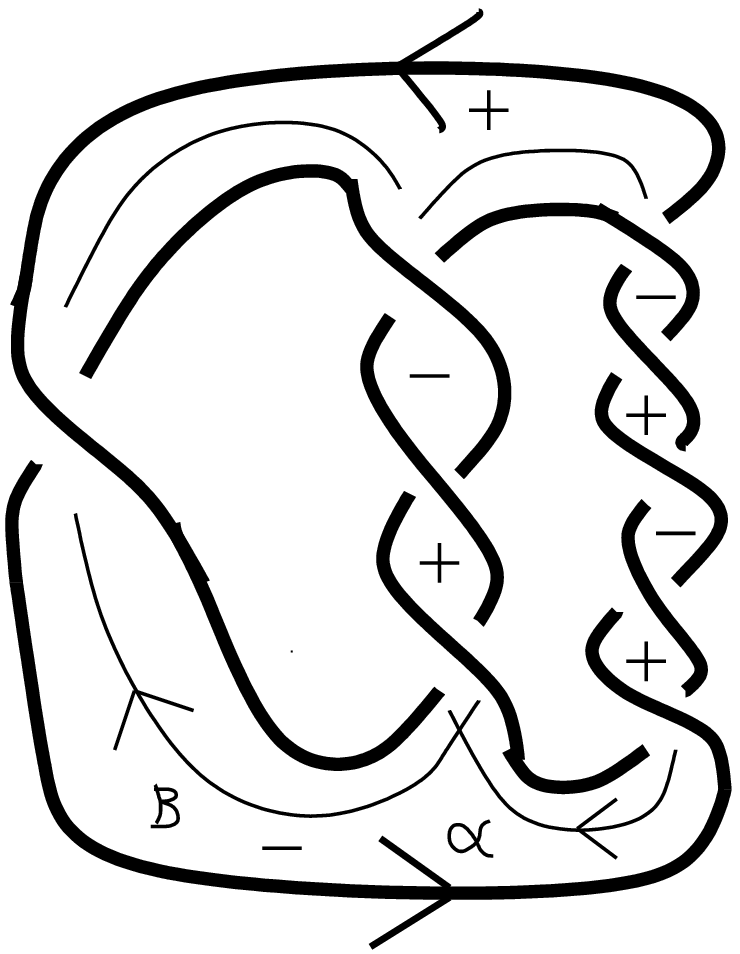}\\
\end{tabular}
\\
Fig. 5.3; $P_{1,3,5}$ -- the pretzel knot of type $(1,3,5)$
\end{center}
\end{example}

There is a classical skew-symmetric form on a homology group of an 
oriented surface, called an intersection form, which is 
related to the Seifert form $f$. 

\begin{definition}\label{c3:4.13}
Let $S$ be an oriented surface. For two homology classes 
$x,y\in H_1(S)$ represented by
transversal cycles we define
their algebraic intersection number $\tau(x,y)$ as the sum of the signed 
intersection points where the sign is defined in the following way:
if $x$ meet $y$ transversally at a point $p$ then the sign 
of the intersection at $p$ is equal $+1$ 
if \includegraphics[trim=0mm 0mm 0mm 0mm, width=.05\linewidth]
{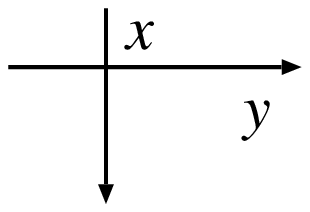} 
and $-1$ if \includegraphics[trim=0mm 0mm 0mm 0mm, width=.05\linewidth]
{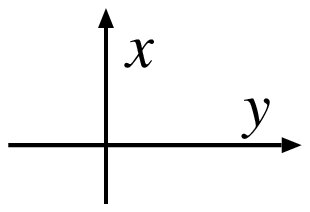}.
\end{definition}

\begin{exercise}\label{c3:4.14}
Prove that $\tau:H_1(S)\times H_1(S)\rightarrow Z$ is bilinear and
skew-symmetric (i.e.~$\tau(x,y) = -\tau(y,x)$).
\end{exercise}

\begin{exercise}\label{c3:4.15}
Prove that the determinant of a matrix of $\tau$ is equal to 1 if 
$\partial S = S^1$, or $\partial S = \emptyset$ and it is equal to $0$ 
otherwise.

Solution. Assume that $S$ has more then one boundary component
 and $\partial_1$ is 
one of them. Then $\partial_1$ is a nontrivial element in $H_1(S)$ with 
trivial intersection number with any element of $H_1(S)$. Thus the matrix 
of $\tau$ is singular and its determinant is equal to zero.  

Assume now that $\partial S = S^1$, or $\partial S = \emptyset$. Let us 
choose loops representing a basis of $H_1(S)$
such as in Fig.~5.4. In this basis the matrix of $\tau$ is as follows
$$\left[
\begin{array}{ccccc}
0&1&&\bigcirc&\\
-1&0&&&\\
&&\ddots&&\\
&&& 0& 1 \\
&\bigcirc&&-1&0\\
\end{array}
\right]
$$ 
Thus, its determinant is equal to 1. Notice also that the determinant 
of a matrix changing a basis of $H_1(S)$ is equal to $1$ or $-1$, thus 
the determinant of the form does not depend on the choice of a basis.

\begin{center}
\begin{tabular}{c} 
\includegraphics[trim=0mm 0mm 0mm 0mm, width=.6\linewidth]
{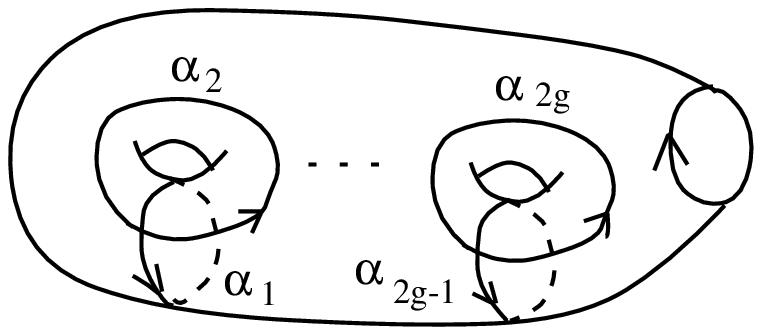}\\
\end{tabular}
\\
Fig. 5.4
\end{center}

%

\end{exercise}

\begin{exercise}\label{IV:4.16}
Prove that, if $S$ is a Seifert surface of a link
then $\tau (x,y) = f(x,y) - f(y,x)$. 

Solution. 
It follows from Lemma \ref{c3:4.3} that the crossing change between 
two oriented disjoint curves $J$ and $K$ in $S^3$ changes the linking 
number between them by $1$ or $-1$, diagrammatically we have:
 $lk(\parbox{1.1cm}{\psfig{figure=FigIV431.eps,height=0.7cm}}) - 
lk(\parbox{1.1cm}{\psfig{figure=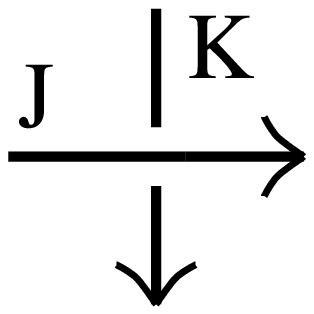,height=0.7cm}}) = 1$. 
If $J$ and $K$ are two, possibly intersecting, oriented curves on an 
oriented surface we see that the pair $(J^+,K)$ differ from the pair 
$(J^-,K)$ by crossing changes at crossings of $J$ and $K$. Furthermore 
the convention we use is that 
$sgn(\parbox{0.9cm}{\psfig{figure=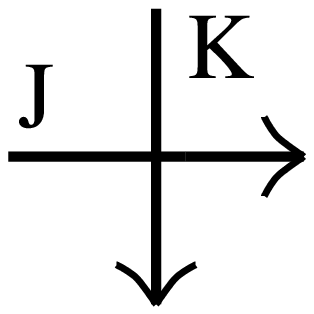,height=0.6cm}}) = -1$.\\ 
Thus $f(J,K)-f(K,J) = \lk(J^+,K)- \lk(K^+,J) = \lk(J^+,K) - \lk(J^-,K)= 
\sum_{p\in J\cap K)}\sgn p = \tau(J,K)$. 
The solution is completed\footnote{The equality  
\parbox{2.8cm}{\psfig{figure=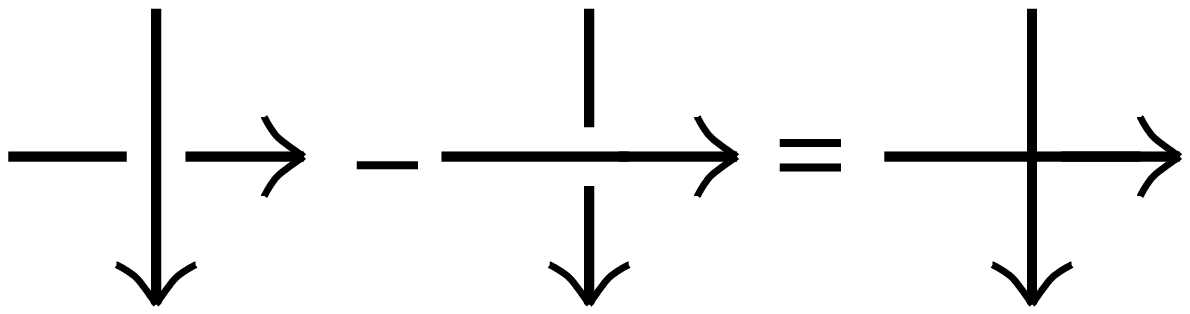,height=0.7cm}} is a defining  
relation of Vassiliev-Gusarov invariants or skein modules of links; compare 
Chapter IX of \cite{P-Book}. This relation, combined with 
\parbox{2.2cm}{\psfig{figure=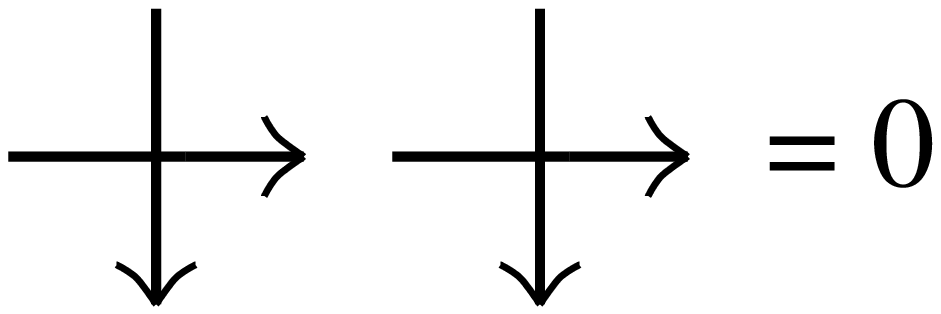,height=0.7cm}} (that is,  
the value of a link with at least two singular crossings is equal to zero), 
leads to the (global) linking number, described as Vassiliev-Gusarov invariant of degree 1.}.
\end{exercise}

\begin{corollary}\label{IV:4.17}
The Seifert matrix $V$ of a knot $K$ in $S^3$ satisfies 
the following equation:
$$\det (V-V^T) = 1.$$
\end{corollary}

Proof. We note that $V-V^T$ is a matrix of $\tau$ (Exercise \ref{IV:4.16})
and its determinant is equal to 1 (Exercise \ref{c3:4.15}).

A Seifert matrix is not an invariant of a knot or a link, but it can be used 
to define some well known invariants, including the Alexander polynomial.


Now we describe relation between Seifert matrices
of (possibly different) Seifert surfaces of a given link.

\begin{definition}\label{Definition H-5.19}
We call matrices $S$-equivalent if one can be obtained from the other
by a finite number of the following modifications: 
\begin{enumerate}

\item [(1)] $A \Leftrightarrow PAP^T$ where $\det P = \mp 1$. 

\item [(2)] $$A \Leftrightarrow \left[ \begin{array}{ccc} A&\alpha&0\\0&0&1\\
0&0&0\\ \end{array}\right] 
\mbox{ and } A\Leftrightarrow \left[
\begin{array}{ccc}
A&0&0\\
\beta&0&0\\
0&1&0\\
\end{array}
\right]
$$
where $\alpha$ is a column and $\beta$ is a row.
\end{enumerate}
\end{definition}

\begin{theorem}\label{c3:4.19}
Let us assume that $L_1$ and $L_2$ are isotopic links and 
$F_1$, respectively,  $F_2$, are their Seifert surfaces. 
If $A_1$ and $A_2$ are their Seifert matrices computed
in some basis  $B_1$ and, respectively, $B_2$ then
$A_1$ is $S$-equivalent to $A_2$.
\end{theorem}

We perform the proof in two steps. Namely, we will prove 
the following two claims:
\begin{enumerate}

\item [(1)] If we attach a handle to $F_1$ then the resulting surface
(boundary of which is again $L_1$) has its Seifert surface $S$-equivalent 
to the Seifert surface of $F_1$.

\item [(2)] We can assume that $L_1=L_2$. Then there exists a Seifert surface for $L_1$ 
which can be reached (modulo isotopy) from both $F_1$ and $F_2$ by the operation of
attaching handles.
\end{enumerate}

First we prove (1).

Let $A_1$ be a Seifert matrix of $F_1$ (in some basis of $H_1(F)$). 
By $\gamma$ and $\mu$ let us denote two new
generators of $H_1(F\cup\mbox{ handle})$ --- see Fig.~5.5.

\begin{center}
\begin{tabular}{c} 
\includegraphics[trim=0mm 0mm 0mm 0mm, width=.5\linewidth]
{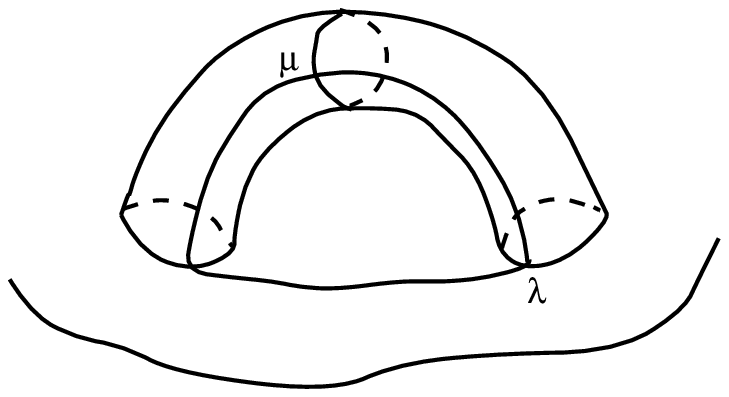}\\
\end{tabular}
\\
Fig. 5.5.
\end{center}

%

Let us recall that the Seifert form 
$f:H_1(F)\times H_1(F)\rightarrow Z$
was defined by the formula
$f (x,y) = \lk(x^+,y)$, 
where $\lk$ denotes linking number in $S^3$ and $x^+$ 
is obtained by pushing the cycle $x$ out of $F$ in the normal direction of
$F$.

If the pushing moves the cycle $\mu$ outside of the handle (that is $\mu^+$ 
is outside the handle) then 
the resulting Seifert matrix is
$$\left[ \begin{array}{ccc}
A&\alpha&0\\
\beta&w_0&0\\
0&\pm 1&0\\
\end{array}
\right]
$$ 
which is $S$-equivalent to the matrix $A$ ($\alpha$ and $\omega_0$ can be 
converted to $0$ matrices by type (1) operations; similarly, 
$\pm 1$ can be converted to $1$ by a type (1) operation). 
In the matrix, $\beta$ is 
a row vector determined by linking numbers of $\lambda^+$ with the basis 
of $H_1(F)$, $\alpha$ is a column vector determined by linking numbers of 
the basis $H_1(F)$ with $\lambda^-$ and $\omega_0 = lk(\lambda^+,\lambda^-)$. 

Otherwise (i.e. $\mu^+$ is inside the handle) we get the matrix:

$$\left[ \begin{array}{ccc}
A& \alpha&0\\
\beta&\omega_0&\pm 1\\
0&0&0\\
\end{array}
\right]
$$ 
which is $S$-equivalent to $A$ as well.

Proof of (2). Assume that  the Seifert surface $F_1$ 
intersects $F_2$ transversally (modulo the boundary $L_1$; in the neighborhood of $L_1$ they may 
 be assumed to be disjoint outside $L_1$). 
Now we will use the following

\begin{lemma}\label{c3:4.20}
Let $M$ be compact connected 3-manifold
and let $F_1$, $F_2$ be such submanifolds of $\partial M$ 
that $\partial M =F_1\cup F_2$ 
and $F_1\cap F_2 = \partial F_1 = \partial F_2$. 
Then there exists a surface $F$ in $M$ such that
$\partial F = \partial F_1$ and $F$ 
can be obtained from $F_1$ as well as from $F_2$
by attaching 1-handles. 
More precisely: $F$ cuts $M$ into two 3-submanifolds $M_1$, containing $F_1$ 
and $M_2$ containing $F_2$. Furthermore $M_i$ can be obtained from $F_i$ 
(more precisely $F_i \times [0,1]$) by attaching $1$-handles\footnote{We 
attach {\it $k$-handle} to an $(n+1)$-dimensional manifold $M$ along 
an open subset, $N$, 
of the boundary by choosing a disk $D^{n+1}= D^k\times D^{n+1-k}$ and 
the embedding $\phi: \partial D^k \times D^{n+1-k} \to N$, and gluing 
$M$ with $D^{n+1}$ using $\phi$. In our case $n=2$.} to 
$int(F_i)$. We have $F_i\cup F = \partial M_i$; in particular $F$ is 
obtained from $F_i$ by 1-surgeries.
\end{lemma}

{\it Sketch of the proof.} The presented proof is based on the proof 
of existence of Heegaard decomposition of a 3-manifold from triangulation 
(e.g. \cite{Hemp,J-P}). 
Let $X$ be a triangulation of $(M,F_1,F_2)$. 
In particular $L$ is in the 1-skeleton of triangulation $\Gamma_1$.
Let $\Gamma_2$ denote the dual
1-skeleton. That is, $\Gamma_2$ is the maximal 1-subcomplex
of the first baricentric subdivision $X'$ of $X$, such that
$\Gamma_2$ is disjoint with $\Gamma_1$.
Let $V_i\  (i = 1,2)$ be a regular neighborhood of
$\Gamma_i$ associated to the second baricentric subdivision of $X$.
Then $X = V_1\cup V_2$ and $V_i$ is obtained from $F_i$ 
by attaching (solid) 1-handles. 
Therefore $(F_1\cup V_1)\cap(F_2\cup V_2)$ 
is the surface $F$ that we look for. 

The proof of claim (2) is inductive with respect to
the number of circles in the intersection $F_1\cap F_2$:

\begin{description}

\item[(1)] Suppose that $F_1\cap F_2 = L_1$. 
Then we apply Lemma \ref{c3:4.20} to a part of $S^3$ which is 
bounded by the closed surface $F_1\cup F_2$.

\item[(n+1)] Inductive step.\ \ 
Suppose that (2) holds if the number of components of $F_1\cap F_2$ is smaller
than $n+1$.\\
Now, assume that $F_1\cap F_2$ consists of $n+1$ circles.
Then $F_1\cup F_2$ cuts $S^3$ into a number of connected
components and moreover different ``sides'' of $F_1$ and $F_2$ 
bound different components. Let $M$ be a component such that
$F'_1 = F_1\cap \partial M$ and $F'_2 = F_2\cap \partial M$.
Now we apply Lemma \ref{c3:4.20} to the triple $(M,F'_1,F'_2)$ and consequently
let $F'_0$ be the surface provided by the lemma. That is, $F'_0$ is obtained
by attaching solid 1-handles to either $F'_1$ or $F'_2$.\\
Let $F^{0}_1$ and $F^{0}_2$ be obtained from $F_1$ and $F_2$ by replacing $F'_1$ and $F'_2$ by $F'_0$.
Then by moving slightly surfaces $F^{0}_1$ and $F^{0}_2$
we can obtain a smaller number of components
of their intersection and thus we can apply the inductive assumption.
This concludes the proof of (2) and of Theorem \ref{c3:4.19}.
\end{description}

An elementary, diagrammatic proof of Theorem \ref{c3:4.19}, based on 
Reidemeister moves and the fact that any link has a special diagram 
(compare Exercises \ref{Exercise 1.9} and \ref {Exercise 1.10} or Proposition 13.15 of \cite{B-Z}), 
is given in \cite{BFK}.

\section{From Seifert form to Alexander polynomial and signatures}\label{IV.5}
\markboth{\hfil{\sc Chapter IV. Goeritz and Seifert matrices }\hfil}
{\hfil{\sc  Alexander polynomial and signature of a link}\hfil}

The Conway's potential function is defined as a 
normalized version of the  Alexander polynomial 
using Seifert matrix, as follows \cite{K-1}:
\begin{lemma}\label{c3:4.21}
Let $A$ be a Seifert matrix of an oriented link $L$ and define 
the potential function $\Omega_L(x) = \det(xA-x^{-1}A^T)$. 
Then $\Omega_L(x)$ does not depend on the choice of a Seifert surface and 
its Seifert matrix. 
In particular, if $T_1$ is the trivial knot then $\Omega_{T_1}(x) =1$.
\end{lemma}

Proof. The result follows from Lemma \ref{c3:4.19}. 
Indeed, simple computations show that if we replace the matrix
$A$ with another $S$-equivalent matrix then $\Omega_L(x)$ remains the same.
We use the following identity
$$\det(x\left[\begin{array}{cc} 0&1\\0&0\\ \end{array}\right] -
x^{-1}\left[\begin{array}{cc} 0&0\\1&0\\ \end{array}\right] ) = \det\left[
\begin{array}{cc}
0&x\\
-x^{-1}&0\\
\end{array}\right] = 1.$$
The same identity is used in the computations for the trivial knot.


If we choose $x=-\sqrt t$ then the potential function is the 
normalized Alexander polynomial (i.e. Alexander-Conway polynomial). 
The transposition of a matrix is preserving its determinant thus the substitution 
$x \to -x^{-1}$ (or $\sqrt t  \to \frac{1}{\sqrt t}$ is preserving the potential function 
and Alexander-Conway polynomial.
Furthermore, we can put $z=x^{-1}-x = \sqrt t -\frac{1}{\sqrt t}$. As follows from Theorem 6.2,
 we obtain, after the substitution, the
Conway polynomial $\nabla_L(z)$ 
(terminology maybe sometimes confused, as 
 $\nabla_L(z)$ is also often called  Alexander-Conway polynomial).
\begin{theorem}[Kauffman \cite{K-1}]\label{Theorem 5.2}\ \\
$\Omega_L(x)= \Delta_L(t) = \nabla_L(z)$, where $x=-\sqrt t$, $z=x^{-1}-x = \sqrt t -\frac{1}{\sqrt t}$.
\end{theorem}
Proof (hint). We have to show that $\Omega_{L_+}(x)- \Omega_{L_-}(x)= 
(x^{-1}-x)\Omega_{L_0}(x)$. In order to demonstrate it we use the properly 
chosen Seifert surfaces $F_+,F_-,F_0$ for $L_+,L_-$ and $L_0$ respectively.

We give all details in the analysis of the more general case of
 the behavior of Seifert matrices under $\bar t_{2k}$-moves,
which generalize the crossing change, which is $\bar t_{\pm 2}$-move.

\begin{definition}(\cite{P-2})\label{Definition 4.23}
 The $\bar t_{2k}$-move  (introducing $k$ full twists on 
anti-parallel oriented arcs) is the elementary operation on an oriented 
diagram $L$ resulting in $\bar t_{2k}(L)$ as illustrated in Figure 6.1.
\end{definition} 
Notice that $\bar t_2$-move is a crossing change from a positive to negative 
crossing ($L_-= \bar t_2(L_+)$). We can choose Seifert surfaces $F(L)$, 
$F(\bar t_{2k}(L))$, and $F(L_{\infty})$ for 
$L=L_{\psfig{figure=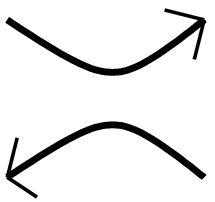,height=0.4cm}}$,
  $\bar t_{2k}(L)$, and 
$L_{\infty}= L_{\psfig{figure=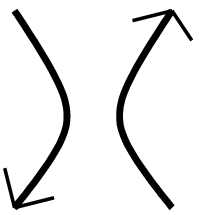,height=0.4cm}}$, 
respectively, to look locally as in Figure 6.1.

\vspace*{0.3in} \centerline{\psfig{figure=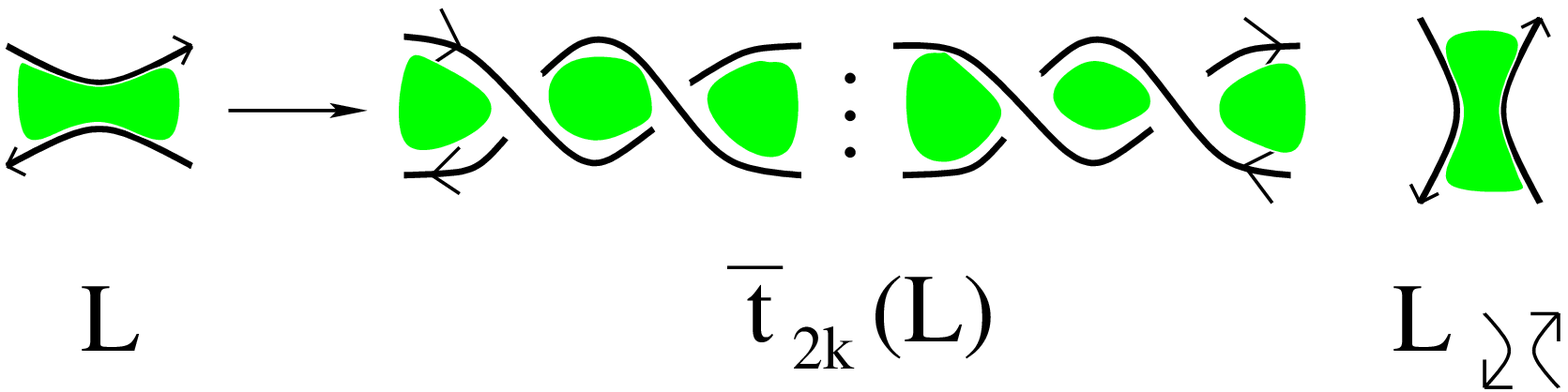,height=3.6cm}}
\begin{center}
Fig. 6.1; Oriented links $L$, $\bar t_{2k}(L)$, and $L_{\infty}$, and their Seifert surfaces  
\end{center}

Let us choose a basis for 
$H_1(F(L_{\psfig{figure=infty-antiparal.eps,height=0.4cm}}))$ and add 
one, standard, element, $e_{\psfig{figure=antiparal.eps,height=0.4cm}}$
 to obtain a basis for $H_1(F(L_{\psfig{figure=antiparal.eps,height=0.4cm}}))$, 
and $e_{\bar t_{2k}(L)}$ to 
get a basis of $H_1(F(\bar t_{2k}(L))$. Denote the 
Seifert matrix of $L_{\psfig{figure=infty-antiparal.eps,height=0.4cm}}$ in 
the chosen basis by 
$A_{L_{\psfig{figure=infty-antiparal.eps,height=0.4cm}}}$. In these bases 
we have immediately.

\begin{lemma}\label{Lemma H-6.4}\ \\

$A_{L_{\psfig{figure=antiparal.eps,height=0.4cm}}}=\left[
  \begin{array}{cc}
  A_{L_{\psfig{figure=infty-antiparal.eps,height=0.4cm}
   }} & \alpha \\
   \beta & q
   \end{array}
   \right],$
\ \\ \ \\

$A_{\bar t_{2k}(L})=\left[
  \begin{array}{cc}
  A_{L_{\psfig{figure=infty-antiparal.eps,height=0.4cm}}} & \alpha \\
   \beta & q+k
   \end{array}
   \right],$

where $\alpha$ is a column given by linking numbers of 
$e^+_{\psfig{figure=antiparal.eps,height=0.4cm}}$ (or 
$e^+_{\bar t_{2k}(L)}$) with basis elements of 
$H_1(F(L_{\psfig{figure=infty-antiparal.eps,height=0.4cm}}))$,
$\beta$ is a row given by linking numbers of basis elements of
$H_1(F(L_{\psfig{figure=infty-antiparal.eps,height=0.4cm}}))$ with 
$e^-_{\psfig{figure=antiparal.eps,height=0.4cm}}$ (or
$e^-_{\bar t_{2k}(L)}$), and $q$ is a number equal to 
$lk(e^+_{\psfig{figure=antiparal.eps,height=0.4cm}}, 
e_{\psfig{figure=antiparal.eps,height=0.4cm}})$ 
(compare \cite{K-1,P-T-2} or \cite{P-1}).
\end{lemma}

\begin{corollary}\label{Corollary H-6.5}
\begin{enumerate}
\item[(i)] If two oriented links are $\bar t_{2k}$ equivalent (that is 
they differ by a finite number of $\bar t_{2k}$-moves) then their 
Seifert matrices are $S$-equivalent modulo $k$.
\item[(ii)] The potential function satisfies:
$$\Omega_{\bar t_{2k}(L)} - 
\Omega_{\psfig{figure=antiparal.eps,height=0.4cm}} = k(x-x^{-1})
\Omega_{\psfig{figure=infty-antiparal.eps,height=0.4cm}}.$$
In particular the case $k=-1$ gives: 
$\Omega_{L_+}(x)- \Omega_{L_-}(x)=
(x^{-1}-x)\Omega_{L_0}(x)$.

\end{enumerate}
\end{corollary}
\begin{proof} (i) It follows from the fact we noted in Lemma \ref{Lemma H-6.4} 
that for properly chosen Seifert surfaces and basis of their homology, 
the entries of Seifert matrices for $\bar t_{2k}$ and 
$\psfig{figure=antiparal.eps,height=0.4cm}$ are congruent modulo $k$.\\
$$(ii) \ \ \Omega_{\bar t_{2k}(L)}= 
\det (xA_{\bar t_{2k}(L)} - x^{-1}A^T_{\bar t_{2k}(L)}) 
= \det \left[
  \begin{array}{cc}
  A_{L_{\psfig{figure=infty-antiparal.eps,height=0.4cm}}} & 
x\alpha - x^{-1}\beta^T \\
   x\beta - x^{-1}\alpha^T& (x-x^{-1})(q+k)
   \end{array}
   \right],$$ 
$$ and \ \ \Omega_{\psfig{figure=antiparal.eps,height=0.4cm}}=
\det \left[
  \begin{array}{cc}
  A_{L_{\psfig{figure=infty-antiparal.eps,height=0.4cm}}} &
x\alpha - x^{-1}\beta^T \\
   x\beta - x^{-1}\alpha^T& (x-x^{-1})q
   \end{array}
   \right].$$
Thus the difference is equal to 
$k(x-x^{-1})\Omega_{\psfig{figure=antiparal.eps,height=0.4cm}}$.

\end{proof}

\begin{example}
We can use Corollary 6.5 to compute the potential (and Alexander-Conway) polynomial of 
the pretzel link $L=P_{2k_1+1,2k_2+1,...,2k_m+1}$ (see Fig. 5.3 or 6.2). Namely, we apply the 
formula of Corollary \ref{Corollary H-6.5}(ii) for any column of a pretzel link. For $z=x^{-1}-x$ we get 
$$\Omega_L(x)= \nabla_L(z) = \sum_{j=0}^{m-1} s_{m,j}z^j\nabla_{T_{2,m-j}}(z) = $$  
$$ z^{m-1}({m-1 \choose 0} + s_{m,1}{m-2 \choose 0}+ s_{m,2}{m-3 \choose 0}+...) + $$ 
$$z^{m-3}({m-2 \choose 1} + s_{m,1}{m-3 \choose 1}+ s_{m,2}{m-4 \choose 1}+...) +...=$$ 
$$\sum_{j=0}^{\lfloor(m-1)/2\rfloor} \Bigg(\sum_{i=0}^{m-1-2j} {m-1-j-i \choose j}s_{m,j}\Bigg)z^{m-1-2j},$$
Where $s_{m,j}$ is an elementary symmetric polynomial in variables $k_1,..,k_m$ of degree $j$, 
that is $\Pi_{i=1}^m (z+k_i) = \sum_{j=0}^m s_{m,j}z^{m-j}$ and $\nabla_{T_{2,m-j}}(z)=
\nabla_{P_{1,1,...,1}}(z)$ is the 
Alexander-Conway polynomials of the torus links of type $(2,m-j)$, in particular, it satisfies 
Chebyshev type\footnote{We have  $\nabla_{T_{2,n}}(z) = i^{1-n}S_{n-1}(iz)$.} 
(compare Example \ref{Example IV.1.15}) relations 
$\nabla_{T_{2,n}}(z) = z\nabla_{T_{2,n-1}}(z)+ \nabla_{T_{2,n-2}}(z)$ 
(with initial data $\nabla_{T_{2,0}}(z)=0$ and $\nabla_{T_{2,1}}(z)=1$). In particular, 
 $\Omega_{T_{2,n}}(x) = \nabla_{T_{2,n}}(x^{-1}-x) = \frac{x^{-n}-(-1)^nx^n}{x^{-1}+x}= 
{n-1 \choose 0}z^{n-1}	+ {n-2 \choose 1}z^{n-3} +...+ {n-1-i \choose i}z^{n-1-2i} +... =
\sum_{i=0}^{\lfloor(n-1)/2\rfloor} {n-1-2i \choose i}z^{n-1-2i}$.
\end{example}

\begin{center}
\begin{tabular}{c}
\includegraphics[trim=0mm 0mm 0mm 0mm, width=.35\linewidth]
{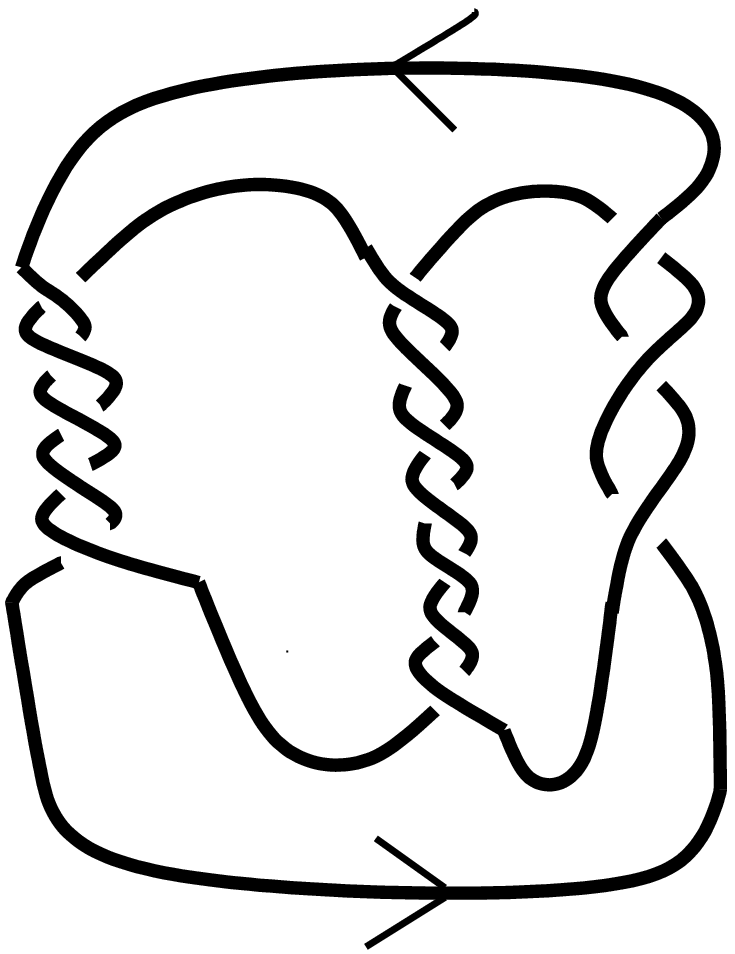}\\
\end{tabular}
\\
Fig. 6.2; $P_{5,7,-3}$ -- the pretzel knot with the trivial Alexander-Conway polynomial
\end{center}

\subsection{Tristram- Levine signature}
We generalize definition of the classical (Trotter-Murasugi) 
signature after Tristram and Levine (see \cite{Gor,Lev,P-T-2,Tr}).

Recall that a symmetric Hermitian form $h: C^n \times C^n \to C$ 
is a map which satisfies $h(a+b,c)=h(a,c)+h(b,c)$, 
$h(\lambda a,b) = \lambda h(a,b)$, and $h(a,b) = \overline{h(b,a)}$.
The matrix $H$ of a symmetric Hermitian form in any basis is called   
a Hermitian matrix (i.e. $H = \bar H^T$). A symmetric Hermitian form 
has a basis in which the matrix is diagonal with $1$, $-1$ or $0$ 
entries. The numbers, $n_1$ of $1$'s, $n_{-1}$ of $-1$'s and $n_0$ of $0$'s 
form a complete invariant of a symmetric Hermitian form (the Sylvester 
law of inertia).
The number $n_0$ is called the nullity of the form and 
$\sigma = n_1 - n_{-1}$ is called the signature of the form. Recall also 
that if we count eigenvalues of $H$ (with multiplicities) then 
$n_1$ is the number of positive eigenvalues of $H$ and $n_{-1}$ is 
the number of negative eigenvalues.

\begin{definition}(\cite{Tr,Lev})\label{Definition IV.5.7}
Let $A_L$ be a
Seifert matrix of a link $L$. For each complex number $\xi$ ($\xi\neq 1$)
consider the Hermitian matrix $H_L(\xi)=(1-\bar\xi)A_L+(1-\xi )A_{L}^{T}$. 
The signature of this matrix
 is called the Tristram-Levine signature of the link $L$. If the parameter $\xi$ is considered, 
we denote the signature by  $\sigma_L(\xi )$, if we consider $\psi = 1-\xi$ as a parameter, 
we use notation $\sigma_{\psi}(L)$
The classical signature $\sigma$ satisfies $\sigma(L)=\sigma_1(L)=\sigma_L(0)=\sigma_L(-1)$.
Also, by well justified convention, we put $\sigma_L(1)=0$ (see Remark \ref{Remark IV.5.8}).
\end{definition}

Tristram-Levine signature is a well defined link invariant
 as it is an invariant of 
$S$-equivalence of Seifert matrices.  Checking this is similar to the calculation 
for the potential 
function (we leave a pleasure exercise of verifying it to the reader).

\begin{remark}\label{Remark IV.5.8}
The signature of a Hermitian matrix is unchanged when matrix is multiplied 
by a positive number\footnote{The Hermitian matrix $H$ is Hermitian similar 
to $\lambda H$ for any real positive number $\lambda$; $\lambda H = (\sqrt \lambda Id)H (\sqrt \lambda Id)$.}, 
we can (and will) often assume that $\xi$ in $\sigma_L(\xi)$ and $\psi$ in $\sigma_{\psi}(L) $
are of unit length. With such assumptions we have Tristram-Levin 
signature functions, $\sigma_L(\xi), \sigma_{\psi}(L): S^1 \to Z$. 
$\sigma_L(\xi)$ is the signature function tabulated in \cite{Ch-L}, and $\sigma_{\psi}(L)$  is 
used in Examples in this book. $S^1$ will be usually 
parameterized by $arg(\psi)\in [-\pi,\pi]$.\footnote{In \cite{Ch-L}, $S^1$ is parameterized by
$\frac{arg \xi}{\pi}$.} Generally, we have $\sigma_L(\xi)= \sigma_{1-\xi}(L)$ 
but when restricted to the unit  circle, we have to write  $\sigma_L(\xi)= \sigma_{(1-\xi)/(|1-\xi|}(L)$.
Notice that for $\psi = \frac{1-\xi}{|1-\xi|}$, we have 
$\psi^2= \frac{(1-\xi)(1-\xi)}{(1-\xi)(1-\bar\xi)}=\frac{1-\xi}{1-\bar\xi}=-\xi$ (and $(i\psi)^2= \xi$).
Therefore,  $\sigma_{\psi}(L) = \sigma_L((i\psi)^2)= \sigma_L(\xi)$, for $Re (\psi)\geq 0$. 
As we show in Corollary \ref{Corollary IV.5.13}, 
$\sigma_{i}(L) = 0$, which justifies the convention\footnote{In the literature 
on Tristram-Levine signature of knots, often used normalization of the Hermitian matrix 
$(1-\bar\xi)A_L+(1-\xi )A_{L}^{T}$ is to take $|\xi|=1$ ($\xi\neq 1)$. When 
one writes the function $\sigma_K(\xi)$ then usual assumption about the parameter $\xi$ is that 
it is on the unit circle. Then one has 
$\det ((1-\bar\xi)A_L+(1-\xi )A_{L}^{T})= 
\det {\Large (}(\xi -1)(\frac{1-\bar\xi}{\xi -1}A -A^T){\Large )}= 
\det {\Huge  (}(\xi -1)(\bar \xi A -A^T){\huge )} \stackrel{.}{=}(\xi -1)^n\Delta(\bar\xi)$, where 
$ \stackrel{.}{=}$ denotes equality up to $\pm t^i$, \cite{Gor} (compare Lemma \ref{IV.4.33}). 
When dealing with links, we found more convenient (see \cite{P-T-2,P-1} to consider $\psi= 1-\xi$ 
and assume that $|\psi| =1$. Then we have 
$det (i(\bar \psi A + \psi A^T)) = det(i\bar \psi A - i\psi A^T) = 
\Omega(i\bar\psi) = \Omega(i\psi)=\nabla (-i(\bar\psi +\psi))$ (compare Lemma \ref{IV.4.33}). 
Therefore, for any knot $\sigma_{\psi}(K) = \sigma_K((i\psi)^2)= \sigma_K((\xi)$.} 
that $\sigma_L(1)=0$.
\end{remark}

\begin{corollary}[\cite{P-2}]\label{Corollary H-6.9}\ \\
\begin{enumerate}
\item[(i)]
For any $\bar t_{2k}$-move and $Re(1-\xi) \geq 0$ (i.e. $|arg(\psi)|\leq \pi/2$) we have:
$$ 0 \leq \sigma_{\bar t_{2k}(L)}(\xi) - \sigma_{L}(\xi) \leq 2.$$
In particular (\cite{P-T-2}, for 
$Re(1-\xi) \geq 0$, we have 
$-2\leq \sigma_{L_+}(\xi) -\sigma_{L_-}(\xi)\leq 0$.
\item[(ii)] Furthermore, for any $\xi$ and $k$ we have:
$$ 0 \leq |\sigma_{L_{\psfig{figure=infty-antiparal.eps,height=0.4cm}}}(\xi) - 
\sigma_{\bar t_{2k}(L)}(\xi)| \leq 1. $$ 
In particular, $0\leq |\sigma_{L_+}(\xi) - \sigma_{L_0}(\xi)|  \leq 1$.
\end{enumerate}
\end{corollary}

\begin{proof} Applying Lemma \ref{Lemma H-6.4} we obtain

  $H_{\bar t_{2k}(L}(\xi)=\left[
  \begin{array}{cc}
  H_{L_{\psfig{figure=infty-antiparal.eps,height=0.4cm}
   }}(\xi) & a \\
   a^{-T} & m+k(2-\xi-\bar\xi)
   \end{array}
   \right],$

    $H_{L}(\xi)=\left[
  \begin{array}{cc}
  H_{L_{\psfig{figure=infty-antiparal.eps,height=0.4cm}
   }}(\xi) & a\\
   a^{-T} & m
   \end{array}
   \right],$

where $a=(1-\bar\xi)\alpha+(1-\xi)\beta^T$ and
$m=((1-\bar\xi)+(1-\xi))q$. Because $2-\xi-\bar\xi\geq 0$,
so $0\leq\sigma (H_{\bar t_{2k}(L)}(\xi))-\sigma(H_{L}(\xi))\leq
   2$, and the proof of (i) is finished\footnote{It holds, in general, that 
if two $n \times n$ Hermitian matrices $H$ and $H'$ differ only at 
one entry, $a'_{nn} > a_{nn}$ then $0\leq \sigma(H') - \sigma(H) \leq 2$. 
Furthermore, if $\det H \det H' >0$ then $\sigma(H') = \sigma(H)$ and 
if $\det H \det H' <0 $ then $\sigma(H') = \sigma(H)+2$.}. Part (ii) 
follows from the easy observation that deleting the last row and column 
of a Hermitian matrix can change the signature at most by $\pm 1$.
\end{proof}

We can use results of computations in Examples \ref{Example H-5.11}, \ref{Example H-5.12}, and \ref{c3:4.11} to find 
the Tristram-Levine signature for the trefoil knot, the figure eight knot, 
and the pretzel knot $P_{2k_1+1,2k_2+1,2k_3+1}$.

\begin{example}\label{Example 4.28}
 Using the Seifert matrix for the right-handed trefoil 
knot ($\bar 3_1$) computed in Example \ref{Example H-5.11} we find that:
$$H_{\bar 3_1}(\xi) = 
\left[\begin{array}{cc}\xi + \bar\xi -2&1-\xi\\
1-\bar\xi&\xi + \bar\xi -2\\ 
\end{array}\right]$$
Therefore 
$$\sigma_{\bar 3_1}(\xi)=
\left\{
\begin{array}{ccc} 
-2 & if & Re(1- \xi) > \frac{1}{2}     \\
-1 & if & Re(1- \xi) = \frac{1}{2}     \\
 0 & if &  - \frac{1}{2} < Re(1- \xi) < \frac{1}{2}  \\
 1 & if &  Re(1- \xi) = -\frac{1}{2} \\
 2 & if &  Re(1- \xi) < \frac{1}{2}  
\end{array}\right. $$
\end{example}
Part of the regularity of the Tristram-Levine signature can be explained by 
the observation that for $\xi_2 = 2-\xi_1$ (i.e. $1-\xi_2 = -(1-\xi_1)$) 
we have $H_L(\xi_2)= -H_L(\xi_1)$ and $\sigma_L(\xi_2)= -\sigma_L(\xi_1)$.
 
\begin{example} Using the Seifert matrix for the figure eight knot 
($4_1$) computed in Example \ref{Example H-5.12} we find that:
$$H_{4_1}(\xi)= 
\left[\begin{array}{cc} 2- \xi - \bar\xi &\bar\xi -1\\ \xi - 1& 
       \xi + \bar\xi -2 \\ 
\end{array}\right]$$
For any $\xi \neq 1$, we have 
$\det H_{4_1}(\xi)= -(2- \xi - \bar\xi )^2 - (1-\xi)(1-\bar\xi) < 0$, 
thus $\sigma_{4_1}(\xi) = 0$.
\end{example}
The  observation that for the figure eight knot the Tristram-Levin 
signature is always equal to zero is not that unexpected because 
the figure eight knot is an ampchiheiral knot ($4_1 =\bar 4_1$) and we have:
\begin{corollary}\label{Corollary 4.30} 
If $\bar L$ is the mirror image 
of a link $L$ then the Seifert matrix $A_{\bar L} = -A_{L}$, 
$H_{\bar L}(\xi) = -H_{L}(\xi)$, $\sigma_{\bar L}(\xi) = -\sigma_L(\xi)$, and 
$\sigma_{\psi}(\bar L)= -\sigma_{\psi}(L)$.
 In particular, the Tristram-Levine 
signature of an  ampchiheiral link is equal to zero. 
\end{corollary}
We can also observe that $i$ ($i=\sqrt{-1}$) 
 times the matrix of $\tau$ from Exercise 5.15
is a Hermitian matrix of the signature equal to $0$ thus for a knot, 
$\sigma_i(K)=0$. This holds also for links as signature is unchanged by 
adding to the matrix rows and columns of zeros:
\begin{corollary}\label{Corollary IV.5.13}
For any link $L$ we have $\sigma_i(L)=\sigma_{-i}(L)= 0$.
\end{corollary}

It is useful to summarize our observations about the Tristram-Levin 
 signature of links using $\psi = 1-\xi$ and $|\psi|=1$.
\begin{corollary}\label{Corollary IV.5.14}
 When we change $\psi$ from $1$ to $i$, the signature $\sigma_{\psi}(L)$
changes from the classical (Trotter-Murasugi) $\sigma(L)$ to $0$ 
(equivalently, if $\xi$ changes from $1$ to $-1$, then $\sigma_L(\xi)$ 
changes from 0 to $\sigma(L)$). 
Furthermore, $\sigma_{\psi}(L)= \sigma_{\bar \psi}(L)= -\sigma_{-\psi}(L)= -\sigma_{\psi}(\bar L)$.
\end{corollary}

\begin{example} Using the Seifert matrix of the pretzel knot $P_{2k_1+1,2k_2+1,2k_3+1}$
computed in Example \ref{c3:4.11} we find that:
$$ H_{P_{2k_1+1,2k_2+1,2k_3+1}}=
\left[\begin{array}{cc}
-(\psi+\bar\psi)(k_1+k_2+1)    &    k_2\bar\psi + (k_2+1)\psi\\
 (k_2+1)\bar\psi + k_2\psi   &    -(\psi+\bar\psi)(k_2+k_3 +1)\\ 
\end{array}\right]
$$
Furthermore, 
$$\det H_{P_{2k_1+1,2k_2+1,2k_3+1}}= (\psi + \bar\psi)^2(1+k_1+k_2+k_3 + k_1k_2+k_1k_3+k_2k_3) -1.$$
Therefore the Tristram-Levine signature of a pretzel knot with $1+k_1+k_2+k_3 + k_1k_2+k_1k_3+k_2k_3 >0$ 
(e.g. a  positive pretzel knot) 
satisfies (in lieu of Corollary 6.12 we consider only $Re(\psi)\geq 0$):
$$\sigma_{\psi}(P_{2k_1+1,2k_2+1,2k_3+1}))=
\left\{
\begin{array}{ccc}
-2 & if & Re(\psi) > \frac{1}{2\sqrt {1+k_1+k_2+k_3 + k_1k_2+k_1k_3+k_2k_3}} \\
-1 & if & Re(\psi) = \frac{1}{2\sqrt {1+k_1+k_2+k_3 + k_1k_2+k_1k_3+k_2k_3}} \\
 0 & if & 0 \leq Re(\psi) < \frac{1}{2\sqrt {1+k_1+k_2+k_3 + k_1k_2+k_1k_3+k_2k_3}}
\end{array}\right. $$
Notice, that in the example of Seifert of $P_{5,7,-3}$, Figure 6.2,  
we have $ \det H_{P_{5,7,-3}}=-1$ and 
$\sigma_{\psi}(P_{5,7,-3})\equiv 0$. We utilize the result of this calculation in \cite{P-Ta}.
\end{example}

\subsection{Potential function and Tristram-Levine signature} 
Lemma \ref{c3:4.21} and Definition \ref{Definition IV.5.7} suggest that there is a relation between 
the potential function and the Tristram-Levine signature of links. 
In fact we have:

\begin{lemma}\label{IV.4.33}
Assume that the potential function at  $i\psi$ is 
different from zero. Then
$$i^{\sigma_{\psi}(L)} = \frac {\Omega_L(i\psi)}{|\Omega_L(i\psi)|}=
\frac {\Delta_L(t_0)}{|\Delta_L(t_0)|}= 
\frac {\nabla_L(-i(\psi +\bar\psi)}{|\nabla_L(-i(\psi +\bar\psi|},$$
where $\Delta_L(t_0)$ is the Alexander-Conway polynomial and 
$t_0=-\psi^2$ ($\sqrt {t_0} = -i\psi$).\ 
In particular, the Tristram-Levine signature is determined modulo $4$ by 
the appropriate value of the potential function 
(or Alexander-Conway polynomial); compare Chapter III of \cite{P-Book}.
\end{lemma}

\begin{proof} The idea is to compare the formulas for the potential functions and 
the signature, that is:
$$ \Omega_L(i\psi) = \det (i\psi A_L - (i\psi)^{-1}A^T_L) \ =$$
$$i^n \det (\psi A_L + \bar\psi A^T_L)\ \ \  and$$
$$\sigma_{\psi}(L)= \sigma( \bar\psi A_L + \psi A_L^T)$$
In more detail, we write our proof as follows:\\ 
Let $H$ be a non-singular Hermitian matrix of dimension $n$ and 
$\lambda_1,\lambda_2,...,\lambda_n$ its eigenvalues (with multiplicities). 
Then $det(iH)= i^n \det H = i^n \lambda_1 \lambda_2 \cdots \lambda_n = 
i^n(-1)^{n_{-1}}|\det H| = i^{n-2n_{-1}}|det H|= i^{n_{1}-n_{-1}}|det H|= 
i^{\sigma(H)}|det H|$.
Therefore, $\frac{det(iH)}{|det(iH)|} = i^{\sigma(H)}$. By 
applying this formula for 
$H= \psi A_L + \bar\psi A_L^T$, $|\psi|=1$, and remembering that 
$\sigma(\bar H) = \sigma (H)$, we obtain the formula of Lemma \ref{IV.4.33}.

\end{proof}
\begin{example}\label{Example IV. 5.17}
We can use Lemma \ref{IV.4.33} to compute quickly Tristram-Levine signature\footnote{It is, essentially, 
the same proof we used in Chapter III of \cite{P-Book} to show that a signature is a skein equivalence invariant: 
The Alexander-Conway polynomial determines signature modulo $4$ and the 
Murasugi type inequalities ($|\sigma_{\psi}(L_+)-\sigma_{\psi}(L_0)|\leq 1$ 
and for $Re(\psi)\geq 0$, $0\leq (\sigma_{\psi}(L_-)- \sigma_{\psi}(L_+)\leq 2$)
 gives the direction, and limit the size of the signature change, compare also Corollary \ref{Corollary H-6.5}.} of 
the torus link of type $(2,n)$, $T_{2,n}$. We use the fact that we already computed 
the classical signature, and Alexander-Conway (and potential) polynomial to be (for $k\neq 0$):
$$\sigma(T_{2,n})=1-n,\ \ \Delta_{T_{2,n}}(z)=\Omega_{T_{2,n}}(x)= 
\frac{x^{-n}-(-1)^nx^n}{x^{-1}+x}=  t^{\frac{1-n}{2}}\frac{t^n+(-1)^{n+1}}{t+1},$$ 
where $z=x^{-1}-x = t^{1/2}- t^{-1/2}$. In particular $\sigma_{\psi}(T_{2,n})$ can change 
only if $x=i\psi$ is a root of the potential function, and because 
$\Omega_{T_{2,n}}(i\psi) = i^{1-n}\frac{\psi^n - \psi^{-n}}{\psi - \psi{-1}}$, the only 
changes holds at $\psi$ satisfying $\psi^{2n}=1$ and $\psi\neq \pm 1$. \\
 We have for $Re \psi \geq 0$, $k\neq 0$, $0\leq j \leq n-1$: 
$$\sigma_{\psi}(T_{2,n})=
\left\{
\begin{array}{ccc}
1-n & if &  Re(\psi) > Re (e^{\pi/n}) \\
1-n +2j & if & Re (e^{j\pi/n}) > Re(\psi) >  Re (e^{(j+1)\pi/n}), \ j>0\\
2-n +2j & if & Re(\psi) = Re (e^{j\pi/n}),\ j>0. 

\end{array}\right. $$

\end{example}

\begin{corollary}\label{Corollary H-6.18} 
The classical (Trotter-Murasugi) signature $\sigma(L)=\sigma_1(L) = \sigma_L(-1)$, satisfies:
$$i^{\sigma(L)} = i^{\sigma (A_L + A_L^T)} = \frac{\Omega_L(i)}{|\Omega_L(i)} = 
\frac{\Delta_L(-1)}{|\Delta_L(-1)|} = \frac{Det_L}{|Det_L|}= \frac{\nabla (-2i)}{|\nabla (-2i)|}, \ \ assuming \ \ 
Det_L \neq 0;$$ 
here $\Delta_L(-1)$ denotes $\Delta_L(t)$ for $\sqrt t = -i$. Recall, that $Det_L = 
\Delta_L(-1)= \Omega_L(i) = det(i(A_L + A_L^T)$ is called 
the determinant\footnote{We should mention here that $|Det_L|$ is  
equal to $|det(G_L)|$ where $G_L$ is a Goeritz matrix of $L$. Furthermore, 
 if $D_L$ is a special diagram of an oriented link $L$ 
then $G_{D_L} = A_L + A_L^T$ for a properly chosen basis of 
$H_1(S)$ where $S$ is the Seifert surface of $D_L$ constructed according to 
Seifert algorithm. Thus not only $Det_L = det(iG_{D_L}))$ but also 
$\sigma(L) = \sigma (G_L)$; compare Corollary \ref{Corollary H-2.7}.} of a link $L$.

\end{corollary}
\begin{example}\label{Example IV.5.19} We compute here the Tristram-Levine signature 
of the knot $K=6_2$ using Lemma \ref{IV.4.33} and discuss the standard convention and notation.\\
We have:
$$\sigma_{\psi}(6_2)=
\left\{
\begin{array}{ccc}
-2 & if & Re(\psi) > \frac{1}{2}\sqrt{\frac{1+\sqrt 5}{2}} \\
-1 & if & Re(\psi) = \frac{1}{2}\sqrt{\frac{1+\sqrt 5}{2}} \approx 0.636... \\
 0 & if & 0 \leq Re(\psi) < \frac{1}{2}\sqrt{\frac{1+\sqrt 5}{2}}.
\end{array}\right. $$
\ \\
Step 1. We compute the the Conway polynomial $\nabla_{6_2}(z)= 1-z^2+z^4$; we use resolution 
in Figure 6.3 to find this value and also observe that changing a crossing at $p$ results in 
the trivial knot and smoothing at $p$ results in a connected sum of the right handed trefoil knot 
and the left handed Hopf link ($K^p_0= \bar 3_1 \# H_-$). In particular the unknotting number 
$u(6_2)=1$.\\
Step 2. $Det_K= \nabla_K(-2i)=-11$, thus $\delta (K)\equiv 2 \mod 4$, and because $K$ can be unknotted 
by changing one positive crossing, thus $-2 \leq \sigma (K) \leq 0$, and finally $\sigma (K)=-2$.\\
Step 2. Roots of $\nabla_{6_2}(z)$ are at $z^2=\frac{-1\pm \sqrt 5}{2}$. Thus for $t_0=\xi=-\psi^2$, 
$z=-i(\psi + \bar\psi)$, 
we have $\xi +\bar \xi = (i\psi)^2 +\bar{ (i\psi)^2}= z^2+2=\frac{3\pm \sqrt 5}{2}$. Because 
$|\psi|=|\xi|=1$, therefore $-2\leq \xi +\bar \xi \leq 2$ and $\xi +\bar \xi = \frac{3 - \sqrt 5}{2}$ 
($Re (\xi)= \frac{3 - \sqrt 5}{4}$). Finally, assuming $Re (\psi)\geq 0$ we get 
$\psi = \frac{1}{2}\sqrt{\frac{1+\sqrt 5}{2}})$.\\
Step 3. For $Re (\psi)\geq 0$, the value $Re (\psi) = \frac{1}{2}\sqrt{\frac{1+\sqrt 5}{2}}$ is the only 
place where the Tristram-Levine signature $\sigma_{\psi}(6_2)$ can be changing, 
and because we know already that $\sigma_1(6_2)=-2$ and $\sigma_i(6_2)=0$  we conclude that 
$\sigma_{\psi}(6_2)=-2$ if $Re(\psi) > \frac{1}{2}\sqrt{\frac{1+\sqrt 5}{2}}$ and 
$\sigma_{\psi}(6_2)=0$  if $0 \leq Re(\psi) < \frac{1}{2}\sqrt{\frac{1+\sqrt 5}{2}}$.\\
Step 4. It remains to show that $\sigma_{\psi}(6_2)=-1$ for $Re(\psi) = \frac{1}{2}\sqrt{\frac{1+\sqrt 5}{2}}$.
Here we argue that, because the considered $\psi$ is the singular root of the Alexander polynomial 
(precisely $t_0=-\psi^2$), therefore the value of signature at this point cannot differ by more than one 
from the neighboring values (so from $0$ and from $-2$). More detailed analysis of the Hermitian matrix 
$\bar\psi A + \psi A^T$, leads to the conclusion that if $t_0=-\psi^2$ is a singular root 
of the Alexander polynomial of a knot $K$ then $\sigma_{\psi}(K) = \frac{\sigma_{\psi_-}(K) + \sigma_{\psi_+}(K)}{2}$, 
where $\psi_-$ and $\psi_+$ are parameters just before $\psi$ and just after $\psi$ on 
the unit circle \cite{Mat}. 

In the convention of \cite{Gor,Ch-L}
one defines the Tristram-Levine signature function of variable $\xi$ ($|\xi|=1$) as
$\sigma_L(\xi)=\sigma((1-\bar\xi) A + (1-\xi) A^T))$. For $Re (\psi) \geq 0$, one has
$\sigma_{\psi}(L) = \sigma_L(\xi)$, where $\xi=-\psi^2$ ($\psi =\frac{1-\xi}{|1-\xi|}$).
In {\it knotinfo} Web page \cite{Ch-L}, the parameter $s$
satisfying $\xi=e^{\pi i s}$ is used. In particular, $\sigma_{6_2}(\xi)=-1$
for $Re (\xi) = \frac{3-\sqrt 5}{4}=1-\cos (\pi/5) \approx 0.191 $, and 
$s= arc cos (1-\cos (\pi/5))/\pi \approx 0.44$ (compare Remark \ref{Remark IV.5.8}).
   
\end{example}
\ \\
\centerline{{\psfig{figure=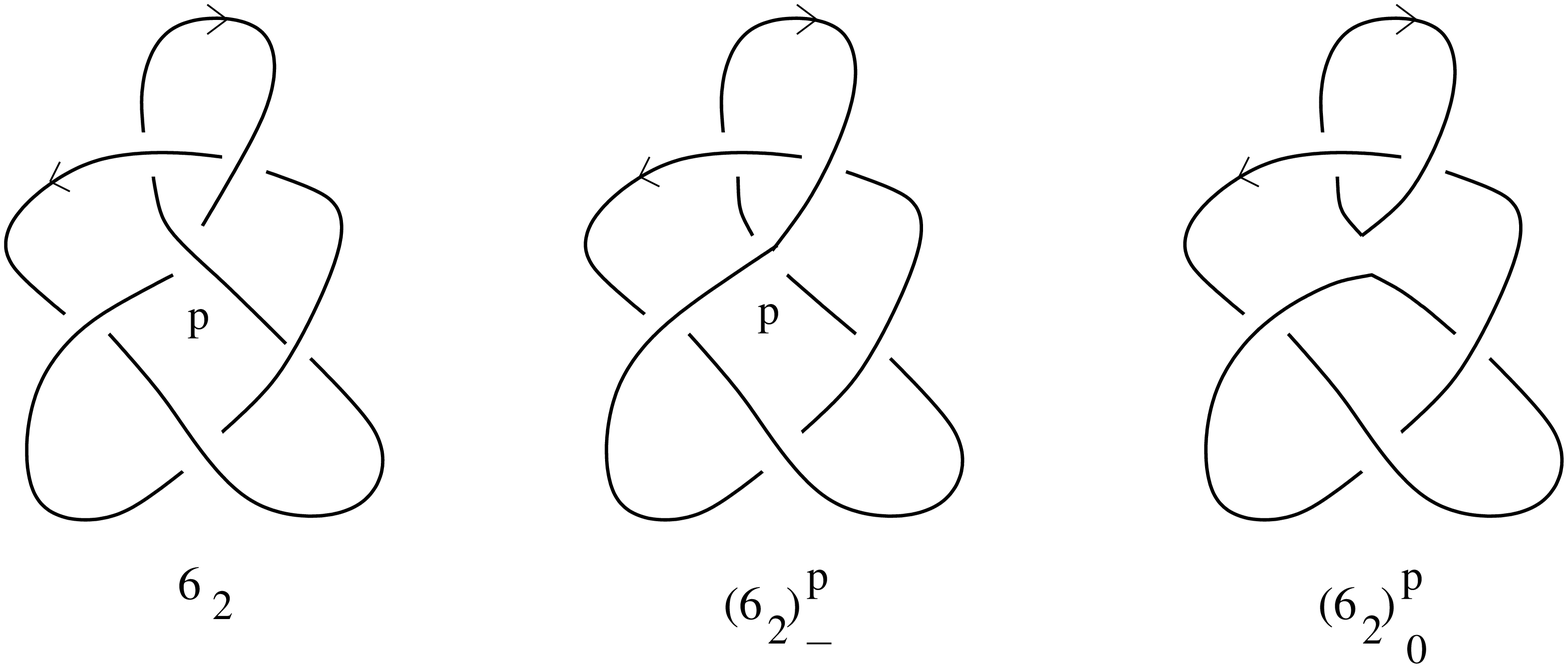,height=4.1cm}}}
\begin{center}
Fig.~6.3; Computing the Conway polynomial of the knot $K=6_2$.\\ 
$\nabla_{6_2}(z) = \nabla_{T_1}(z) +\nabla_{\bar 3_1 \# H_-}(z)= 1+(1+z^2)(-z)=1-z^2-z^4$
\end{center}

\begin{example}\label{Example IV.4.35}
The knot $9_{42}=$ 
{\parbox{1.9cm}{\psfig{figure=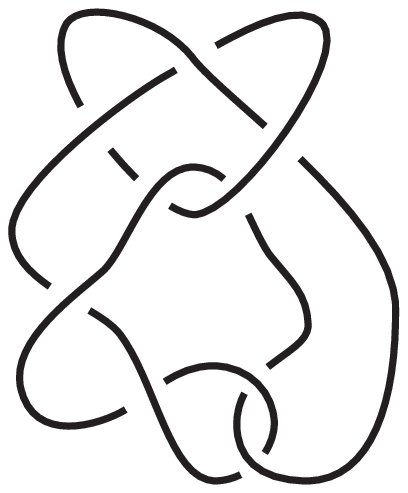,height=1.7cm}}}
is the smallest knot which is not amphicheiral but the Jones, Homflypt, and 
Kauffman polynomials are symmetric (e.g. $V_{9_{42}}(t) = V_{9_{42}}(t^{-1})$).
The non-amphicheirality of $9_{42}$ is detected by signature: 
$\sigma(9_{42}) = -2 = -\sigma(\overline{9_{42}})$. This description can 
leave however an impression that the fact that $9_{42}$ is not ambient 
isotopic to its mirror image cannot be checked by the Jones polynomial 
alone. However, it follows from Corollary 5.17 that 
$(-1)^{\sigma (K)/2}= sign (V_K(-1))$ for any knot $K$, thus if a knot is 
amphicheiral then $V_K(-1) >0$. For $9_{42}$ we have $V_{9_{42}}(-1)=
Det_{9_{42}} = -7 <0$ thus $9_{42}$ is not amphicheiral. Furthermore, 
because $9_{42}$ can be unknotted by changing one positive crossing, 
we can deduce that $\sigma(9_{42})=-2$.
\end{example}
In fact, the absolute value of the determinant $|Det_K|=|V_K(-1)|=|\Delta_K(-1)|= |\nabla (-2i)|$
suffices to show that  the knot $9_{42}$ is not amphicheiral. K.~Murasugi proved 
in \cite{M-10} (Theorem 5.6), the following result:
\begin{theorem}\label{Theorem IV. 5.19} For any knot $K$
$$\sigma_{K} \equiv |Det_K| - 1 \mod 4 $$

\end{theorem}
\begin{proof} We use the fact that $Det_K= \nabla (-2i) \equiv 1 \mod 4$.
Therefore, $|Det_K| \equiv 1 \mod 4$ if $Det_K > 0$ and $|Det_K| \equiv -Det_K \equiv -1 \mod 4$ if 
$Det_K < 0$. Furthermore, from Corollary \ref{Corollary H-6.18} follows that 
$Det_K = (-1)^{\sigma(K)/2}|Det_K|$. 
Therefore, 
$$\sigma_{K} \stackrel{\mod 4}{\ \ \equiv} \ 
\left\{ \begin{array}{ccc} 0 & if & |Det_K| \equiv 1 \mod 4 \\ 
2 & if & |Det_K| \equiv 3 \mod 4 
\end{array}\right. $$
and Theorem \ref{Theorem IV. 5.19} follows.
\end{proof}
Murasugi's Theorem leads to a curious formula:
\begin{corollary} For any knot $K$
$$Det_K = (-1)^{ (|Det_K|-1)/2} |Det_K|. $$
\end{corollary}

J.~Milnor proved that the signature of a knot with the Alexander polynomial equal to one is equal to 
zero \cite{Mil}. In fact, it follows directly from Lemma 5.16 that the Tristram-Levin signature 
can change only at roots of unit length of Alexander polynomial; therefore a link which has 
the Alexander polynomial 
without any root on the unit circle has constant Tristram-Levin signature function. Thus:
\begin{corollary} [\cite{Mil}] 
If the Alexander polynomial $\Delta_L(t)$ is different from 
zero on the unit circle then for any $\psi$, $(|\psi|=1)$, we have $\sigma_{\psi}(L)=0$.
\end{corollary}

If we assume only that the determinant of a knot is equal to 1 then we get as a 
conclusion that the signature is divisible by eight (compare \cite{M-9}, page 149 after 
Exercise 7.5.4):
\begin{proposition}\label{Proposition IV.5.22} 
If the determinant of a knot $K$ is equal to $1$ then $\sigma (K) \equiv 0 \mod 8$.
\end{proposition}
\begin{proof}
$Det_K =1$ means that for a Seifert matrix $A$ of a knot $K$, $\det (A +A^T) =1$; 
The matrix/form $A+A^T$ is often called the Trotter form. The diagonal entries of the 
Trotter form are even because the diagonal of $A+A^T$ is twice a diagonal of $A$. We can 
summarize these conditions by saying that  the Trotter form is  even and unimodular; 
recall that unimodularity means that $\det  (A +A^T) $ is invertible (here equal to $\pm 1$). 
The form is even if $x (A +A^T) x^T$ is always an even number.
Finally, every even unimodular form over $Z$ has its signature divisible by 8; 
see Theorem II.5.1 in \cite{M-H}. 
\end{proof}

\section{A combinatorial formula for the signature of alternating diagrams; \\
Quasi-alternating links}\label{Section 7. }
\markboth{\hfil{\sc Chapter IV. Goeritz and Seifert matrices }\hfil}
{\hfil{\sc IV.7.\  Combinatorial formula}\hfil}

Corollary \ref{Corollary H-6.18} has various interesting consequences. P.~Traczyk used it 
back in 1987 \cite{Tra}
 to find the combinatorial formula for the signature of alternating links,
starting from analysis of the condition $\sigma(L_+)=\sigma(L_0)-1$ (and $\sigma(L_-)=\sigma(L_0)+1$) 
and observing that it holds for any essential crossing of an alternating diagram.
The property was refined by Manolescu, Ozsvath, and Szabo 
and used to define quasi-alternating links \cite{O-S}, whose Khovanov \cite{Kho} and Heegaard 
Floer homology share with alternating links many interesting properties 
\cite{M-O,C-K} (compare Chapter X of \cite{P-Book}).
The property, of links which  Manolescu, Ozsvath, and Szabo observe to be important, and which 
always holds for alternating links, is the following (compare Subsection 1.4): 
$$|Det_{\parbox{0.6cm}{\psfig{figure=L+nmaly.eps}}}| = 
|Det_{\parbox{0.5cm}{\psfig{figure=L0nmaly.eps}}}|  +
|Det_{\parbox{0.5cm}{\psfig{figure=Linftynmaly.eps}}}| $$
The following result combines the above properties (compare \cite{M-O}):

\begin{theorem}\label{Theorem H-7.1}\ \\
The following two conditions are equivalent, providing that determinants 
of $L_0$ and $L_{\infty}$ are not equal to zero\footnote{In (a) one deals with a
 Kauffman skein triple of unoriented links; in (b) 
one chooses any orientation of $L_+$ 
(e.g. \parbox{0.7cm}{\psfig{figure=L+maly.eps}})
and related orientation of $L_0$ 
(\parbox{0.5cm}{\psfig{figure=L0maly.eps}}), and
any orientation of $L_{\infty}$ (e.g. 
{\parbox{0.5cm}{\psfig{figure=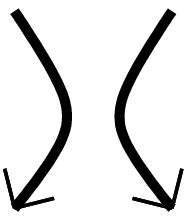,height=0.5cm}}} or 
{\parbox{0.5cm}{\psfig{figure=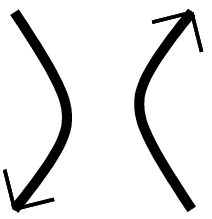,height=0.5cm}}}).}\\
(a)\ \  
$|Det_{L_+}| = |Det_{L_0}| + |Det_{L_{\infty}}|$  \\
(b) \ \ $\sigma(L_+) = \sigma(L_0) - 1$ and $\sigma(L_+) = 
\sigma(L_{\infty})- \frac{1}{2}(w(L_0) - w(L_{\infty}))$. \\
The similar equivalence also holds for a negative crossing:\\
($a'$) \ \ $|Det_{L_-}| = |Det_{L_0}| + |Det_{L_{\infty}}|$  \\
($b'$)  \ \ $\sigma(L_-) = \sigma(L_0) + 1$ and $\sigma(L_-) =
 \sigma(L_{\infty}) + \frac{1}{2}(w(L_0) - w(L_{\infty}))$. 
\end{theorem}
\begin{proof} ((a) $\Leftrightarrow$ (b)): \ We apply the formula $Det (L) = i^{\sigma(L)}|Det(L)|$ and 
use the relation between the Jones polynomial, and its Kauffman bracket variant, with 
the signature. Recall,  that the Jones polynomial 
$V_L(t)$ of an oriented link $L$ is normalized to be one for 
the trivial knot and satisfies the skein relation
$ t^{-1}V_{{\psfig{figure=L+maly.eps}}}(t) -
tV_{{\psfig{figure=L-maly.eps}}}( t) =
(t^{\frac{1}{2}} - t^{-\frac{1}{2}})
V_{{\psfig{figure=L0maly.eps}}}(
t)$. For $t=-1$ (or, more precisely, $\sqrt t = i$)  we obtain exactly 
the skein relation of the determinant: 
$ Det_{{\psfig{figure=L+maly.eps}}} -
Det_{{\psfig{figure=L-maly.eps}}} =
-2i
Det_{{\psfig{figure=L0maly.eps}}}$.
Thus $Det_L = V_L(-1)$; $\sqrt t = i$.
Recall also that the Kauffman bracket polynomial of unoriented link diagrams, 
$\langle D \rangle \in Z[A^{\pm 1}]$,  
is defined by the following properties \cite{K-6}:
\begin{enumerate}
\item[(i)] $\langle \bigcirc\rangle = 1$

\item[(ii)] $\langle\bigcirc\sqcup D\rangle = -(A^2+A^{-2})\langle D\rangle$

\item[(iii)] $\langle$ {\parbox{0.5cm}{\psfig{figure=L+nmaly.eps}}}
$\rangle =
A\langle $\parbox{0.5cm}{\psfig{figure=L0nmaly.eps}}
$\rangle
+
A^{-1}\langle$
{\parbox{0.5cm}{\psfig{figure=Linftynmaly.eps}}}
$\rangle$
\end{enumerate}
Furthermore,
if $\vec D$ is an oriented diagram with underlying unoriented diagram $D$ 
then $V_{\vec D}(t) = (-A^3)^{w(\vec D)}\langle D \rangle$.
Thus for $A^2=-i$ ($A^4=-1$) we get: \\
 $Det (D)= (-A^3)^{-w(D)}<D>= A^{w(D)}<D>$.  
Recursive formula for the Kauffman bracket \
$<D_+> = A<D_0> + A^{-1}<D_{\infty}>$ leads to \\
$(-A^3)^{w(D_+)}Det(D_+) = 
A(-A^3)^{ w(D_0)}Det(D_0) + A^{-1}(-A^3)^{ w(D_{\infty})}Det(D_{\infty})$\\
then leads to $A^{-w(D_+)}Det(D_+) = A^{1-w(D_0)}Det(D_0)+ 
A^{-1-w(D_{\infty})}Det(D_{\infty})$, \\ then leads to\\ 
$A^{-w(D_+)}i^{\sigma(D_+)}|Det(D_+)| = 
A^{1-w(D_0)}i^{\sigma(D_0)}|Det(D_0)| +
A^{-1-w(D_{\infty}}i^{\sigma(D_{\infty}}|Det(D_{\infty})| $\\ 
and eventually to \\ 
$ |Det(D_+)| = $
$$A^{w(D_+) - w(D_0) +1}i^{\sigma(D_0)-\sigma(D_+)}|Det(D_0)| + 
A^{w(D_+) - w(D_{\infty}) -1}i^{\sigma(D_{\infty}) - \sigma(D_+)}
|Det(D_{\infty})|=$$ 
$$i^{\sigma(D_0)-\sigma(D_+)-1}|Det(D_0)| + 
i^{\sigma(D_{\infty}) - \sigma(D_+) - 1/2(w(D_0)-w(D_{\infty}))}
|Det(D_{\infty})|.$$ 
When we compare this formula with that of Theorem \ref{Theorem H-7.1}(a) we see that (a) holds iff 
$i^{\sigma(D_0)-\sigma(D_+)-1}=1$ and 
$i^{\sigma(D_{\infty}) - \sigma(D_+) - 1/2(w(D_0)-w(D_{\infty}))} = 1$ 
and these conditions are equivalent to conditions \\
$\sigma(D_0)-\sigma(D_+) \equiv 1 \mod 4$ and 
$\sigma(D_{\infty}) - \sigma(D_+) - 
\frac{1}{2}(w(D_0)-w(D_{\infty})\equiv 0 \mod 4 $.\ These conditions are 
equivalent to (b) because by 
Corollary \ref{Corollary IV 1.14}(i), we have generally that $|\sigma (D_+)-\sigma (D_0)| \leq 1$).
Furthermore, in general, we have that $|\sigma(D_+) - \sigma(D_{\infty}) 
+ \frac{1}{2}(w(D_0)-w(D_{\infty})|\leq 2$. The last inequality require 
some explanation and consideration of two cases in which 
{\parbox{0.5cm}{\psfig{figure=L+maly.eps}}} is either a mixed crossing or a 
self-crossing.
\begin{enumerate}
\item[(m)] If {\parbox{0.5cm}{\psfig{figure=L+maly.eps}}} is a mixed 
crossing then let $D_j$ be a component of $D_+$ such that the change 
of the orientation of $D_j$ results in the link $D'_- = $ 
{\parbox{0.5cm}{\psfig{figure=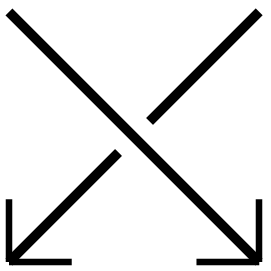,height=0.5cm}}}. \ 
Then by Corollary  \ref{Corollary IV 1.14} 
$|\sigma({\parbox{0.5cm}{\psfig{figure=BookIVLpr-.eps,height=0.5cm}}}) -
\sigma({\parbox{0.5cm}{\psfig{figure=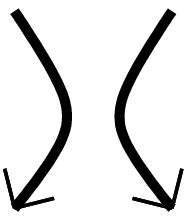,height=0.5cm}}})| 
\leq 1$. Further, by Proposition 2.11(ii),\\
$|\sigma({\parbox{0.5cm}{\psfig{figure=L+maly.eps}}} + 2lk(D_j,D_+-D_j) - 
\sigma({\parbox{0.5cm}{\psfig{figure=BookIVLinf.eps,height=0.5cm}}})|
\leq 1$. Because $4lk(D_j,D_+-D_j) = w(D_+) - w(D'_-) = 
w({\parbox{0.5cm}{\psfig{figure=L0maly.eps,height=0.5cm}}}) - 
w({\parbox{0.5cm}{\psfig{figure=BookIVLinf.eps,height=0.5cm}}})+2$ we 
obtain \\
$|\sigma(D_+) - \sigma(D_{\infty}) + \frac{1}{2}(w(D_0)-w(D_{\infty})) + 1| 
\leq 1$ and finally \\
$-2\leq \sigma(D_+) - \sigma(D_{\infty}) - 
\frac{1}{2}(w(D_0)-w(D_{\infty}))\leq 0$.

\item[(s)] If {\parbox{0.5cm}{\psfig{figure=L+maly.eps}}} is a 
self-crossing then in $D_0=$ {\parbox{0.5cm}{\psfig{figure=L0maly.eps,height=0.5cm}}} 
the two parallel arcs belong to different link components. 
Let $D_j$ component contain the lower arc and let 
$D'_0 = $ {\parbox{0.5cm}{\psfig{figure=antiparal.eps,height=0.5cm}}}
be obtained from $D_0$ by changing the orientation of $D_j$. 
After performing the second Reidemeister move on $D'_0$ we obtain 
a diagram 
{\parbox{1.2cm}{\psfig{figure=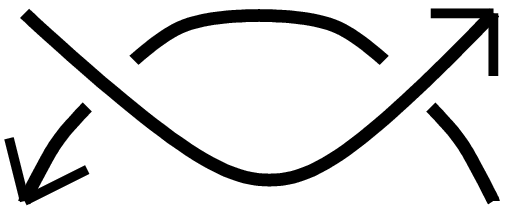,height=0.5cm}}}  
which has two mixed crossings.
We use Corollary  \ref{Corollary IV 1.14}(i) on one of them to get\
$|\sigma ({\parbox{0.5cm}{\psfig{figure=antiparal.eps,height=0.5cm}}}) 
- \sigma ({\parbox{0.5cm}{\psfig{figure=BookIVinf-anti.eps,height=0.5cm}}})|
\leq 1$. Because $\sigma (D'_0) = \sigma (D_0) + 2lk(D_j,D_0-D_j) = 
\sigma (D_0) - \frac{1}{2}(w(D_0)-w(D'_0))
=\sigma (D_0) - \frac{1}{2}(w(D_0)-w(D_{\infty}))$, we obtain\\
$|\sigma ({\parbox{0.5cm}{\psfig{figure=L0maly.eps,height=0.5cm}}})
- \sigma ({\parbox{0.5cm}{\psfig{figure=BookIVinf-anti.eps,height=0.5cm}}}) +
\frac{1}{2}(w({\parbox{0.5cm}{\psfig{figure=L0maly.eps,height=0.5cm}}}) - 
w({\parbox{0.5cm}{\psfig{figure=BookIVinf-anti.eps,height=0.5cm}}}))|\leq 1$,\ then\\
$|(\sigma ({\parbox{0.5cm}{\psfig{figure=L0maly.eps,height=0.5cm}}})
- \sigma ({\parbox{0.5cm}{\psfig{figure=L+maly.eps,height=0.5cm}}})) +
(\sigma ({\parbox{0.5cm}{\psfig{figure=L+maly.eps,height=0.5cm}}})
- \sigma ({\parbox{0.5cm}{\psfig{figure=BookIVinf-anti.eps,height=0.5cm}}})+
\frac{1}{2}(w({\parbox{0.5cm}{\psfig{figure=L0maly.eps,height=0.5cm}}}) -
w({\parbox{0.5cm}{\psfig{figure=BookIVinf-anti.eps,height=0.5cm}}})))|\leq 1$,\\
and finally 
$|\sigma ({\parbox{0.5cm}{\psfig{figure=L+maly.eps,height=0.5cm}}})
- \sigma ({\parbox{0.5cm}{\psfig{figure=BookIVinf-anti.eps,height=0.5cm}}})+
\frac{1}{2}(w({\parbox{0.5cm}{\psfig{figure=L0maly.eps,height=0.5cm}}}) -
w({\parbox{0.5cm}{\psfig{figure=BookIVinf-anti.eps,height=0.5cm}}}))|\leq 2$
as required.
\end{enumerate}
The equivalence $(a') \Leftrightarrow (b')$ follows from  (a) $\Leftrightarrow$ (b)  
by considering mirror 
images of diagrams from (a) and (b). In particular, 
for the diagram $\bar D$ being the mirror image of $D$, we always have that 
$\sigma (\bar D) = - \sigma (D)$, and  $w(\bar D) = -w(D)$).

\end{proof}

It is not difficult to see that any crossing of an alternating diagram 
satisfies properties (a),(a') of Theorem \ref{Theorem H-7.1}.
This follows from the fact that if $D$ is an alternating diagram then also $D_0$ and $D_{\infty}$ 
are alternating, 
and for an alternating diagram $|Det_D|$ can be
interpreted as the number of spanning trees of the underlying Tait graph, $G_b(D)$,
and the number of spanning trees is additive under deleting contracting
rule; see Subsection \ref{Subsection 1.4}. 
These ideas are developed in Chapter V of \cite{P-Book}. Without referring to
it, the properties (a) and (a') of alternating links follow from the
proof of Traczyk formula for the signature of alternating diagrams which
we present below.  First, we have to recall the necessary
terminology. In fact, we use this as an opportunity for introducing 
basic language which unifies the notion of Tait surface and Tait graph 
(Footnote 13) with that of Seifert surface and Seifert graph \cite{Crom}. 
Before general definition let us recall the definition of the Seifert graph.

\begin{definition}\cite{Crom}\label{Definition SG}\ \\
The {\sl Seifert graph} of an oriented  diagram $\vec D$ is a  signed (planar) graph
$\Gamma (\vec{D})$ whose vertices correspond to Seifert circles of the
diagram and edges correspond to crossings of the diagram. The sign
of an edge is determined by the sign of the corresponding
crossing.
\end{definition} 

In the more general setting we allow arbitrary smoothings of crossings of (not necessary 
oriented) diagram $D$. 
\begin{definition}\label{Definition H-7.3}
A Kauffman state $s$
of $D$ is a function from the set of crossings of $D$ to
the set $\{+1,-1\}$. Equivalently,
we assign to each crossing of $D$ a marker
according to the following convention:\\
\ \\
\centerline{\psfig{figure=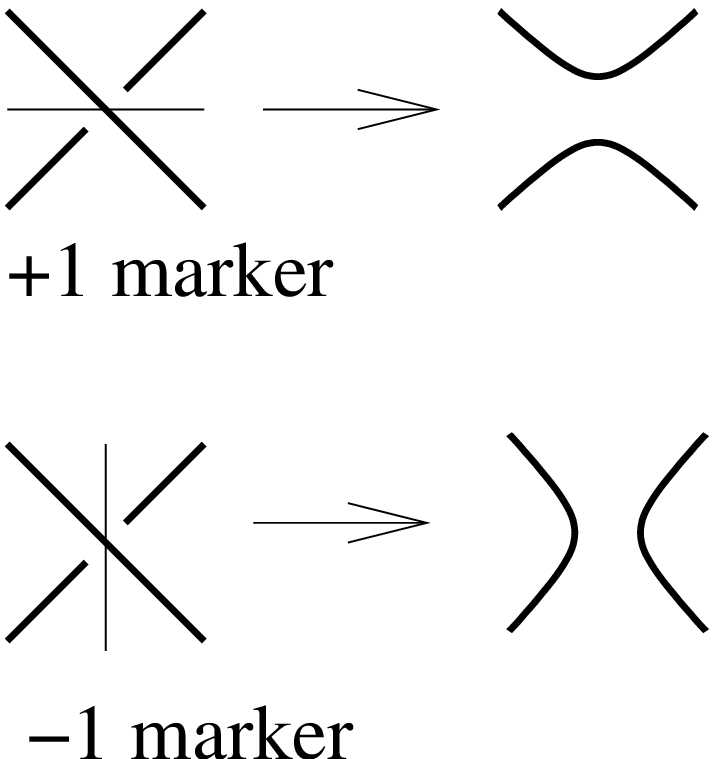,height=4.5cm}}
\begin{center}
Fig. 7.1; markers and associated smoothings
\end{center}

By $D_s$ we denote the system of circles in the diagram
obtained by smoothing all
crossings of $D$ according to the markers of the state $s$, Fig. 7.1. 
$|D_s|$ denotes the number of circles in $D_s$.

\end{definition}

In this notation the Kauffman bracket polynomial of $D$ is given by
the state sum formula:
$$<D> = \sum_s A^{\sigma(s)}(-A^2-A^{-2})^{|D_s|-1},$$
where $\sigma(s) = \sum_p s(p)$ is the number
of positive markers minus the number of negative markers in the state $s$.\\
The state sum formula looks like a useful but not necessarily sophisticated tool, 
however, state sums (and their limits) are basic and deep concepts in the
statistical physics and very likely the next breakthrough in Knot Theory (and more) 
will utilize a connection (still to be discovered) between phase transition of a physical 
system and Khovanov type homology based on closeness of states of the system 
(possibly persistent homology \cite{E-Ha} will play a role here).

But we are straying too far from our local goal of associating graphs and 
surfaces to any Kauffman state $s$.
\begin{definition}(\cite{PPS})\label{Definition PPS}\ \ 
\begin{enumerate}
\item[(i)]
Let $D$ be a diagram of a link and $s$ its Kauffman state. We form
a graph, $G_s(D)$, associated to $D$ and $s$ as follows.
 Vertices of $G_s(D)$ correspond to circles of $D_s$.
Edges of $G_s(D)$ are in bijection with crossings of $D$ and an
edge connects given vertices if the corresponding crossing
connects circles of $D_s$ corresponding to the
vertices. As in the case of the Tait graph, $G_s(D)$ 
is a signed graph where the sign of an edge $e(p)$ is $s(p)$, that 
is the sign of the marker of the Kauffman state $s$ at the crossing $p$.
\item[(ii)] 
In the language of associated graphs we can state the definition
of an s-adequate diagrams as follows:\ the diagram $D$ is s-adequate
 if the graph $G_s(D)$ has no loops (adequacy is studied and utilized in Chapter V of \cite{P-Book}).
\item[(iii)] We associate with every Kauffman state $s$ of a diagram $D$, a 
surface $F_s(D)$ embedded in $R^3$ and with $\partial F_s(D)=D$, 
in a manner similar to Construction \ref{c3:2.4} of a Seifert surface. That is, 
we start from the collection of circles $D_s$. Each of the circles bounds a disk 
in the projection plane. We make the disks disjoint by pushing them slightly up 
above the plane of projection, starting from the innermost disks.
We connect the disks together at the original crossings of the diagram $D$
by half-twisted bands so that the 2-manifold which we obtain
has $D$ as its boundary, see Figure 3.5 (we ignore orientation of the diagram, and 
the resulted surface can be unorientable). Equivalently, we can start a construction of 
$F_s(D)$ from the graph $G_s(D)$ as a spine (strong deformation retract) of the 
constructed surface and proceed as follows:
The graph $G_s(D)$ possesses an additional structure, that is 
 a cyclic ordering of edges at every vertex
following the ordering of crossings at any circle of $D_s$. The
sign of each edge is the label of the corresponding crossing. In
short, we can assume that $G_s(D)$ is a ribbon (or framed) graph,
 and that with every state we associate a surface $F_s(G)$
whose core is the graph $G_s(D)$. $F_s(G)$ is naturally embedded
in $R^3$ with $\partial F_s(G)= D$. If $s$ is the state separating black regions 
of checkerboard coloring of $R^2-D$ then $F_s(G)$ is the Tait surface of 
the diagram described in Exercise 2.8. For $s=\vec{s}$, that is, $D$ is
oriented and markers of $\vec{s}$ agree with orientation of $D$,
$G_s(D)$ is the Seifert graph of $D$ and $F_s(G)$ is the Seifert
surface of $D$ obtained by Seifert construction. We do not use
this additional data in this survey but it may be of great use in
analysis of Khovanov homology (compare \cite{A-P} or Chapter X of \cite{P-Book}).
\end{enumerate}

\end{definition}
The surface $F_s(G)$ is not the only surface associated with the graph 
$G_s(D)$, another such surface is Turaev surface, $M(s)$ \cite{Tu}, which 
for positive ($s_+$) or negative ($s_-$) states of an alternating diagram is a planar surface. 
With some justification Turaev surface can be called a background surface of a diagram.
The construction of $M(s)$ for a given state $s$ of $D$ is illustrated, after \cite{Tu},
in Figure 7.2. That is, $M(s)$ is obtained from a regular neighborhood of a projection of a link 
by modifying (by half-twists) neighborhoods of $s$-wrong edges (see Figure 7.2 and 
compare it to Figure 1.8 to see that any alternating digram has only $s_+$-true edges). 
Notice, that $M(s)$ depends on $s$ and the link projection but not over-under information 
of a link diagram.
Alternatively, we can say that $M(s)$ is a surface realizing the natural cobordism 
between circles of $D_s$ and circles of $D_{-s}$. In \cite{DFKLS} the Turaev genus of a link is 
defined to be the minimal genus of Turaev surface over all diagrams $D$ of a link with $s_+(D)$ states. 
The immediate consequence is 
that alternating link has the Turaev genus equal to zero. Notice also, that if we cup off the circles 
of $D_s$ in $M(s)$ by 2-discs we obtain the surface, $M^+(s)$ with boundary $D_{-s}$ and the graph 
$G_s(D)$ as its spine.
\ \\ \ \\
\centerline{\psfig{figure=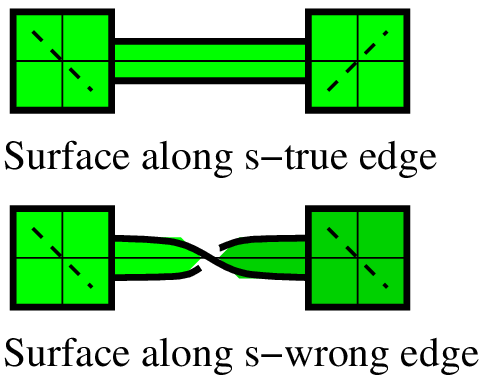,height=4.5cm}}
\begin{center}
Fig. 7.2; Turaev surface $M(s)$ is composed of squares along every crossing of $D$ connected 
by ribbons according to convention illustrated in this Figure. s-true edge and s-wrong edge 
are arcs of the diagram $D$ connecting crossings and the name depends on the label given 
by $s$ to boundary crossings \cite{Tu}  
\end{center}

Going back to Traczyk's combinatorial formula, we recall  the convention for
checkerboard shading of the projection plane. In an alternating
diagram we choose  the {\sl standard shading} as in Figure 7.3(a)
complementary to the shading given in Figure 7.3(b) (this essentially agrees with
Tait's convention of checkerboard coloring, however we do not assume that the outside region is 
white or black).
\ \\
\centerline{{\psfig{figure=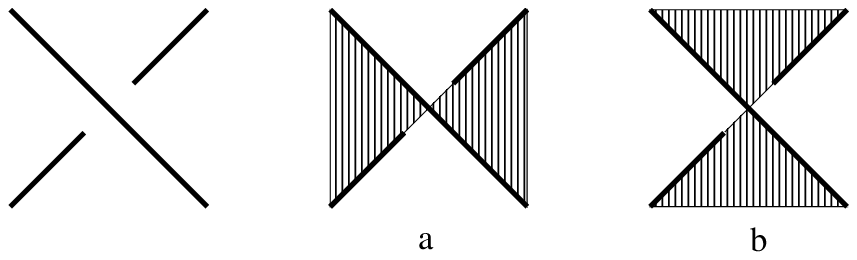,height=2.1cm}}}
\begin{center}
Fig.~7.3; Checkerboard shading of the plane of the projection: (a) Tait's,\ (b) dual to Tait's
\end{center}

\ \\
We denote by $B$ the number of
black (shaded) areas and by $W$ the number of white areas
(for an alternating diagram $D$ we have $B=|D_{s_-}|$ and $W=|D_{s_+}|$). Furthermore for 
an oriented diagram $\vec D$ let $\Gamma (\vec D)$ denote its Seifert graph (Definition \ref{Definition SG}),
$T$ its (signed) spanning tree and $d_+(T)$ (resp. $d_-(T)$) the number of positive (resp. negative) 
edges in $T$. For an alternating diagram the numbers $d_+(T),d_-(T)$ do not depend on $T$ so we can write 
$d_+(\vec D)$,and $d_-(\vec D)$ in this case\footnote{This is the case for more general class of 
homogeneous diagrams introduced in \cite{Crom} and defined as diagrams for which 2-connected 
components of the Seifert graph have all edges of the same sign (i.e. they are homogeneous). 
Alternating diagrams are special cases of homogeneous diagrams; this well known fact follows also 
from Lemma \ref{Lemma H-7.5} as the lemma can be proved for a fixed choice of a spanning tree and 
the left side of the equation does not depend on the choice of a spanning tree.}

\begin{lemma}\label{Lemma H-7.5}
If $\vec D$ is an oriented connected alternating diagram of a link then 
$$ \frac{1}{2}(w(\vec{D}) + |D_{s_+}| - |D_{s_-}|) =d_+(\vec{D}) - d_-(\vec{D})$$
In particular, the left hand side of the equation is unchanged when one goes from 
$\vec D$ to $\vec D^p_0$ for a non-nugatory crossing $p$ (in $\vec D^p_0$ the crossing $p$ is smoothed 
according to orientation of $\vec D$). 
\end{lemma}
\begin{proof} One can easily proof Lemma \ref{Lemma H-7.5} by induction on the number of 
non-nugatory crossings of $\vec D$. First one observes that if $\vec D$ has only nugatory crossings then 
$\Gamma (\vec D)$ is a tree and 
$d_+(\vec D)= c_+(\vec D) = s_+(\vec D)-1$ (and $d_-(\vec D)= c_-(\vec D) = s_-(\vec D)-1$), thus 
the formula in Lemma \ref{Lemma H-7.5} holds. In the inductive step we consider a non-nugatory crossing $p$ of 
$\vec D$ and compare ingredients of the formula for $\vec D$ and $\vec D^p_0$, and having the 
formula for  $\vec D^p_0$ deduct it for $\vec D$. It is worth however to compare $d_+$, $d_-$, 
$c_+$, $c_-$,  $|D_{s_+}|$, and  $|D_{s_-}|)$ in more detail.
\end{proof}

\begin{lemma} Let $p$ be any crossing of an oriented diagram $\vec D$. Then
\begin{enumerate}
\item[(i)] $$\vec s (p) = 
\left\{
\begin{array}{ccc}
s_+(p) & if & p\ is \ positive \\
s_-(p) & if & p \ is \ negative
\end{array}\right. $$
In particular if $\vec D$ is a positive diagram then $\vec s = s_+$, and if 
$\vec D$ is a negative diagram then $\vec s = s_-$.
\item[(ii)] $|(\vec D^p_0)_{\vec s}| = |\vec D_{\vec s}|,$
\item[(iii)] $$|(\vec D^p_0)_{s_+}| =  
\left\{
\begin{array}{ccc}
 |\vec D_{s_+}| & if & p\ is \ positive \\
 |\vec D_{s_+}|-\varepsilon_+  & if & p \ is \ negative
\end{array}\right. $$

 $$|(\vec D^p_0)_{s_-}| = 
\left\{
\begin{array}{ccc}
 |\vec D_{s_-}|-\varepsilon_-  & if & p\ is \ positive \\
  |\vec D_{s_-}|   & if & p \ is \ negative
\end{array}\right. $$ 
Here $\varepsilon_+$ and $\varepsilon_-$ are $+1$ or $-1$. If $p$ is 
a non-nugatory crossing of an alternating diagram then $\varepsilon_+ = \varepsilon_- =1$.
\end{enumerate}
\end{lemma}
\begin{proof} (i) The proof is illustrated in Figure 7.4.\\
The other parts are equally elementary and we leave them as exercise for the reader.
\end{proof}
\vspace*{0.2in} \centerline{\psfig{figure=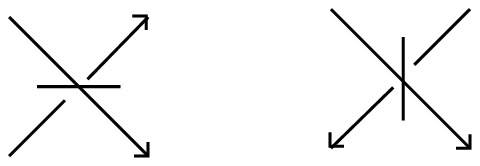,height=3.5cm}}
\begin{center}
Fig. 7.4; $\vec s(p) = s_+(p)$ if $sgn(p)=1$, \ \ and \ \
$\vec s(p) = s_-(p)$ if $sgn(p)=-1$
\end{center}

\begin{lemma}\label{Lemma H-7.7} 
 If $D$ is a connected alternating diagram,  then for 
a complex number $A$ such that $A^4=-1$, we have: 
\begin{enumerate}
\item[(i)] $<D>_{A^4=-1}= A^{B-W}|<D>_{A^4=-1}|$ 
\item[(ii)] For any crossing $p$ of an alternating diagram $D$ one has:
$$|<D>_{A^4=-1}| = |<D^p_0>_{A^4=-1}| + |<D^p_{\infty}>_{A^4=-1}|$$
in other words the absolute value of the determinant of a diagram is additive under 
the Kauffman bracket skein triple. 
\end{enumerate}
\end{lemma}
\ \\
\begin{proof} If all crossings of $D$ are
nugatory, then $D$ represents the trivial knot. Choose an orientation of $D$.
The orientation defines signs of crossings, which are however independent
on chosen orientation. As we noticed in Lemma XX 
in this case $c_+=W-1 $ and $c_-=B-1$. 
Thus $\langle D\rangle = (-A^3)^{w(D)}= (-A^3)^{W-B}$ (for a 
knot $w(D)$ does not depend on the orientation of $D$). For $A^4=-1$,
$\langle D\rangle_{A^4=-1}= (-A^4)^{W-B}(A)^{B-W}=A^{B-W}$ as required.
The inductive step follows easily:
If $p$ is a non-nugatory crossing of $D$, then from the Kauffman bracket skein relation
$$\langle D\rangle = A \langle D_{0} \rangle+ A^{-1}\langle D_{0+}\rangle$$
and by the inductive assumption, for $A^4=-1$, follows that:
$$\langle D\rangle_{A^4=-1} = AA^{B-W-1}|\langle D^p_0\rangle_{A^4=-1}| + 
A^{-1} A^{B-W+1}|\langle D^p_{\infty}\rangle_{A^4=-1}| =$$ 
$$A^{B-W}{\Large (}|\langle D^p_0\rangle_{A^4=-1}| + |\langle D^p_{\infty}\rangle_{A^4=-1}|{\Large )} =  
A^{B-W}|\langle D\rangle_{A^4=-1}|$$
which completes the proof of Lemma 7.7(i). It also establishes Lemma 7.7(ii) for 
a non-nugatory crossing $p$ of a connected diagram $D$. If $p$ is a nugatory crossing, 
then $\langle D^p_0\rangle_{A^4=-1}|$ or $|\langle D^p_{\infty}\rangle_{A^4=-1}|$ is equal to zero 
and (ii) holds immediately. If $D$ is not connected diagram then (ii) holds for any 
connected component of $D$ and (ii) follows because Kauffman bracket (and signature) is 
multiplicative under disjoint sum.

\end{proof}

As a corollary of Theorem \ref{Theorem H-7.1}, \ref{Lemma H-7.5}, and \ref{Lemma H-7.7}, 
we have Traczyk's result.
\begin{theorem}\cite{Tra}\label{Theorem Tr4.2}
If $D$ is a reduced\footnote{{\it Reduced} means that no crossing of
$D$ is nugatory and the crossing $p$ of $D$ is called nugatory if
$D^p_0$ has more (graph) component from $D$.}
alternating diagram of an oriented link,  then
\begin{enumerate}
\item [(1)]
$\sigma(D) = -(c_+ - c_-) + d_+ - d_- =
-w + d_+ - d_- $
\item [(2)]
$\sigma(D) = -{1\over 2}(c_+ - c_-) + {1\over 2}(W-B) =
-{1\over 2}w + {1\over 2}(W-B) =  \\
\ \ -{1\over 2}(w + |D_{s_+}| - |D_{s_-}|)$
\item [(3)] $\sigma(D)= \sigma(D^p_0)- sign(p)$
\end{enumerate}
\end{theorem}


\subsection{Quasi-alternating links}\ \\
Quasi-alternating links introduced by Manolescu, Ozsvath, and Szabo in 
\cite{O-S,M-O,C-K}
are motivated by 
properties (a),(a') of Theorem \ref{Theorem H-7.1}, described in the theorem relations to signature, 
and applications of these properties to the thinness of Khovanov and Heegaard Floer homology: 
\begin{definition}\cite{O-S}\label{Definition H-7.9}\ 
The family of quasi-alternating links is the smallest family of links 
which satisfies:
\begin{enumerate}
\item[(i)] The trivial knot is quasi-alternating.
\item[(ii)] If $L$ is a link which admits a crossing such that \\
(1) \ both smoothings ($L_0$ and $L_{\infty}$) are quasi-alternating, and\\
(2) \ $|Det_L| = |Det_{L_0}| + |Det_{L_{\infty}}|$,\\
 then $L$ is quasi-alternating.
\end{enumerate}
The crossing used in the definition is called a quasi-alternating crossing 
of $L$.
\end{definition}

Notice that a split link has its determinant equal to 0 so it cannot be quasi-alternating 
(determinants of quasi-alternating links are always positive as easily follows by induction 
from Definition \ref{Definition H-7.9}). 
Therefore, we can use condition (b) of Theorem \ref{Theorem H-7.1}
as alternative definition of the family of quasi-alternating links.

One can ask why we choose condition (2) in the definition of quasi-alternating links
and not a weaker first part of conditions (b), (b') from the Theorem \ref{Theorem H-7.1} 
($\sigma(D_+) = \sigma(D_0)-1$ or $\sigma(D_-) = \sigma(D_0)+1$). The first answer 
is purely practical: this is exactly what is needed to have thin Khovanov (and Heegaard) homology 
 (see Chapter X of \cite{P-Book}). One can also argue that condition which refers only to unoriented 
links is sometimes a plus.

We already have proved that non-split alternating links satisfy properties which make them 
quasi-alternating: if $D$ is an alternating diagram then also $D_0$ and $D_{\infty}$ are
alternating, and every non-nugatory crossing of an alternating diagram is
quasi-alternating (satisfies property (ii)(2)) as long as $D$ is a non-split link.

According to \cite{M-O} among the 85 prime knots with up to nine crossings,
$82$ are quasi-alternating ($71$ are alternating), 
$2$ are not quasi-alternating ($8_{19}$ and $9_{42}$),
 and the knot $9_{46}$ still remains undecided. It was showed by 
 A.~Schumakovitch using odd Khovanov homology that $9_{46}$ is not quasi-alternating. 
The classification of quasi-alternating knots up to 11 crossings was 
completed by J.~Greene in \cite{Gr}.
  
It was also determined which pretzel links are quasi-alternating (partial 
classification of quasi-alternating Montesinos links is advanced in \cite{C-K,Gr,J-S,Wid}:
\begin{theorem}\cite{C-K,Gr} [Characterization of quasi-alternating pretzel links]\ \\
The pretzel link $P_{(1,...,1,p_1,...,p_n,-q_1,...,q_m)}$ with $e$ $1$th, $e+n+m\geq 3$,
 and $p_i\geq 2$, $q_i\geq 3$ is 
quasi-alternating if and only if one of the conditions below holds:\\
(1) $e \geq m$, \\
(2) $e=m-1 > 0$, \\
(3) $e=0$, $n=1$, and $p_1> min(q_1,...,q_m)$, \\
(4) $e=0$, $m=1$, and $q_1> min(p_1,...,p_n)$, \\
The same is true  on permuting parameters\footnote{Thus all pretzel links are covered in the theorem.} $p_i$ and $q_j$.
\end{theorem}
The importance of quasi-alternating links rests in the following results of 
Manolescu and Ozsvath:\\
(1) Quasi-alternating links are Khovanov homologically $\sigma$-thin 
(over Z).\\
(2) Quasi-alternating links are Floer homologically $\sigma$-thin 
(over $Z_2$).\\ 
We explain the meaning of the first result in Chapter X of \cite{P-Book} showing also 
how to generalize it to Khovanov homologically $k$-almost thin links. 

To have some measure of complexity or depth of quasi-alternating links we 
introduce the
quasi-alternating computational tree index $QACTI(L)$
 is defined inductively from the definition of quasi-alternating link as follows:\\
\begin{definition}\label{Definition H-7.11}
For the trivial knot $T_1$, $QACTI(T_1)=0$. $QACTI(L)$ is the minimum over all quasi-alternating
crossings $p$  (of any diagram) of  $L$ of $max(QACTI(L^p_0),QACTI(L^p_{\infty}))+1$.

In other words, $QACTI(L)$ is the minimal depth of any binary computational resolving tree of $L$
using only quasi-alternating crossings and having the trivial knot in leaves.
\end{definition}

From Definitions \ref{Definition H-7.9} and \ref{Definition H-7.11},
 and Theorem \ref{Theorem H-7.1} we get approximation on $QACTI(L)$:
\begin{corollary}\label{Corollary 7.12} Let $L$ be a quasi-alternating link then:
\begin{enumerate}
\item[(i)] $|Det(L)|-1 \geq QACTI(L) \geq log_2(|Det(L)|)$
\item[(ii)] $QACTI(L) \geq |\sigma (\vec L)|$, for every orientation of $L$.
\item[(iii)] If $p$ is a quasi-alternating crossing of $L$ then\\
 $QACTI(L)\leq QACTI(L^p_0) +1$, and
$QACTI(L)\leq QACTI(L^p_{\infty})+1$.
\end{enumerate}
\end{corollary}

Let us finish this survey with a nice example of a quasi-alternating 
 knot of 13 crossings due to S.~Jablan and R.~Sazdanovic \cite{J-S}.

\ \\
\centerline{\psfig{figure=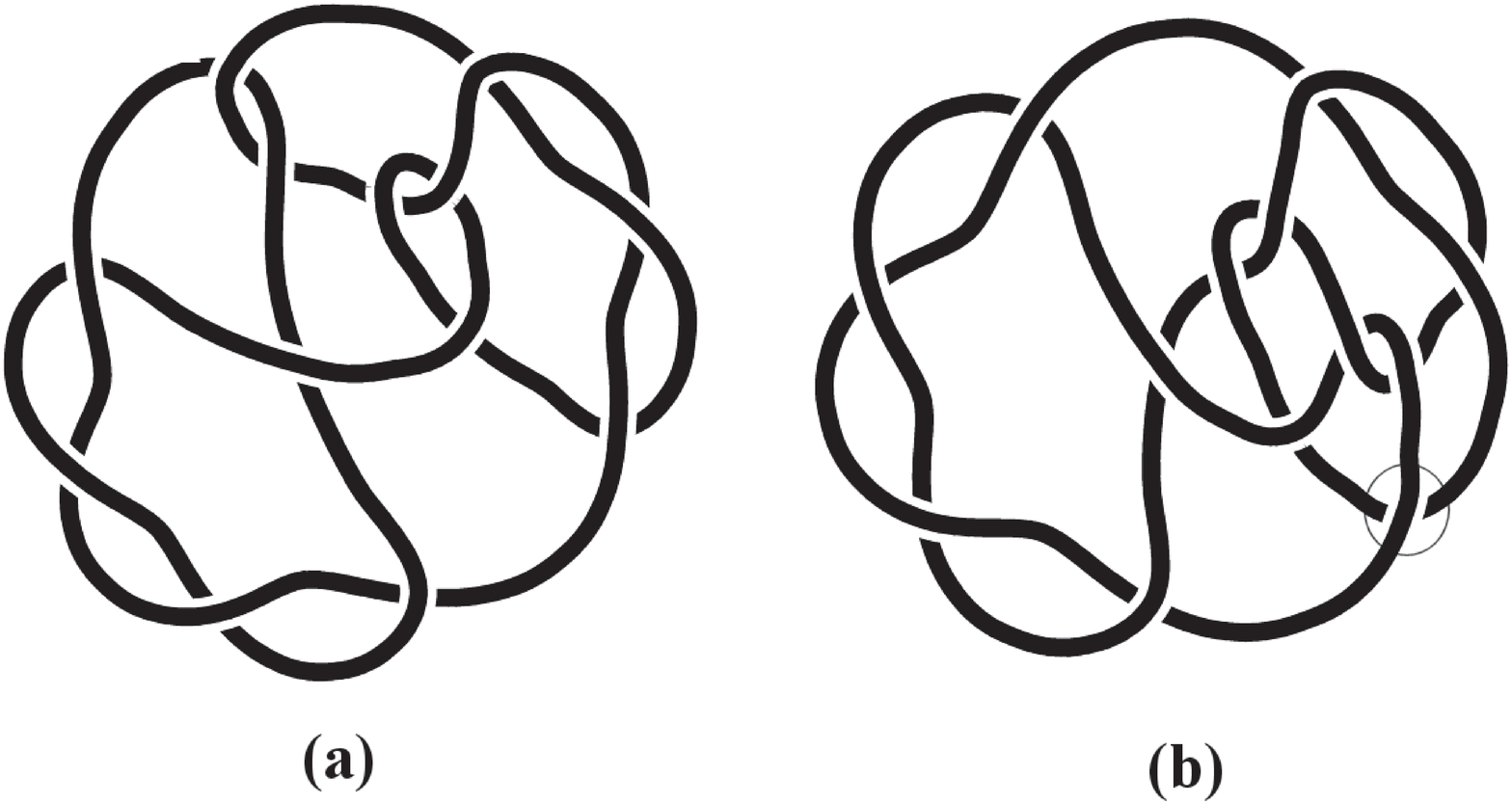,height=8.5cm}}
\begin{center}
  Figure 7.5; A quasi-alternating knot $13_{n_{1659}}$ with 2 diagrams of 
(minimal number) 13 crossings. The first diagram is (Conway) algebraic but 
no crossing is quasi-alternating. The second diagram, which bases on Conway's 
polyhedron $6^*$, has the circled crossing quasi-alternating. The determinant 
of $13_{n_{1659}}$ is equal to $51$ while smoothings of quasi-alternating 
crossing gives the trivial knot and a quasi-alternating link with determinant $50$,
\cite{J-S}.
\end{center}

\ \\ \ \\ \ \\
\noindent \textsc{Dept. of Mathematics,
The George Washington University, Washington, DC 20052};\
e-mail: {\tt przytyck@gwu.edu}

\end{document}